\documentclass[11pt, letterpaper,reqno]{amsart}
\usepackage[left=1in,right=1in,bottom=0.5in,top=0.8in]{geometry}

\usepackage[toc,page]{appendix}

\usepackage{tikz}
\usepackage{amsmath, amsthm, amsfonts, amssymb, bm,bbm}
\usepackage{mathrsfs}

\usepackage{relsize}
\usepackage{booktabs}
\usepackage[numbers]{natbib}
\usepackage{setspace}
\usepackage{graphicx}
\usepackage{enumitem}
\usepackage[labelfont=bf, font={small,it}]{caption}
\usepackage{float}
\usepackage{mathtools}
\usepackage{array}
\usepackage{xcolor}
\usepackage{subcaption}
\usepackage{fancyhdr}
\usepackage{mdframed}
\usepackage{listings}
\usepackage{upgreek}
\usetikzlibrary{patterns} 
\usetikzlibrary{arrows}
\usetikzlibrary{decorations.markings} 
\usetikzlibrary{shapes}
\usepackage{microtype}
\usepackage[hidelinks]{hyperref}

\newtheorem{thm}{Theorem}[]
\newtheorem{cor}[thm]{Corollary}

\newtheorem{lem}[thm]{Lemma}

\newtheorem{prop}[]{Proposition}
\newtheorem{rmk}[]{Remark}
\newtheorem{defi}[thm]{Definition}

\newcommand{\Pb}{\mathbb{P}}
\newcommand{\E}{\mathbbm{E}}
\newcommand{\Id}{\mathbbm{1}}

\newcommand{\I}{{\rm i}}
\newcommand{\s}{\tilde{s}}
\newcommand{\rr}{\tilde{r}}
\newcommand{\A}{\widetilde{A}}
\newcommand{\B}{\widetilde{B}}
\newcommand{\G}{\mathcal{G}}

\newcommand{\R}{\mathbb{R}}
\newcommand{\KK}{\mathcal{K}}
\newcommand{\wKK}{\widehat{\mathcal{K}}}
\newcommand{\tKK}{\widetilde{\mathcal{K}}}

\newcommand{\Z}{\mathbb{Z}}
\renewcommand{\Re}{\mathrm{Re}\,}

\DeclareMathOperator{\Pf}{\mathrm{Pf}}

\newcommand{\zcd}{\raisebox{-0.5\height}{\includegraphics[scale=0.04]{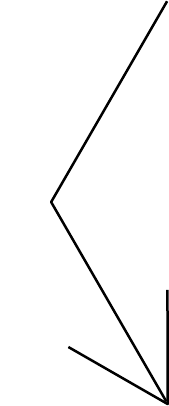}}}
\newcommand{\wcu}{\raisebox{-0.5\height}{\includegraphics[scale=0.04]{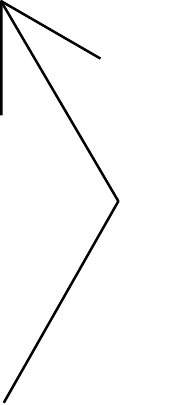}}}

\title{Stationary Half-space geometric Last Passage Percolation without Random Shift}
\author{Alessandra Occelli, Jiyue Zeng}
\date{January 2026}

\begin{document}
\raggedbottom

% The AMS classes require the abstract before \maketitle.
\begin{abstract}
    We derive exact formulas for the diagonal distribution of half-space geometric last-passage percolation with a two-parameter family of stationary initial conditions throughout the phase diagram. Unlike previous results, our formulas remove the random shift at the origin. Under critical scaling, we obtain new formulas for the half-space stationary KPZ fixed-point distribution at $x=0$ in the high-density and maximal-current phases. In the low-density phase, we derive a new representation and show that it matches the corresponding formula from the stationary log-gamma polymer.
\end{abstract}

\maketitle

\setcounter{tocdepth}{1}
\tableofcontents

\section{Introduction}

Models of last passage percolation (LPP) form a central class of exactly solvable systems in the Kardar--Parisi--Zhang (KPZ) universality class \cite{KPZ}. Through their correspondence with interacting particle systems and directed polymers, they provide a powerful framework for the rigorous analysis of KPZ fluctuations. In full-space geometry, a detailed understanding has been achieved for a wide range of initial and boundary conditions, including curved, flat, and stationary cases \cite{Johansson2000,PrSp2002,BaikRains2000}.

In recent years, considerable progress has been made in understanding KPZ fluctuations in half-space geometries under various initial conditions \cite{BBCS,BaikRainsHalfSpace,BeteaFerrariOccelli2020, BEFO22,barraquandKPZBrownian, barraquandKPZ, zeng, zeng2026}, where the presence of a boundary constraining admissible paths gives rise to new phenomena absent in full space.
While integrable half-space LPP models with specific boundary weights exhibit Pfaffian point process structures, the presence of a boundary fundamentally alters both the fluctuation behavior and the algebraic structure of the model. In particular, the stationary case --- corresponding to stationary initial data in the KPZ fixed point --- has remained largely inaccessible to exact analysis.

A natural stationary half-space LPP model can be defined by choosing boundary weights that ensure stationarity of increments \cite{BeteaFerrariOccelli2020,BC23}, in analogy with the full-space construction via Burke-type properties \cite{Balazs2006}. However, this stationary half-space model is not directly integrable: correlations introduced at the boundary destroy the Pfaffian structure that underlies solvability. This difficulty already appears at finite size and persists in scaling limits. %To recover a formula suitable to asymptotic analysis, a shift argument is introduced 

As in the full-space setting, the stationary half-space last passage percolation can be accessed via a shift argument, whereby the stationary distribution is obtained from a one-parameter integrable model by translating the last passage time and performing an analytic continuation in the boundary parameters \cite{BeteaFerrariOccelli2020}.

In a recent work, \cite{zeng} generalized this approach to a family of two-parameter stationary half-space geometric LPP and \cite{zeng2026} solved the one-point distribution of the stationary log-gamma polymer model. In the physics work \cite{barraquandKPZ}, the authors solved the one-point distribution of the half-space KPZ equation with Hariya-Yor stationary initial condition. While these results constitute an important step forward, all of them contain a residual shift that prevents the derivation of a closed formula for the stationary diagonal distribution itself and obstructs a direct asymptotic analysis at the KPZ scale. The main contribution of the present paper is to overcome this obstruction.
\smallskip

\noindent\textbf{\textcolor{blue}{Half space LPP.}} In this paper we consider half-space last passage percolation on the domain
$$\mathbb{H}=\{(i,j)\in\Z^2: 1\leq j\leq i\},$$
where admissible paths move by up-steps (0,1) and right-steps (1,0) and are constrained to remain in $\mathbb{H}$. To each site $(i,j)\in\mathbb{H}$ we assign a random weight, and for any endpoint $(m,n)\in\mathbb{H}$ we define the last passage time $G(m,n)$ as the maximum total weight collected along all admissible paths from $(1,1)$ to $(m,n)$.

We focus on geometric half-space LPP models with boundary inhomogeneities. In the bulk $(i>j)$, weights
are taken to be independent and identically distributed, while along the diagonal and the first row the weights are modified to encode boundary interactions and initial conditions.

\medskip

\noindent\textbf{\textcolor{blue}{Fluctuation regimes for half-space LPP.}}
Depending on the boundary strength and the distance from the diagonal, the half-space LPP exhibits qualitatively different fluctuation behaviors:
\begin{itemize}
\item a boundary-dominated Gaussian regime, where maximizing paths stick macroscopically to the diagonal and fluctuations
are of order $N^{1/2}$;
\item a bulk-dominated KPZ regime, where the boundary is asymptotically negligible and fluctuations follow Tracy--Widom laws on the scale;
\item critical crossover regimes, where boundary and bulk compete and lead to GOE--GSE-type transitions;
\item and a distinguished stationary regime, corresponding to the diagonal characteristic direction, where boundary and bulk contributions balance and KPZ fluctuations persist.
\end{itemize}
The present paper is concerned exclusively with this last regime.

\medskip

\noindent\textbf{\textcolor{blue}{Half-space stationary model.}}
In full space, stationary LPP and stationary KPZ are by now well understood and lead to the Baik--Rains distribution as the universal one-point law \cite{BaikRains2000,BFP10}. 

In half space, the situation is more intricate: for a fixed boundary interaction strength, there exists a one-parameter family of stationary measures indexed by the asymptotic slope of the initial condition. 
These measures can be constructed and proved to be invariant under the dynamics of half-space growth models in \cite{BC23}. 
We first define this two-parameter family of stationary measures.
\begin{defi}\label{def: stationary measure}
    For any fixed $\sqrt{q} \in (0,1),$ $s\in (\sqrt{q},1]$, $r \in (0,1/s],$ we define
    \begin{equation}
        \mathcal{I}_{r,s}(x) := R_2(x) + \max\bigg\{0,\,\, \max_{1\leq k \leq x, k\in \Z} \{R_1(k)- R_2(k-1) - Y\}\bigg\}, \quad \text{ for } x \in \Z_{\geq 1},
    \end{equation}
    where $R_1$ and $R_2$ are independent $\text{Geom}(\sqrt{q}/s)$ and $\text{Geom}(\sqrt{q}s)$ random walks with $R_1(0) = 0$ and $R_2(0) = 0$, and $Y \sim$ $\text{Geom}(rs)$ is independent. We let $\mathcal{I}_{r,s}(0) = 0$. If $rs=1$, we define $Y = \infty$.
\end{defi}
When $r = s$ or $ r= 1/s$, the process $\mathcal{I}_{r,r}(\cdot)$ or $\mathcal{I}_{r,1/r}(\cdot)$ becomes a Geom$(\sqrt{q}/r)$ random walk, which we call the product stationary measure. We call $\mathcal{I}_{r,s}$ the two-parameter stationary measure.

\begin{prop}[\textnormal{Proposition 3.2 \cite{BC23}}]\label{prop: stationary measure}
    For any fixed $\sqrt{q} \in (0,1),$ $s\in (\sqrt{q},1],$ $r \in (0,1/s],$ the process $(\mathcal{I}_{r,s}(x))_{x\in\Z_{\geq 0}}$ is stationary for the half-space geometric LPP model $\eqref{recurrence}$.
\end{prop}

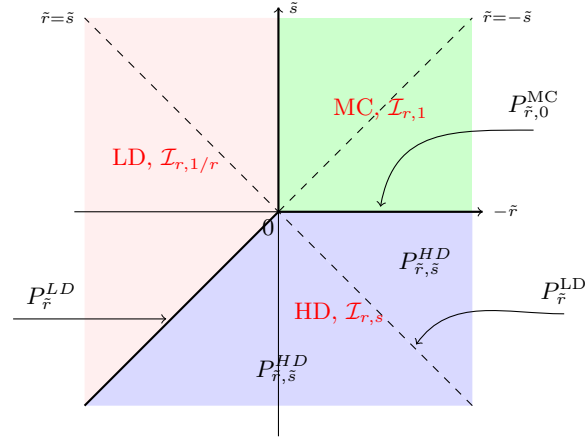
\begin{figure}[t]
\centering
    \begin{tikzpicture}[scale=1.35]
        \fill[green!20] (0,0) -- (1.9,0) -- (1.9,1.9) -- (0,1.9) -- (0,0) -- cycle;
        \fill[red!20,opacity=0.3] (0,0) -- (0,1.9) -- (-1.9,1.9) -- (-1.9,-1.9) -- (0,0) -- cycle;
        \fill[blue!15] (0,0) -- (1.9,0) -- (1.9,-1.9) -- (0,0) -- cycle;
        \fill[blue!15] (0,0) -- (1.9,-1.9) -- (-1.9,-1.9) -- (0,0) -- cycle;
        
        \draw[thin,->] (-2,0) -- (2,0)node[right] {\footnotesize{$\mathsmaller{-\rr}$}};
        \draw[thin,->] (0,-2.2) -- (0,2)node[right] 
        {\footnotesize{$\mathsmaller{\s}$}};

        \draw[thick] (0,0) -- (0,2);
        \draw[thick] (0,0) -- (2,0);
        
        \draw[thick] (-1.9,-1.9) -- (0,0);
        \draw[thin, dashed] (0,0) -- (1.9,1.9) node[right] {\footnotesize{$\mathsmaller{\rr = -\s}$}};
        \draw[thin,dashed] (1.9,-1.9) -- (-1.9,1.9) node[left] {\footnotesize{$\mathsmaller{\rr = \s}$}}; 
        \node at (-0.1,-0.15) {\footnotesize{$0$}};
        
        \node[text = {red}] at (0.6,-1){\footnotesize{HD, $\mathcal{I}_{r,s}$}};
        \node[text = {red}] at (1,1){\footnotesize{MC, $\mathcal{I}_{r,1}$}};
        \node[text = {red}] at (-1.1,0.5){\footnotesize{LD, $\mathcal{I}_{r,1/r}$}};

        \draw[thin,->]
        (2.5,0.8) node[above]{\footnotesize $P^{\mathrm{MC}}_{\rr,0}$}to[out=180,in=80,looseness=1.2] (1,0.05);
    
      % From left, straight into the LD label:
      \draw[thin,->]
        (-2.6,-1.05) -- (-1.1,-1.05)
        node[near start,above]{\footnotesize $P^{LD}_{\rr}$};
    
      % From bottom, curved arrow to the HD region boundary:
      \draw[thin,->]
        (2.8,-1) node[above]{\footnotesize $P^{\mathrm{LD}}_{\rr}$}to[out=180,in=55] (1.35,-1.3);

    \node[text = {black}] at (1.4,-0.5){\footnotesize{ $P_{\rr,\s}^{HD}$}};
    \node[text = {black}] at (0,-1.55){\footnotesize{ $P_{\rr,\s}^{HD}$}};
        
    \end{tikzpicture}
    \captionof{figure}{This picture describes the ergodic theory of the geometric LPP model given boundary parameter and drift parameter. Depending on where $(-\log(r), \log(s))$ lies in the diagram, the process $G(\cdot,N) - G(N,N)$ converges to one of the three spatial processes, $\mathcal{I}_{r,s},$ $\mathcal{I}_{r,1/r},$ or $\mathcal{I}_{r,1},$ which is claimed in Proposition $\ref{conjecture}$. On the full line $r=s$ and the half line $rs=1, r\geq 1, s\leq 1$, the process converges to the \text{Geom}$(\sqrt{q}/r)$ random walk. \textbf{MC}, \textbf{HD}, and \textbf{LD} represent maximal-current (green), high-density (blue), and low-density (red) phases. Consider the critical scaling of parameters $r = 1 + \rr(c_0N/2)^{-1/3}$ and $s = 1 + \s(c_0N/2)^{-1/3}$. {We characterize the asymptotic distribution of $G(N,N)$ in the high-density phase, denoted by $P_{\rr,\s}^{HD}$ and in the maximal-current regime, denoted by $P_{\rr,0}^{MC}$}. We also obtain, by a separate product-stationary construction, the low-density/product-stationary formula $P_{\rr}^{LD}$. Our formulas remove the random shift, thereby yielding geometric LPP that starts precisely with the corresponding stationary measures.}
    \label{phaseDiagram}
\end{figure}

The convergence of LPP to the stationary measure has been proved in a very recent paper \cite{das2026invariant}. We record it in the following proposition.
\begin{prop}[\textnormal{Theorem 1.10 and Theorem 1.12 \cite{das2026invariant}}]\label{conjecture}
    For given $\sqrt{q} \in (0,1),$ $r \in (0,1/\sqrt{q}),$ $\{\mathcal{I}_{r,s}\}_{s \in (\sqrt{q}, \min\{1,1/r\}]}$ constitutes all extremal stationary measures for $\eqref{recurrence}$. For any initial condition $G(N,1)$ such that $\lim_{N\rightarrow\infty} G(N,1)/N = \rho \in R_{>0}$, we have the following. Let \begin{equation}\label{relation}
        s = \frac{(1+\rho)\sqrt{q}}{\rho}.
    \end{equation}
    \begin{itemize}
        \item(high-density) For $\log(s) \leq -\log(r)$ and $\log(s) \leq 0,$ the process
\[
    \left(G(x+N, N) - G(N,N)\right)_{x\in \Z_{\geq 0}}
\]
 converges weakly as $N \rightarrow \infty$ to $(\mathcal{I}_{r,s}(x))_{x\in \Z_{\geq 0}}$.
        \item(maximal-current) For $-\log(r) \geq 0$ and $\log(s) \geq 0$, the process
\[
    \left(G(x+N, N) - G(N,N)\right)_{x\in \Z_{\geq 0}}
\]
 
        converges weakly as $N \rightarrow \infty$ to $(\mathcal{I}_{r,1}(x))_{x\in \Z_{\geq0}}$.
        \item(low-density) For $-\log(r) \leq 0$ and $\log(s) \geq -\log(r),$ the process
\[
    \left(G(x+N, N) - G(N,N)\right)_{x\in \Z_{\geq 0}}
\]
 converges weakly as $N \rightarrow \infty$ to $(\mathcal{I}_{r,1/r}(x))_{x\in \Z_{\geq 0}}$, i.e., a Geom$(\sqrt{q}/r)$ random walk.
    \end{itemize}
\end{prop}

Despite this structural understanding, the lack of explicit integrable formulas has limited progress on identifying fluctuation laws in the stationary half-space setting. Previous work succeeded in treating special cases using analytic continuation and shift arguments, see \cite{BEFO22,BeteaFerrariOccelli2020,zeng,zeng2026}, but in the genuinely two-parameter stationary regime an additional random shift is required to preserve direct Pfaffian integrability and obstructing exact distribution. We resolve this issue.

Although stationarity holds for all $rs<1$, the fluctuation behavior of the stationary last passage time depends sensitively on the location of the parameters within the stationary phase. Of particular interest is the behavior of the model in the vicinity of the critical boundary 
$rs=1$. As is apparent from exact formulas, the diagonal last passage time exhibits qualitatively different scaling behavior depending on how close the parameters are to this critical point.
{In this work, we focus on the critical stationary regime in which $rs\to1$ as $N\to\infty$. More precisely, we tune the boundary parameters on the KPZ scale as
\[
r=1+\rr(c_0N/2)^{-1/3},\qquad
s=1+\s(c_0N/2)^{-1/3},
\]
where
\[
c_0=\frac{2\sqrt q(1+\sqrt q)}{(1-\sqrt q)^3}.
\]
Under this critical scaling, the boundary parameters contribute at the same scale as the $N^{1/3}$ fluctuations of the last-passage time. Our aim is to determine the corresponding limiting distributions in the different stationary regimes.}

Under this scaling, the diagonal last-passage time admits a non-degenerate scaling limit. It is expected that our formula gives the one-point distribution of the stationary half-space KPZ fixed point. This is not solved because there has been no formula for half-space KPZ fixed point with general random initial condition.
\smallskip

\noindent\textbf{\textcolor{blue}{Goals and methods.}}
{The purpose of the present paper is twofold. First, we derive an exact finite-size formula for the stationary diagonal last-passage time in the parameter domain covered by the two-parameter construction, and treat the product-stationary specialization separately. Second, we analyze these formulas as the parameters approach the critical point under KPZ scaling. In the genuinely two-parameter stationary regime the boundary increments are correlated, whereas the product-stationary specialization is treated by the separate construction developed later in the paper.}

Exact one-point formulas have previously been obtained for several two-parameter stationary half-space KPZ models, including half-space stationary geometric LPP \cite{zeng}, the half-space stationary log-gamma polymer \cite{zeng2026}, and the half-space stationary KPZ equation \cite{barraquandKPZ,barraquandKPZBrownian}. We explain the occurrence of a random shift using geometric LPP. In this model, the integrable initial condition \(\mathcal{I}(\cdot)\) encodes the stationary measure through the gap process

$$
\bigl(\mathcal{I}(x)-\mathcal{I}(0)\bigr)_{x\geq 0},
$$

so that the random value \(\mathcal{I}(0)\) appears as an additional additive offset in the corresponding one-point formula. In the product-stationary case, this offset is independent of the rest of the model and can therefore be removed completely by a shift argument. For genuinely two-parameter stationary models, however, the offset is correlated with the remaining initial data and cannot be separated in this way.

The key idea of the present paper is to construct a representation in which the contribution of this offset becomes independent of the remaining LPP environment. To achieve this, we embed the model into a larger inhomogeneous half-space geometric LPP model with carefully chosen parameters. In the enlarged model, the stationary diagonal last-passage time is represented as the maximum of two jointly integrable last-passage times minus an independent geometric random variable. This geometric variable isolates and removes the random offset; because it is independent of the two last-passage times, its contribution can be extracted explicitly from the distribution. Conditioning on this variable reduces the stationary distribution to a geometric mixture of joint distribution functions admitting explicit Fredholm Pfaffian formulas. The max--minus identity in Proposition \(\ref{prop:maxkernel}\) then evaluates this mixture and yields a closed finite-size formula for the unshifted two-parameter stationary observable. This gives, to our knowledge, the first exact finite-size Fredholm Pfaffian formula for the stationary diagonal last-passage time itself, valid throughout the high-density and maximal-current regimes.

Our construction is related in spirit to recent work on critical fluctuations in LPP with thick boundaries \cite{EPM2025}, where embedding into a larger model is also used to recover integrability. The obstructions in the two settings are nevertheless different. In the thick-boundary problem, the embedding addresses the additional boundary geometry, whereas in the present work it resolves the random-offset obstruction arising from stationary half-space initial data. Our embedding is therefore specifically adapted to the gap structure of stationary half-space models and produces an exact formula for the stationary observable itself.

The resulting finite-size Pfaffian formula is also well suited for asymptotic analysis. In Theorem \ref{thm: high-density phase theorem}, we perform a steepest-descent analysis in the critical regime \(r,s\to 1\) under KPZ \(N^{1/3}\)-scaling. We prove that the centered and scaled stationary diagonal last-passage time converges to a Fredholm Pfaffian distribution governed by an Airy-type matrix kernel. We identify this limiting law with the one-point distribution of the half-space Airy\(_{\mathrm{stat}}\) process at the diagonal.

The same method also yields a more direct formula in the product-stationary specialization. In this case, the distribution is expressed as a single Fredholm Pfaffian, rather than as the difference of two Fredholm Pfaffians that arises from the shift argument. We verify that this simpler formula agrees with the previously obtained, more complicated expression derived from the half-space stationary log-gamma polymer model.

From a probabilistic perspective, our results complement earlier work on stationary KPZ fluctuations in full space and on half-space KPZ growth with non-stationary boundary conditions. Although the present paper is restricted to diagonal one-point distributions, the embedding method is not intrinsically limited to this observable. We expect that it will also be useful for studying off-diagonal and multipoint distributions, as well as more general stationary scaling limits in half-space KPZ models.

\subsection{Main results} \phantom{AAAAAAAAAAA}
\smallskip

\noindent\textbf{Exact finite-size Pfaffian formula.}
Our first result expresses the stationary diagonal last-passage time as the maximum of two jointly integrable last-passage times minus an independent geometric random variable. This representation yields a geometric mixture of joint distribution functions, which can be evaluated as a single Fredholm Pfaffian with an explicit \(2\times2\) matrix kernel given by contour integrals. {The resulting formula is valid at finite size in the parameter range established by the analytic-continuation argument below. Its simpler structure also allows us to take the limit \(s\to1\), thereby obtaining a formula in the maximal-current phase, which is difficult to access directly from the earlier, more complicated formulas.}
\smallskip 

\noindent\textbf{Removal of the residual random shift.} A key novelty of the present work is the complete removal of the residual random shift that appears in earlier representations of stationary half-space LPP. Unlike previous approaches based on deconvolution, our embedding method absorbs the shift into a joint integrable structure, yielding a formula directly for the stationary observable itself. This resolves a structural obstruction that previously prevented both a closed-form expression and a direct asymptotic analysis. 

\smallskip

{\noindent\textbf{KPZ scaling limit.} We perform a steepest-descent analysis of the kernel and prove convergence of the Fredholm Pfaffian to a limiting Pfaffian distribution expressed in terms of Airy-type integrals.}

\smallskip
\noindent\textbf{Matching of formulas.}
We establish a direct matching between two finite-$N$ distribution formulas
in the product-stationary case through a kernel transformation. These
formulas arise from the half-space log-gamma polymer and half-space
geometric LPP, whose Pfaffian kernels have markedly different structures.
This matching is established on the degeneration to exponential LPP. It shows that the shift argument can be encoded directly
in the Fredholm Pfaffian representation. We further show that our
high-density formula reduces to the low-density formula as
$\tilde{r}\to\pm\tilde{s}$, corresponding to the degeneration of the
two-parameter stationary initial condition to a product-stationary
initial condition.

\subsection{Organization of the paper.}
Section~\ref{paper:main-results} states our one-point limiting formulas in the high-density, maximal-current, and product-stationary/low-density regimes.
Section~\ref{paper:hd-distribution} explains the Fredholm Pfaffian structure of the half-space geometric LPP model and provides a reformulation of the kernel.
Section~\ref{paper:stationary-limit} performs an analytic continuation to reach the two-parameter stationary model.
Section~\ref{paper:hd-asymptotics} proves the asymptotic formula for the high-density phase and Theorem~\ref{thm: high-density phase theorem}.
Section~\ref{sec: max current} provides the analogous analysis in the maximal-current phase and proves Theorem~\ref{thm: maximal-current phase theorem}.
Section~\ref{section: low} proves the finite-size and asymptotic formulas for the product-stationary model, including Theorem~\ref{thm: low-density phase theorem}.
Section~\ref{sec:matching-low-density} matches the product-stationary formula with the formula obtained from the log-gamma polymer.
Finally, Section~\ref{sec:boundary-matching} proves that the high-density formula converges to the low-density formula at both $\rr=-\s$ and $\rr=\s$.

\section{Main results}\label{paper:main-results}
Given any initial condition $(\mathbf{G}(N,1))_{N\in \Z_{\geq 1}}$ and random variables $\omega_{N,M}$ for all $N\geq M$, $N,M \in \Z_{\geq 2}$, we define the LPP recurrence relation for $\mathbf{G}(N,M)$ as follows:
\begin{equation}\label{recurrence}
    \begin{cases}
        \mathbf{G}(N,M) = \omega_{N,M} + \max(\mathbf{G}(N-1,M),\, \mathbf{G}(N,M-1)), &\text{ for }N>M\geq 2,\\
        \mathbf{G}(N,N) = \omega_{N,N} + \mathbf{G}(N,N-1), &\text{ for } N\geq 2.
    \end{cases}
\end{equation}
We say a random variable $\omega$ is distributed as $\mathrm{Geom}(q)$ if $\Pb(\omega = k) = (1-q)q^k,$ for all $k \in \Z_{\geq 0}.$
Let $\sqrt{q}\in (0,1)$ and $r\in (0,1/\sqrt{q})$ be two parameters, where $\sqrt{q}$ is called the bulk parameter and $r$ is called the boundary parameter. We assume the weights are all independent and distributed as follows:
\begin{equation}\label{weights}
\begin{cases}
    \omega_{N,M} \sim \text{Geom}(q), &\text{ for } N> M\geq 2,\\
    \omega_{N,N} \sim \text{Geom}(r\sqrt{q}), &\text{ for } N\geq 2.
\end{cases}
\end{equation}

\begin{thm}[High-density phase]\label{thm: high-density phase theorem}
For any $\sqrt{q} \in (0,1)$, $\tilde{s} <0$ and $\tilde{r} \in (-\infty, \tilde{s})\cup (\tilde{s}, -\tilde{s})$, we consider the scaling $s = 1 + \s/(c_0N/2)^{1/3}$ and $r = 1 + \rr/(c_0N/2)^{1/3}$. Let $\mathbf{G}$ denote the solution to \eqref{recurrence} with the initial condition $\mathbf{G}(1+\cdot,1) - \mathbf{G}(1,1) \stackrel{(d)}{=} \mathcal{I}_{r,s}(\cdot)$. Then we have
    \begin{equation}\label{paper:hd-limit}
    \lim_{N\rightarrow \infty}\Pb\left( \frac{\mathbf{G}({N+2,N+2})-\frac{2\sqrt{q}}{1-\sqrt{q}}N}{(c_0N/2)^{1/3}} \leq \tilde{d}\right) = -(\rr + \s)\int_{0}^{\infty}\Pf(J- \widetilde{K})_{\mathbb{L}^2(\widetilde{T})} e^{(\rr + \s)\tilde{y}}d\tilde{y},
\end{equation}
where $\widetilde{T} = \{1\} \times (\tilde{d}, \infty) \cup \{2\} \times (\tilde{d}+\tilde{y}, \infty).$
\end{thm}
{\noindent\textbf{Parameter-boundary note.} The values $\rr=\s$ and $\rr=-\s$ are excluded from this theorem. They correspond to product-stationary specializations/boundaries. The matching of formulas under these two cases is proved in Section \ref{sec:boundary-matching}.}

We define the kernel $\widetilde{K}$ used in the above theorem. Let $\widetilde{K} = \begin{pmatrix}
    \widetilde{K}_{11} & \widetilde{K}_{12}\\
    \widetilde{K}_{21} & \widetilde{K}_{22}\\
\end{pmatrix}.$
For $n,n' \in \{1,2\}$ and $X,Y\in \mathbb{R}$, we define
\begin{equation}\label{paper:hd-ab}
    \begin{aligned}
        &\A(1,\zeta) := (\zeta - \s), \quad \A(2,\zeta) := (\zeta - \s)(\zeta + \s),\\
        &\B(1,\zeta) := -(\zeta + \s), \quad \B(2,\zeta) := (\zeta + \s)(\zeta - \s),
    \end{aligned}
\end{equation}
    \begin{equation}\label{paper:hd-kernel}
    \begin{aligned}
            \widetilde{K}_{11}(n,X;n',Y) := &-  \!\int\limits_{{}_{}\wcu\, {}_{0}}\! \frac{d\omega}{2\pi\I}\!\int\limits_{{ }_{0}\zcd } \!\frac{d\zeta}{2\pi\I}e^{\frac{\zeta^3}{3} - \frac{\omega^3}{3}  - X\zeta + Y\omega} \A(n,\zeta)\B(n',\omega)\frac{(\zeta + \omega)(\zeta - \tilde{r})(\omega + \tilde{r})}{4\zeta\omega(\zeta - \omega)}, \\
            \widetilde{K}_{12}(n,X;n',Y) := &-  \!\int\limits_{{}_{\tilde{r},-\s,\s}\wcu\, {}_{}}\! \frac{d\omega}{2\pi\I}\!\int\limits_{{ }_{0, \omega}\zcd } \!\frac{d\zeta}{2\pi\I} e^{\frac{\zeta^3}{3} - \frac{\omega^3}{3}  - X\zeta + Y\omega} \frac{\A(n,\zeta)}{\A(n',\omega)}\frac{(\zeta + \omega)(\zeta - \tilde{r})}{(2\zeta)(\omega - \tilde{r})(\zeta - \omega)} \\
            &- \Id_{\{n = 1, n' = 2, X\leq Y\}} e^{-\s(Y-X)},\\
            \widetilde{K}_{21}(n,X;n',Y) := & -\widetilde{K}_{12}(n',Y;n,X),\\
            \widetilde{K}_{22}(n,X;n',Y) := & \,\,\widetilde{h}_{22}(n,X;n',Y) + \widetilde{R}(n,X;n',Y)\Id_{X\leq Y} - \widetilde{R}(n',Y;n,X)\Id_{X> Y}.
    \end{aligned}
\end{equation}
We define functions $\widetilde{h}_{22}$ and $\widetilde{R}$ as follows.
\begin{equation}\label{paper:hd-h22}
    \begin{aligned}
        &\widetilde{h}_{22}(n,X;n',Y) :=\!\!\!\!\!\int\limits_{{}_{}\wcu\, {}_{\tilde{r},-\s,\zeta,0}}\! \!\!\!\!\frac{d\omega}{2\pi\I}\!\!\!\int\limits_{{ }_{}\zcd {}_{\s,-\tilde{r}, -\tilde{s}} } \!\!\!\!\!\frac{d\zeta}{2\pi\I} \widetilde{F}(\zeta, \omega)+\!\!\!\!\int\limits_{\Gamma_{\tilde{r},\s,-\s}}\!\!\!\! \frac{d\omega}{2\pi\I}\!\!\!\!\int\limits_{{ }_{\s,-\tilde{r}, -\tilde{s},\omega}\zcd {}_{} } \!\!\!\!\!\!\frac{d\zeta}{2\pi\I} \widetilde{F}(\zeta, \omega)+ \widetilde{G}(n,X;n',Y) + \widetilde{O}(n,X;n',Y),
    \end{aligned}
\end{equation}
where
\begin{equation}\label{eq: integrand of 22 limit}
    \widetilde{F}(\zeta, \omega) := e^{\frac{\zeta^3}{3} - \frac{\omega^3}{3}  - X\zeta + Y\omega}\frac{1}{\A(n',\omega)\B(n,\zeta)}\frac{(\zeta + \omega)}{(\zeta - \omega)(\zeta + \tilde{r})(\omega - \tilde{r})}.
\end{equation}
\begin{equation}\label{eq: definitino of G tilde}
    \begin{aligned}
    \widetilde{G}(n,X;n',Y):= &-\frac{1}{\s^2-\rr^2}e^{-\frac{\rr^3}{3} + Y \rr - \frac{\s^3}{3} +X\s}\left(\Id_{n'=1} + \frac{1}{\s+\rr} \Id_{n'=2}\right)\left(\Id_{n=1} + \frac{1}{2\s} \Id_{n=2}\right)\\
        &+ \frac{1}{\s^2-\rr^2}e^{-\frac{\rr^3}{3} + X \rr - \frac{\s^3}{3} +Y\s}\left(\Id_{n'=1} + \frac{1}{2\s} \Id_{n'=2}\right)\left(\Id_{n=1} + \frac{1}{(\s+\rr)} \Id_{n=2}\right).\\
    \end{aligned}
\end{equation}

\begin{equation}\label{eq: def of O}
    \begin{aligned}
        \widetilde{O}(n,X;n',Y) :=&\Id_{n=2} \frac{1}{(\s-\rr)^2(-2\s)}e^{-\frac{\rr^3}{3} +Y\rr + \frac{\s^3}{3} - X\s}\left(\Id_{n'=1} + \frac{1}{\s+\rr} \Id_{n'=2}\right)\\
        &+\Id_{n'=2} \frac{1}{(\s-\rr)^2(2\s)}e^{-\frac{\rr^3}{3} +X\rr + \frac{\s^3}{3} - Y\s}\left(\Id_{n=1} + \frac{1}{\s+\rr} \Id_{n=2}\right).
    \end{aligned}
\end{equation}

\begin{equation}\label{paper:hd-residue}
    \begin{aligned}
        \widetilde{R}(2,X;2,Y) :=& \frac{e^{\rr(Y-X)}}{(\s - \rr)^2(\s+\rr)^2}+ \frac{e^{\s(Y - X)}}{(\s - \rr)(\s+\rr)^2(-2\s)}- \frac{e^{\s(X-Y)}}{(\s-\rr)^2(\s+\rr)(2\s)}\\
        &+\frac{e^{\s(X-Y)}}{(\s-\rr)(\s+\rr)(2\s)}\left(X -\s^2\right) +  \frac{e^{\s(Y-X)}}{(\s-\rr)(\s+\rr)(2\s)}\left(-X +\s^2\right),\\
        \widetilde{R}(1,X;2,Y) :=& \frac{e^{\rr(Y-X)}}{(\s+\rr)^2(\s-\rr)}-\frac{e^{\s(Y-X)}}{(\s+\rr)(2\s)(\s-\rr)} -{\frac{e^{\s(X-Y)}}{(\s-\rr)^2(2\s)}}
        -\frac{e^{\s(X-Y)}}{(\s-\rr)(\s+\rr)}\left(-X+\s^2\right),\\
        \widetilde{R}(2,X;1,Y) :=& \frac{e^{\rr(Y-X)}}{(\s-\rr)^2(\s+\rr)}-\frac{e^{\s(Y-X)}}{(\s-\rr)(\s+\rr)^2}-\frac{e^{\s(Y-X)}}{(\s-\rr)(\s+\rr)}\left(X-\s^2\right),\\
        \widetilde{R}(1,X;1,Y) :=& \frac{e^{\rr(Y-X)} - e^{\s(Y-X)}}{(\s -\rr)(\s + \rr)}.
    \end{aligned}
\end{equation}

\begin{thm}[Maximal-current phase]\label{thm: maximal-current phase theorem}
{Fix any $\sqrt{q} \in (0,1)$ and $\tilde{r} \in (-\infty,0)$. Set $s=1$ and consider the scaling $r=1+\rr(c_0N/2)^{-1/3}$.} Let $\mathbf{G}$ denote the solution to \eqref{recurrence} with the initial condition $\mathbf{G}(1+\cdot,1) - \mathbf{G}(1,1) \stackrel{(d)}{=} \mathcal{I}_{r,1}(\cdot)$. Then we have
    \begin{equation}
    \lim_{N\rightarrow \infty}\Pb\left( \frac{\mathbf{G}({N+2,N+2})-\frac{2\sqrt{q}}{1-\sqrt{q}}N}{(c_0N/2)^{1/3}} \leq \tilde{d}\right) = -\rr\int_{0}^{\infty}\Pf(J- \widetilde{K}_{M})_{\mathbb{L}^2(\widetilde{T})} e^{\rr\tilde{y}}d\tilde{y},
\end{equation}
where $\widetilde{T} = \{1\} \times (\tilde{d}, \infty) \cup \{2\} \times (\tilde{d}+\tilde{y}, \infty).$
\end{thm}
We define the kernel $\widetilde{K}_M$ used in the above theorem. Let $\widetilde{K}_{M} = \begin{pmatrix}
    \widetilde{K}_{11}^M & \widetilde{K}_{12}^M\\
    \widetilde{K}_{21}^M & \widetilde{K}_{22}^M\\
\end{pmatrix}$.  We define
\begin{equation}
    \begin{aligned}
        &\widetilde{A}_M(1,\zeta) = \zeta, \quad \widetilde{A}_M(2,\zeta) = \zeta^2,\quad
        \widetilde{B}_M(1,\zeta) = -\zeta, \quad \widetilde{B}_M(2,\zeta) = \zeta^2.\\
    \end{aligned}
\end{equation}
\begin{equation}
    \begin{aligned}
        \widetilde{K}_{11}^M(n,X;n',Y) := &-  \!\int\limits_{{}_{}\wcu\, {}_{0}}\! \frac{d\omega}{2\pi\I}\!\int\limits_{{ }_{0}\zcd } \!\frac{d\zeta}{2\pi\I}e^{\frac{\zeta^3}{3} - \frac{\omega^3}{3}  - X\zeta + Y\omega} \A_M(n,\zeta)\B_M(n',\omega)\frac{(\zeta + \omega)(\zeta - \tilde{r})(\omega + \tilde{r})}{4\zeta\omega(\zeta - \omega)}, \\
            \widetilde{K}_{12}^M  (n,X;n',Y) := &-  \!\int\limits_{{}_{\tilde{r},0}\wcu\, {}_{}}\! \frac{d\omega}{2\pi\I}\!\int\limits_{{ }_{0, \omega}\zcd } \!\frac{d\zeta}{2\pi\I} e^{\frac{\zeta^3}{3} - \frac{\omega^3}{3}  - X\zeta + Y\omega} \frac{\A_M(n,\zeta)}{\A_M(n',\omega)}\frac{(\zeta + \omega)(\zeta - \tilde{r})}{(2\zeta)(\omega - \tilde{r})(\zeta - \omega)}, \\
            &- \Id_{n = 1, n' = 2, X\leq Y} ,\\
            \widetilde{K}_{21}^M(n,X;n',Y) := & -\widetilde{K}_{12}^M(n',Y;n,X).\\
    \end{aligned}
\end{equation}

\begin{equation}\widetilde{K}_{22}^M(n,X;n',Y) := \widetilde{h}_{22}^M(n,X;n',Y) + \widetilde{\mathcal{A}}^M(n,X;n',Y)\Id_{X\leq Y} - \widetilde{\mathcal{A}}^M(n',Y;n,X)\Id_{X> Y}
\end{equation}
where
\begin{equation}
    \begin{aligned}
        \widetilde{h}_{22}^M(n,X;n',Y) := &\int\limits_{{}_{}\wcu\, {}_{\tilde{r},\zeta,0}}\! \frac{d\omega}{2\pi\I}\!\int\limits_{{ }_{}\zcd {}_{0,-\tilde{r}} } \!\frac{d\zeta}{2\pi\I} e^{\frac{\zeta^3}{3} - \frac{\omega^3}{3}  - X\zeta + Y\omega}\frac{1}{\A_M(n',\omega)\B_M(n,\zeta)}\frac{(\zeta + \omega)}{(\zeta - \omega)(\zeta + \tilde{r})(\omega - \tilde{r})}\\
        &+\int\limits_{{}_{\tilde{r},0}\wcu\, {}_{}}\! \frac{d\omega}{2\pi\I}\!\int\limits_{{ }_{-\tilde{r}, 0,\omega}\zcd {}_{} } \!\frac{d\zeta}{2\pi\I} e^{\frac{\zeta^3}{3} - \frac{\omega^3}{3}  - X\zeta + Y\omega}\frac{1}{\A_M(n',\omega)\B_M(n,\zeta)}\frac{(\zeta + \omega)}{(\zeta - \omega)(\zeta + \tilde{r})(\omega - \tilde{r})},\\
    \end{aligned}
\end{equation}
and 
\begin{equation}
    \begin{aligned}
        \widetilde{\mathcal{A}}^M(n,X;n',Y) &= \Id_{n=n'=1}\widetilde{R}_{11}^M(X,Y) + \Id_{n=1,n'=2}\widetilde{R}_{12}^M(X,Y)\\
        &\qquad{}+ \Id_{n=2,n'=1}\widetilde{R}_{21}^M(X,Y) + \Id_{n=n'=2}\widetilde{R}_{22}^M(X,Y),
    \end{aligned}
\end{equation}
with
\begin{equation}
    \begin{aligned}
        \widetilde{R}_{11}^M(X,Y) &:= \frac{1}{\rr^2} -\frac{e^{\rr(Y-X)}}{\rr^2} + \frac{1}{\rr^2}\left(e^{-\frac{\rr^3}{3} + Y\rr} - e^{-\frac{\rr^3}{3} + X\rr}\right),\\
        \widetilde{R}_{12}^M(X,Y) &:= \frac{1}{\rr^3}e^{-\frac{\rr^3}{3} +Y\rr} + e^{-\frac{\rr^3}{3} +X\rr}\left(\frac{1}{\rr^3} - \frac{Y}{\rr^2}\right) + \frac{e^{\rr(Y-X)} +1}{-\rr^3}  +\frac{(Y-2X)}{\rr^2},\\
        \widetilde{R}_{21}^M(X,Y) &:= \frac{e^{\rr(Y-X)}}{\rr^3} + \frac{1}{\rr^3} + \frac{X}{\rr^2} -\frac{1}{\rr^3}e^{-\frac{\rr^3}{3} +X\rr} + e^{-\frac{\rr^3}{3} +Y\rr}\left(-\frac{1}{\rr^3} + \frac{X}{\rr^2}\right),\\
        \widetilde{R}_{22}^M(X,Y) &:= \frac{e^{\rr(Y-X)}}{\rr^4} -\frac{1}{\rr^4}\left(1+\rr(X-Y) + \rr^2X(X-Y)\right)\\
        &\qquad{}+ e^{-\frac{\rr^3}{3} +Y\rr} \left(\frac{-1}{\rr^4} + \frac{X}{\rr^3}\right) + e^{-\frac{\rr^3}{3} +X\rr} \left(\frac{1}{\rr^4} - \frac{Y}{\rr^3}\right).
    \end{aligned}
\end{equation}

\begin{thm}[Low-density phase]\label{thm: low-density phase theorem}
For any $\sqrt{q} \in (0,1)$, $\tilde{r} \leq 0$, we consider the scaling $r = 1 + \rr/(c_0N/2)^{1/3}$. Let $\mathbf{G}$ denote the solution to \eqref{recurrence} with the initial condition $\mathbf{G}(j+1,1) - \mathbf{G}(j,1) \stackrel{(d)}{=} \text{Geom}(\sqrt{q}/r)$ for $j \geq 1$ and $\mathbf{G}(1,1) = 0$. Then we have
    \begin{equation}\label{paper:ld-limit}
    \lim_{N\rightarrow \infty}\Pb\left( \frac{\mathbf{G}({N,N})-\frac{2\sqrt{q}}{1-\sqrt{q}}N}{(c_0N/2)^{1/3}} \leq \tilde{d}\right) = \Pf(J- \tKK)_{\mathbb{L}^2(\widetilde{T})},
\end{equation}
where $\widetilde{T} = (\tilde{d}, \infty).$
\end{thm}
We define the kernel $\tKK$ used in the above theorem.
\begin{equation}\label{paper:ld-kernel}
    \begin{aligned}
        \tKK_{11}(X,Y) &:=\int\limits_{{}_{}\wcu\, {}_{0}}\! \frac{d\omega}{2\pi\I}\!\int\limits_{{ }_{0}\zcd {}_{}}\!\frac{d\zeta}{2\pi\I} e^{\frac{\zeta^3}{3} - \frac{\omega^3}{3}  - X\zeta + Y\omega}\frac{(\zeta^2 - \rr^2)(\omega^2 - \rr^2)(\zeta + \omega)}{4\zeta \omega(\zeta - \omega)},\\
        \tKK_{12}(X,Y) &:=-\int\limits_{{}_{\rr,-\rr}\wcu\, {}_{}}\! \frac{d\omega}{2\pi\I}\!\int\limits_{{ }_{\omega}\zcd {}_{}}\!\frac{d\zeta}{2\pi\I} e^{\frac{\zeta^3}{3} - \frac{\omega^3}{3}  - X\zeta + Y\omega}\frac{(\zeta + \rr)(\zeta + \omega)(\zeta - \rr)}{2\zeta(\omega^2 - \rr^2)(\zeta - \omega)},\\
        \tKK_{22}(X,Y) &:= \widetilde{\mathcal{Q}}(X,Y) + \widetilde{\mathsf{F}}(X,Y)\Id_{X\leq Y} - \widetilde{\mathsf{F}}(Y,X)\Id_{X> Y},
    \end{aligned}
\end{equation}
where
\begin{equation}\label{paper:ld-residue}
    \begin{aligned}
        \widetilde{\mathcal{Q}}(X,Y) &:= -\int\limits_{{}_{}\wcu\, {}_{\zeta}}\! \frac{d\omega}{2\pi\I}\!\int\limits_{{ }_{}\zcd {}_{-\rr,\rr}}\!\frac{d\zeta}{2\pi\I} e^{\frac{\zeta^3}{3} - \frac{\omega^3}{3}  - X\zeta + Y\omega}\frac{(\zeta + \omega)}{(\zeta^2 - \rr^2)(\omega^2 - \rr^2)(\zeta - \omega)}\\
         &\quad\,\, -\int\limits_{\Gamma_{-\tilde{r}, \tilde{r}}}\! \frac{d\omega}{2\pi\I}\!\int\limits_{{ }_{\omega}\zcd {}_{}}\!\frac{d\zeta}{2\pi\I} e^{\frac{\zeta^3}{3} - \frac{\omega^3}{3}  - X\zeta + Y\omega}\frac{(\zeta + \omega)}{(\zeta^2 - \rr^2)(\omega^2 - \rr^2)(\zeta - \omega)},\\
        \widetilde{\mathsf{F}}(X,Y) &:= \frac{X -\rr^2}{2\rr}\left(e^{\rr(Y-X)}- e^{\rr(X-Y)}\right).
    \end{aligned}
\end{equation}

\section{Diagonal distribution of the stationary LPP under high-density regime}\label{paper:hd-distribution}
\begin{figure}
\centering
\begin{tikzpicture}
\draw[very thin, gray] (-0.5,-0.5) grid (4.5,4.5);
    \foreach \x in {1,...,4} {
        \foreach \y in {0,1,...,\x} {
            \fill (\x,\y) circle (2pt);
        }
    }
    \fill (0,0) circle (2pt);
    \draw[thick,-] (0,0) -- (4,0);
    \draw[thick,-] (0,0) -- (4,4);
    \draw[thick,-] (1,1) -- (4,1);
    \fill (0,0) circle (1pt);
    \node at (0,-0.25) {$0$};
     
    \node at (1,-0.25) {$-X$};
    \node at (-0.5,0) {\footnotesize{$(1,1)$}};
    \node at (0.5,1) {\footnotesize{$(2,2)$}};
    \node at (1.2,0.82) {\small{$X$}};
    \fill (1,1) circle (1pt);
    \node at (3,-0.28) {\small{$\sqrt{q}/s$}};
    \node at (3,0.72) {\small{$\sqrt{q}s$}};
    \node at  (1.55,2.1) {\small{$r\sqrt{q}$}};
    \node at (3.1,1.75) {\small{$q$}};
\end{tikzpicture}
\caption{Geometric LPP model as described in $\eqref{eq:stat_wts}$.}
\label{fig: two param stationary version}
\end{figure}
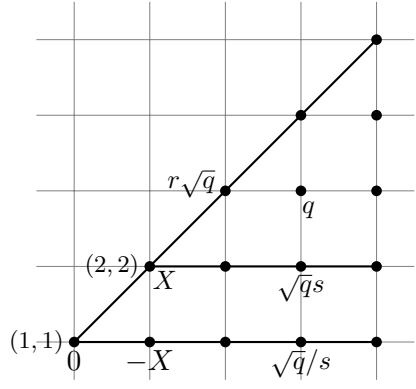

\begin{figure}[htbp]
\centering
\begin{tikzpicture}
\draw[very thin, gray] (-0.5,-0.5) grid (5.5,5.5);
    \foreach \x in {1,...,5} {
        \foreach \y in {0,1,...,\x} {
            \fill (\x,\y) circle (2pt);
        }
    }
    \fill (0,0) circle (2pt);
    \draw[thick,-] (0,0) -- (5,0);
    \draw[thick,-] (0,0) -- (5,5);
    \draw[thick,-] (4,4) -- (4,0);
    \draw[thick,-] (5,5) -- (5,0);
    \draw[thick,red] (3,3) -- (5,3);
    \draw[thick,red] (5,3) -- (5,0);
    \draw[thick,red] (0,0) -- (3,3);
    \fill (0,0) circle (1pt);
    \fill[green] (4,3) circle (1pt);
    \fill[green] (5,3) circle (1pt);

    \node at (-0.5,0) {\footnotesize{$(1,1)$}};
    \node at (6.1,5) {\footnotesize{$(N+2,N+2)$}};
    \node at (2.4,3) {\footnotesize{$(N,N)$}};
    \node at (5.8,3) {\footnotesize{$(N+2,N)$}};
    \node at (4.8,5) {\footnotesize{$0$}};

    \node at (3.8,4) {\small{$X$}};
    \fill (1,1) circle (1pt);
    \node at (5.4,4) {\small{$-X$}};
    \fill (1,1) circle (1pt);
    \node at (5.3,2) {\small{$\sqrt{q}t$}};
    \fill (1,1) circle (1pt);

    \node at (4.3,2) {\small{$\sqrt{q}s$}};
    \fill (1,1) circle (1pt);

    \node at  (1.55,2.1) {\small{$r\sqrt{q}$}};
    \node at (3.1,1.75) {\small{$q$}};
\end{tikzpicture}
\caption{Two parameter geometric LPP model.}
\label{fig: two param stationary flipped version}
\end{figure}
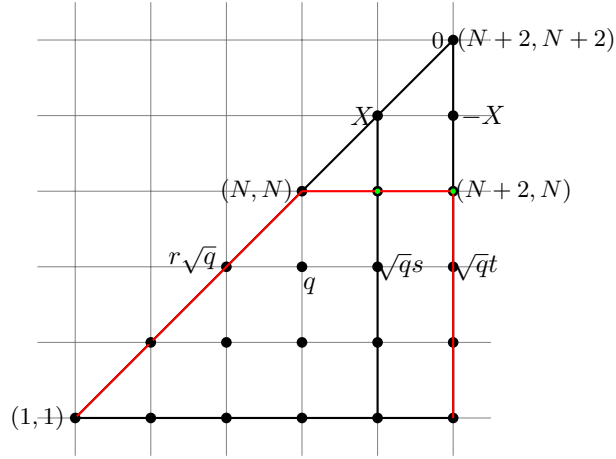

\subsection{Inhomogeneous half-space geometric last passage percolation model}
We provide the definition (path version) of the half-space geometric last passage percolation (LPP) model.
Let $q_0, q_1,q_2,
\dots$ be positive real parameters such that $q_0q_i \in (0,1)$ for all $i \geq 1$ and $q_iq_j \in (0,1)$ for all $i> j \geq 1.$ We define a set $\mathbb{H} = \{ (i,j) \in \Z^2 | 1 \leq j \leq i\}.$ Let $(\omega_{i,j})_{(i,j)\in \mathcal{D}}$ be independent random variables with $\omega_{i,j}\sim \text{Geo}(q_iq_j)$ for $i>j\geq 1$ and $\omega_{i,i}\sim \text{Geo}(q_0q_i)$ for $i \geq 1$. We call $q_0$ the diagonal parameter and $q_i$ the ith-row parameter. An \emph{up-right path} $\pi$ from a point $A$ to a point $E$ is a sequence of points $(\pi(0),\pi(1),\ldots,\pi(n))$ in $\Z^2$ such that \mbox{$\pi(k+1)-\pi(k)\in \{(0,1),(1,0)\}$}, with $\pi(0)=A$ and $\pi(n)=E$, and where $n$ is called the length $\ell(\pi)$ of $\pi$.
We define the half-space last passage time from $(1,1)$ to the point $(N,M)$ for $N\geq M$, denoted as $G(N,M)$, to be
\begin{equation}
G(N,M) = \max_{\pi:(1,1)\to (N,M)} \sum_{(i,j)\in\pi} \omega_{i,j}
\end{equation}
where the maximum is taken over all up-right paths in $\mathbb{H}$ from $(1,1)$ to $(N,M)$. 
To relate this definition to the recurrence relation \eqref{recurrence}, ${G}(N,M)$ is the solution to \ref{recurrence} with initial condition $G(1,1) = \omega_{1,1}$ and $G(N,1) - G(N-1,1) = \omega_{N,1}$ for $N \geq 2.$

\subsection{From inhomogeneous model to the two parameter stationary case}
Fix
$\sqrt{q} \in (0,1)$, $s \in (\sqrt{q}, 1)$, and $r\in (0,1/s),$ $r\neq \sqrt{q}$.
We set $q_0 = r$, $q_1 = 1/s$, $q_2 = s$ and $q_i = \sqrt{q}$ for $i\geq 3.$ Let $X$ be an independent $\text{Geo}(rs)$ random variable. We set $\omega_{1,1} = 0,$ $\omega_{2,1} = -X,$ and $\omega_{2,2} = X.$ Then the two-parameter stationary LPP model has the following weights:
\begin{equation}  \label{eq:stat_wts}
  \omega_{i, j} = \begin{cases}
      0, & \textrm{if}\ i=j=1, \\
     -X, & \textrm{if}\ i=2,j=1, \\
      X, & \textrm{if }i=j=2,\\
      \mathrm{Geo}\left( r\sqrt{q} \right), & \textrm{if }i=j\geq 3,\\
    \mathrm{Geo}\left( \sqrt{q}/s \right), & \textrm{if }j=1, i\geq 3, \\
    \mathrm{Geo}\left( s\sqrt{q} \right), & \textrm{if }j=2, i\geq 3, \\
    \mathrm{Geo}(q), &\textrm{otherwise}.
  \end{cases}
\end{equation} This model is shown in Figure $\ref{fig: two param stationary version}.$ Let $G^{stat}_{r,s}(N,N)$ denote the LPP from $(1,1)$ to $(N,N)$ of the above model.

The main challenge of solving the diagonal distribution of this LPP model is the dependence of weights at $(2,1)$ and $(2,2)$. It breaks the Pfaffian structure that the half-space geometric LPP model has. Previously, \cite{zeng} solved the off-diagonal distribution of two parameter stationary model using a deconvolution method that removes weights at $(1,1)$ and $(2,1)$. However, the shortcoming of this method is that the random shift at $(2,2)$ is still preserved. We introduce the following method that removes such a random shift. As a result, we are able to obtain the exact one point distribution of half-space geometric LPP model with two parameter stationary initial condition, which directly translates to the one point distribution for stationary KPZ fixed point. Our method is inspired by \cite{EPM2025}.

We define a flipped model and specify its parameters. Let $q_0 = r,$ $q_1=\dots = q_N = \sqrt{q}$, $q_{N+1} = s$, and  $q_{N+2} = t.$ We define weights as follows:
\begin{equation} \label{eq: weigths flipped}
  \omega_{i, j} = \begin{cases}
      0, & \textrm{if}\ i=j=N+2, \\
     -X, & \textrm{if}\ i=N+2,j=N+1, \\
      X, & \textrm{if }i=j=N+1,\\
      \mathrm{Geo}\left( r\sqrt{q} \right), & \textrm{if }1\leq i=j\leq N,\\
    \mathrm{Geo}\left( s\sqrt{q} \right), & \textrm{if }i = N+1, 1\leq j\leq N, \\
    \mathrm{Geo}\left( \sqrt{q}t \right), & \textrm{if }i = N+2, 1\leq j\leq N, \\
    \mathrm{Geo}(q), &\textrm{otherwise},
  \end{cases}
\end{equation}
where $X\sim \text{Geo}(rs)$ and parameters satisfy $\sqrt{q} \in (0,1)$, $s \in (\sqrt{q}, 1)$, $r\in (0,1/s)$, $t\in (0,1/\sqrt{q})$. This model is shown in Figure $\ref{fig: two param stationary flipped version}.$ Let $G(N,M)$ denote the LPP value from $(1,1)$ to $(N,M)$ under the model \eqref{eq: weigths flipped}. It is easy to observe that
\begin{equation}
    G(N+2,N+2) \Rightarrow G_{r,s}^{stat}(N+2,N+2) \text{ weakly as }t\rightarrow 1/s. 
\end{equation}

\subsection{Distribution as Fredholm Pfaffian}
We define
\begin{equation}
    \begin{aligned}
        &A(1,x) : = (1-s/x), \quad A(2,x) := (1-s/x)(1-t/x),\quad H(x) := \frac{(1-\sqrt{q}/x)^N}{(1-\sqrt{q}x)^N}\\
        &B(1,x) := (1-sx),\quad B(2,x) := (1-sx)(1-tx).
    \end{aligned}
\end{equation}

\begin{prop}\label{prop:maxkernel}
Let $X$ be an independent Geo$(rs)$ random variable. The distribution function of $G_{r,s}^{stat}(N+2,N+2) $ is given by the two-point joint distribution:
    \begin{equation}\label{eq: equal dis}
        \Pb(G_{r,s}^{stat}(N+2,N+2) \leq d) = \lim_{t\rightarrow 1/s}\Pb\left( \max\left(G(N+1,N), G(N+2,N)-X) \right) \leq d\right),
    \end{equation}
    and we have the Fredholm Pfaffian formula
    \begin{equation}
        \Pb\left( G({N+1,N}) \leq d, G({N+2,N}) - X \leq d \right) = \sum_{y=0}^{\infty}\Pf(J- K)_{\ell^2(T)} (1-rs)(rs)^y,
    \end{equation}
    where
    \begin{equation}
    \begin{aligned}
        &\Pf\left( J - K\right)_{\ell^2(T)} = 1 + \sum_{m = 1}^{\infty}\frac{(-1)^m}{m!} \sum_{j_1,\dots, j_m = 1}^2 \sum_{x_1 = d_{j_1}}^{\infty}\cdots \sum_{x_m = d_{j_m}}^{\infty} \Pf\left(K(j_u,x_u; j_v,x_v)\right)_{u,v =1}^{m}, \\
        &T =\{1\}\times \{d+1,d+2, \dots \}\cup \{2\}\times \{d+y+1,d+y+2, \dots \},
    \end{aligned}
    \end{equation}
    and with $2\times 2$ matrix correlation kernel $K$ given by:
\begin{equation}
\begin{aligned}
    K_{11}(n,k; n^{\prime}, \ell) &= \oint \frac{dw}{2\pi\I} \oint \frac{dz}{2\pi\I}  \frac{w^{\ell-1}}{z^{k}}\frac{H(z)A(n,z)B(n',w)}{H(w)}\frac{(zw-1)(z-r)(1-wr)}{(z^2-1)(1-w^2)(z-w)},\\
    K_{12}(n,k; n^{\prime}, \ell) &= \oint \frac{dw}{2\pi\I} \oint \frac{dz}{2\pi\I}\frac{w^{\ell-1}}{z^{k}} \frac{H(z)A(n,z)}{H(w)A(n',w)}\frac{(zw-1)(z-r)}{(z^2-1)(z-w)(w-r)},\\
    K_{22}(n,k; n^{\prime}, \ell) &= \oint \frac{dw}{2\pi\I} \oint \frac{dz}{2\pi\I}\frac{w^{\ell-1}}{z^{k}} \frac{H(z)}{H(w)A(n',w)B(n,z)}\frac{(zw-1)}{(z-w)(1-zr)(w-r)}\\
    &+ \oint \frac{dz}{2\pi\I} \frac{1}{A(n',z)B(n,z)}\frac{(1-z^2)z^{\ell - k - 1}}{(1-zr)(z-r)},
  \end{aligned}
\end{equation}
where
\begin{itemize}
   \item for $K_{11}$, $\max\{t,s,\sqrt{q}\}<|w| <1$, $1 < |z|< \min \{1/t,1/s,1/\sqrt{q}\}$;
   \item for $K_{12}$, $\max \{ r,s,t,\sqrt{q}\} < |w|,$ and $1 < |z| < \min\{1/s,1/t,1/\sqrt{q}\}$ together with the condition $|w| < |z|$ if $n \geq n^{\prime}$ and $|z| < |w|$ if $n < n^{\prime}$;
   \item for $K_{22}$, $\max \{r,s,t,\sqrt{q}\} < |w|< |z| < \min\{1/r,1/s,1/t,1/\sqrt{q}\}$.
  \end{itemize}
\end{prop}

\begin{proof}
    We start with the kernel defined in \cite[Proposition 4.2]{BBCS} adapted to our choice of boundary and bulk parameters defined in \eqref{eq: weigths flipped}. Then apply the change of variable $w \rightarrow w^{-1}$.
\end{proof}

\begin{lem}
For $n,n' \in \{1,2\},$ the $11$-entry can be rewritten as
\begin{equation}\label{eq: K11 original}
    \begin{aligned}
        K_{11}(n,k; n^{\prime}, \ell) &= -\oint\limits_{\Gamma_{\sqrt{q}}} \frac{dw}{2\pi\I} \oint\limits_{\Gamma_{1/\sqrt{q}}} \frac{dz}{2\pi\I}  \frac{w^{\ell-1}}{z^{k}}\frac{H(z)A(n,z)B(n',w)}{H(w)}\frac{(zw-1)(z-r)(1-wr)}{(z^2-1)(1-w^2)(z-w)},
    \end{aligned}
\end{equation}
and the $12$-entry can be reformulated as
\begin{equation}\label{eq: K12 original}
    \begin{aligned}
        K_{12}(n,k; n^{\prime}, \ell) &= K_{12}^{reg}(n,k;n',\ell) - \Id_{\{k\leq \ell\}}\Id_{\{n = 1, n' = 2\}} t^{\ell - k},
    \end{aligned}
\end{equation}
where 
\begin{equation}
    K_{12}^{reg}(n,k;n',\ell) = -\oint\limits_{\Gamma_{r,s,t,\sqrt{q}}} \frac{dw}{2\pi\I} \oint\limits_{\Gamma_{1/\sqrt{q}}} \frac{dz}{2\pi\I}\frac{w^{\ell-1}}{z^{k}} \frac{H(z)A(n,z)}{H(w)A(n',w)}\frac{(zw-1)(z-r)}{(z^2-1)(z-w)(w-r)}.
\end{equation}
\end{lem}

\begin{proof}
    The integrand of $K_{11}$ behaves like $w^{N+\ell+n'-1}$ as $w \rightarrow 0$ and has no pole at $w = 0$ for any $\ell \geq 1.$ The integrand behaves like $z^{-(N+k+1)}$ as $z \rightarrow \infty$ and has no pole at $z = \infty.$ So we can deform the contour of $w$ to include $\sqrt{q}$ and deform the contour of $z$ to include $1/\sqrt{q}.$

    For $K_{12}$, same analysis shows that there is no pole at $w=0$ and $z = \infty$. When $n\geq n',$ we do not cross any pole. However, when $n = 1,$ $n'=2,$ we shrink the $w$-contour and deform it to lie inside the $z$-contour. Then we must add the residue $w=z$, which yields
\begin{equation}
    \begin{aligned}
        \text{Res}(w=z)= -\oint \limits_{\Gamma_{0,t}}\frac{dz}{2\pi\I} \frac{z^{\ell-k}}{(z-t)} = -t^{\ell - k}\Id_{k \leq \ell}.
    \end{aligned}
\end{equation}
Then we deform the contour of $w$ to enclose $\sqrt{q}, s, t,r$ and deform the contour of $z$ to enclose $1/\sqrt{q}.$ The pole $w = t$ is present only when $n' = 2$.
\end{proof}

To state $K_{22}$, we need to define
\begin{equation}
    F(n,k;n',\ell) := \frac{w^{\ell-1}}{z^{k}} \frac{H(z)}{H(w){A}(n',w){B}(n,z)}\frac{(zw-1)}{(z-w)(1-zr)(w-r)}.
\end{equation}
\begin{equation}
    \begin{aligned}
        h_{22}(n,k;n',\ell) :=& -\oint\limits_{\Gamma_{\sqrt{q},r}} \frac{dw}{2\pi\I} \!\!\!\oint\limits_{\Gamma_{1/\sqrt{q},1/r,1/s,1/t}} \!\!\!\!\!\!\!\!\!\frac{dz}{2\pi\I}F(n,k;n',\ell)
        -\oint\limits_{\Gamma_{s}} \frac{dw}{2\pi\I} \!\!\!\oint\limits_{\Gamma_{1/\sqrt{q},1/r}} \!\!\!\!\!\!\frac{dz}{2\pi\I}F(n,k;n',\ell)\\
        &-\oint\limits_{\Gamma_{t}} \frac{dw}{2\pi\I}\!\!\!\!\!\! \oint\limits_{\Gamma_{1/\sqrt{q},1/r}}\!\!\!\!\!\! \frac{dz}{2\pi\I}F(n,k;n',\ell),\\
    \end{aligned}
\end{equation}
and 
\begin{equation}
    g_{22}(n,k;n',\ell) := \oint \limits_{\Gamma_{0,s,t,r}}\frac{dz}{2\pi\I} \frac{1}{A(n',z)B(n,z)}\frac{(1-z^2)z^{\ell - k - 1}}{(1-zr)(z-r)}
\end{equation}
where
\begin{equation}\label{eq: residue evaluation}
    \begin{aligned}
        g_{22}(2,k;2,\ell) &= \frac{r^{\ell-k+1}}{(r-s)(r-t)(1-sr)(1-tr)}\\
        &\qquad{}+ \frac{1}{(1-st)(t-s)}\left(\frac{t^{\ell- k+1}}{(1-tr)(t-r)} - \frac{s^{\ell-k+1}}{(1-sr)(s-r)}\right),\\
        g_{22}(1,k;2,\ell) &=\frac{1}{(s-r)(1-sr)}\left(\frac{r^{\ell-k+1}}{(t-r)} - \frac{s^{\ell - k + 1}}{(t-s)}\right)\\
        &\qquad{}+ \frac{1}{(1-st)}\frac{(1-t^2)t^{\ell -k +1}}{(t-s)(1-tr)(t-r)},\\
        g_{22}(2,k;1,\ell) &= \frac{r^{\ell- k}}{(r-s)(1-sr)(1-tr)} + \frac{s^{\ell - k}}{(1-st)(1-sr)(s-r)},\\
        g_{22}(1,k;1,\ell) &= \frac{s^{\ell-k}-r^{\ell-k}}{(s-r)(1-sr)}.
    \end{aligned}
\end{equation}

\begin{lem}
The $22-$entry can be reformulated as
\begin{equation}
    \begin{aligned}
        K_{22}(n,k; n^{\prime}, \ell) =h_{22}(n,k;n',\ell)+ f_{22}(n,k;n',\ell)+ g_{22}(n,k;n',\ell)\Id_{k \leq \ell} - g_{22}(n',\ell;n,k)\Id_{k > \ell},
    \end{aligned}
\end{equation}
where 
\begin{equation}
    \begin{aligned}
        &f_{22}(2,k;2,\ell) = \frac{s^{\ell +1}t^{k+1} - t^{\ell+1}s^{k+1}}{(t-s)(t-r)(s-r)(1-st)H(s)H(t)},\\
        &f_{22}(1,k;2,\ell) = \frac{s^k t^{\ell+1}}{H(s)H(t)(1-st)(s-r)(t-r)},\\
        &f_{22}(2,k;1,\ell) = -f_{22}(1,\ell; 2,k)  = -\frac{s^{\ell} t^{k+1}}{H(s)H(t)(1-st)(s-r)(t-r)},\\
        &f_{22}(1,k;1,\ell) = 0.
    \end{aligned}
\end{equation}
\end{lem}

\begin{proof}
When \(n=n'=1\), the kernel \(K_{22}(1,k;1,\ell)\) is independent of \(t\). 
Indeed, the pole \(z=1/t\) appears only when \(n=2\), through the factor
\(B(2,z)=(1-sz)(1-tz)\), while the pole \(w=t\) appears only when
\(n'=2\), through the factor \(A(2,w)=(1-s/w)(1-t/w)\).

We now compare the poles appearing in \(K_{22}\) with those already included in
\(h_{22}\). The definition of \(h_{22}\) contains all relevant pairs of poles
except the four pairs
\[
    (z,w)=(1/s,s),\qquad (1/t,t),\qquad (1/t,s),\qquad (1/s,t),
\]
whenever these poles are present. The first two pairs do not contribute because of
the numerator \(zw-1\).
The only nonzero omitted contributions come from the two pairs
\[
    (z,w)=(1/t,s),\qquad (z,w)=(1/s,t).
\]
They are special because their
residue evaluations contain the factor \(1/(1-st)\), which becomes singular as
\(t\to 1/s\). We evaluate these contributions case by case.

When \(n=n'=2\), both off-diagonal pairs are present. Hence
\begin{equation}
    \begin{aligned}
        &-\oint_{\Gamma_s}\frac{dw}{2\pi\I}
          \oint_{\Gamma_{1/t}}\frac{dz}{2\pi\I}
          F(2,k;2,\ell)
        -\oint_{\Gamma_t}\frac{dw}{2\pi\I}
          \oint_{\Gamma_{1/s}}\frac{dz}{2\pi\I}
          F(2,k;2,\ell)
        = f_{22}(2,k;2,\ell).
    \end{aligned}
\end{equation}

When \(n=1\) and \(n'=2\), the pole \(z=1/t\) is absent, while the pole
\(w=t\) is present. Therefore only the pair
\((z,w)=(1/s,t)\) contributes and $f_{22}(1,k;2,\ell)$.

When \(n=2\) and \(n'=1\), the pole \(z=1/t\) is present, while the pole
\(w=t\) is absent. Therefore only the pair
\((z,w)=(1/t,s)\) contributes $f_{22}(2,k;1,\ell)$.

Finally, when \(n=n'=1\), neither \(z=1/t\) nor \(w=t\) is present.
Therefore $f_{22}(1,k;1,\ell)=0.$

It remains to evaluate the residue term \(g_{22}(n,k;n',\ell)\). Assume first
that \(k\leq \ell\). When $n = n' = 2$, there are poles at $z = s,t,r$. When $n =1, n' = 2$, there are poles at $z = s,t,r$. When $n =2, n' = 1$, there are poles at $z = s,r$. Direct evaluations of $g_{22}(n,k;n',\ell)$ is provided in \eqref{eq: residue evaluation}. By anti-symmetry, we get the $g_{22}$ for $k> \ell.$
\end{proof}

\section{Limit to the two-parameter stationary case in high-density regime}\label{paper:stationary-limit}
We define
\begin{equation}
    \begin{aligned}
        &\widehat{A}(1,x) : = (1-s/x), \quad \widehat{A}(2,x) := (1-s/x)(1-1/sx),\\
        &\widehat{B}(1,x) := (1-sx),\quad \widehat{B}(2,x) := (1-sx)(1-x/s)
    \end{aligned}
\end{equation}
as limits of $A(i,x)$ and $B(i,x)$ when $t\rightarrow 1/s$ for $i\in \{1,2\}.$

We define some functions for the limit of the kernel. For $n,n'\in\{1,2\},$  \begin{equation}\label{eq: stat K11}
\begin{aligned}
    \widehat{K}_{11}(n,k;n',\ell)&:= -\oint\limits_{\Gamma_{\sqrt{q}}} \frac{dw}{2\pi\I} \oint\limits_{\Gamma_{1/\sqrt{q}}} \frac{dz}{2\pi\I}  \frac{w^{\ell-1}}{z^{k}}\frac{H(z)\widehat{A}(n,z)\widehat{B}(n',w)}{H(w)}\frac{(zw-1)(z-r)(1-wr)}{(z^2-1)(1-w^2)(z-w)},\\
\end{aligned}
\end{equation}

\begin{equation}\label{eq: stat K12}
\begin{aligned}
    \widehat{K}_{12}(n,k; n', \ell) &:= -\oint\limits_{\Gamma_{r,s,1/s,\sqrt{q}}} \frac{dw}{2\pi\I} \oint\limits_{\Gamma_{1/\sqrt{q}}} \frac{dz}{2\pi\I}\frac{w^{\ell-1}}{z^{k}} \frac{H(z)\widehat{A}(n,z)}{H(w)\widehat{A}(n',w)}\frac{(zw-1)(z-r)}{(z^2-1)(z-w)(w-r)}\\
    &- \Id_{k\leq \ell}\Id_{n = 1, n' = 2}s^{k-\ell},\\
\end{aligned}
\end{equation}

\begin{lem}\label{lem: K11 and K12}
For any $\sqrt{q} \in (0,1),$ $s\in (\sqrt{q},1),$ $r \in (0,s)\cup(s,1/s)$, and $\epsilon \ll 1,$ $K_{11}$ and $K_{12}$ are analytic for $t\in (1/s -\epsilon, 1/s+\epsilon)$ and their $t\rightarrow 1/s$ limits are as follows:
\begin{equation}
    \begin{aligned}
        &\lim_{t\rightarrow 1/s} K_{11}(n,k; n', \ell) = \widehat{K}_{11}(n,k;n',\ell),\quad
        \lim_{t\rightarrow 1/s} K_{12}(n,k; n', \ell) = \widehat{K}_{12}(n,k;n',\ell).
    \end{aligned}
\end{equation}
Moreover, there exists $C>0$ independent of $k,\ell$ such that for $t \in (1/s-\epsilon, 1/s+\epsilon)$
\begin{equation}
\begin{aligned}
    &|K_{11}(n,k;n',\ell)| \leq Ce^{(k+\ell)\log(\sqrt{q}+\epsilon)},\\
    &|K_{12}(n,k;n',\ell)| \leq Ce^{k\log(\sqrt{q}+\epsilon)+\ell\log(\max\{s,r,1/s\} + \epsilon)} + e^{(\ell-k)\log(1/s+\epsilon)}\Id_{k\leq \ell,n=1,n'=2}.
\end{aligned}
\end{equation}
\end{lem}

\begin{proof}
    For $K_{12}$, we can deform the $w-$contour $\Gamma_{t}$ to $\Gamma_{1/s}$. Since the contour for $z$ can be chosen to lie arbitrarily close to $1/\sqrt{q}$, the deformation of $w-$contour does not affect the $z-$contour. For $K_{11}$, we can directly take the $t\rightarrow 1/s$ limit because it does not affect any poles or contours.
    The upper bound for $K_{11}$ is obtained by taking the $z-$contour to be $|z-1/\sqrt{q}| = \epsilon/2$ and the $w-$contour to be $|w-\sqrt{q}| = \epsilon/2$.
    The upper bound for the double integrals of $K_{12}$ is given by two pairs of poles $(z,w) = (1/\sqrt{q},1/s),$ $(z,w) = (1/\sqrt{q},r)$ and by setting the $z-$contour to be $|z-1/\sqrt{q}| = \epsilon/2$.
\end{proof}

The $22$-entry requires exact evaluation since it contains terms with $(1-st)^{-1}.$ We define
\begin{equation}
    \begin{aligned}
        &\widehat{R}(2,k;2,\ell)= \frac{s^2r^{\ell - k +1}}{(s-r)^2(1-sr)^2}
        +\frac{s^{2-k+\ell}}{(s-r)(1-sr)^2(1-s^2)} - \frac{rs^{2+k-\ell}}{(s-r)^2(1-sr)(1-s^2)}\\
        &+ \frac{ks^{2+k-\ell} - (1+k)s^{2-k+\ell}}{(s-r)(1-sr)(1-s^2)}
        + \frac{N\sqrt{q} s^{1-k+\ell} - N \sqrt{q}s^{1+k-\ell}}{(1-\sqrt{q}/s)(s-r)(1-rs)(1-s^2)}
        + \frac{N\sqrt{q} s^{3-k+\ell} - N \sqrt{q}s^{3+k-\ell}}{(1-\sqrt{q}s)(s-r)(1-rs)(1-s^2)},\\
        &\widehat{R}(1,k;2,\ell)
        = \frac{(1-s^2)sr^{\ell - k + 1} - (1-sr)s^{\ell - k + 2}}{(1-sr)^2(1-s^2)(s-r)}
        + \frac{rs^{k-\ell}}{(s-r)^2(1-sr)} - \frac{ks^{k-\ell}}{(s-r)(1-sr)}\\
        &+ \frac{s^{k-\ell}}{(s-r)(1-sr)(1-s^2)}+\frac{N\sqrt{q}s^{k-\ell-1}}{(1-\sqrt{q}/s)(s-r)(1-sr)} + \frac{N\sqrt{q}s^{k-\ell+1}}{(s-r)(1-\sqrt{q}s)(1-rs)},\\
        &\widehat{R}(2,k;1,\ell) = - \frac{sr^{\ell-k}}{(s-r)^2(1-sr)}
        -\frac{s^{\ell - k}}{(s-r)(1-sr)^2} + \frac{(1+k)s^{\ell - k}}{(s-r)(1-sr)} - \frac{N\sqrt{q}s^{\ell - k -1}}{(1-\sqrt{q}/s)(s-r)(1-rs)}\\
        &- \frac{N\sqrt{q}s^{\ell - k + 1}}{(1-\sqrt{q}s)(s-r)(1-rs)},\\
        &\widehat{R}(1,k;1,\ell)
        = \frac{s^{\ell-k}-r^{\ell-k}}{(s-r)(1-sr)}.
        \end{aligned}
    \end{equation}

\begin{equation}
    \begin{aligned}
        \widehat{{C}}(n,k;n',\ell) &:= \frac{r^{\ell}s^k\!\left(\Id_{n'=1}+\frac{sr}{sr-1}\Id_{n'=2}\right)\!\!\left(\Id_{n=1}+\frac{s^2}{s^2-1}\Id_{n=2}\right)}{H(s)H(r)(s-r)(1-sr)}\\
        &\qquad{}+ \frac{r^{k}s^{\ell}\!\left(\Id_{n'=1}+ \frac{s^2}{s^2-1}\Id_{n'=2}\right)\!\!\left(\Id_{n=1}+ \frac{sr}{sr-1}\Id_{n=2}\right)}{H(s)H(r)(1-sr)(r-s)},\\
        \widehat{{D}}(n,k;n',\ell) &:= \Id_{n=2}\frac{r^{\ell}H(s)}{s^{k-1}H(r)}\frac{\left(\Id_{n'=1}+ \frac{sr}{sr-1}\Id_{n'=2} \right)}{(s-r)^2(1-s^2)}\\
        &\qquad{}- \Id_{n'=2} \frac{r^kH(s)}{s^{\ell-1}H(r)}\frac{\left(\Id_{n=1} + \frac{sr}{sr-1}\Id_{n=2}\right)}{(s-r)^2(1-s^2)}.
    \end{aligned}
\end{equation}

\begin{equation}
\begin{aligned}
    \widehat{h}_{22}(n,k;n',\ell) &:= -\!\!\oint\limits_{\Gamma_{\sqrt{q},r,s,1/s}} \!\!\frac{dw}{2\pi\I} \!\oint\limits_{\Gamma_{1/\sqrt{q}}} \!\!\frac{dz}{2\pi\I}\widehat{F}(z,w)-\oint\limits_{\Gamma_{\sqrt{q}}} \frac{dw}{2\pi\I} \!\!\!\oint\limits_{\Gamma_{1/r,1/s,s}} \!\!\!\frac{dz}{2\pi\I}\widehat{F}(z,w),\\
\end{aligned}
\end{equation}
where 
\begin{equation}\label{eq: 22 integrand}
    \widehat{F}(z,w) :=\frac{w^{\ell-1}}{z^{k}} \frac{H(z)}{H(w)\widehat{A}(n',w)\widehat{B}(n,z)}\frac{(zw-1)}{(z-w)(1-zr)(w-r)}.
\end{equation}
Given functions introduced above, we define
\begin{equation}\label{eq: stat K22}
    \widehat{K}_{22}(n,k;n',\ell) := \widehat{h}_{22}(n,k;n',\ell) + \widehat{C}(n,k;n',\ell) + \widehat{D}(n,k;n',\ell) + \widehat{R}(n,k;n',\ell)\Id_{k\leq \ell} - \widehat{R}(n',\ell;n,k)\Id_{k>\ell}
\end{equation}

\begin{lem}\label{lem: K22}
    Fix any $s\in (\sqrt{q},1),$ and $r \in (0,s)\cup(s,1/s)$, $\epsilon \ll 1$.There exists a constant $C>0$ independent of $k,\ell$ such that for $t\in (1/s-\epsilon, 1/s+\epsilon)$
    \begin{equation}
        |K_{22}(n,k;n',\ell)| \leq Ce^{(k+\ell)\log(1/s+\epsilon)}.
    \end{equation}
    The analytic continuation of $K_{22}$ under the limit $t\rightarrow 1/s$ is
    \begin{equation}
        \lim_{t\rightarrow 1/s} K_{22}(n,k;n',\ell) =\widehat{K}_{22}(n,k;n',\ell).
    \end{equation}
\end{lem}

\begin{proof}
The convergence of $h_{22}(n,k;n',\ell)$ to $\widehat{h}_{22}+\widehat{C} + \widehat{D}$ is analogous to that of $K_{11}$ and $K_{12}$. We deform the $w-$contour $\Gamma_{t}$ to $\Gamma_{1/s}$ when the pole $w = t$ is present and deform the $z-$contour $\Gamma_{1/t}$ to $\Gamma_{s}$ when the pole $z = 1/t$ is present. These deformations do not cause the $z-$
and $w-$contours to cross because $s\in (\sqrt{q},1)$ and $r\neq s,1/s$. It is clear that
\begin{equation}
    \begin{aligned}
        &\lim_{t\rightarrow 1/s}  -\left(\oint\limits_{\Gamma_{r}}\frac{dw}{2\pi\I}\!\! \oint\limits_{\Gamma_{1/s}}\!\! \frac{dz}{2\pi\I} +\oint\limits_{\Gamma_{s}}\frac{dw}{2\pi\I}\!\! \oint\limits_{\Gamma_{1/r}}\!\! \frac{dz}{2\pi\I}\right) {F}(z,w) = \widehat{C}(n,k;n',\ell),\\
        &\lim_{t\rightarrow 1/s} -\left(\oint\limits_{\Gamma_{r}}\frac{dw}{2\pi\I} \oint\limits_{\Gamma_{s}} \frac{dz}{2\pi\I} +\oint\limits_{\Gamma_{1/s}}\frac{dw}{2\pi\I}\!\! \oint\limits_{\Gamma_{1/r}}\!\! \frac{dz}{2\pi\I}\right) F(z,w) = \widehat{D}(n,k;n',\ell).
    \end{aligned}
\end{equation}
A direct evaluation using L'Hospital's rule yields
\[
    \lim_{t\to 1/s}
        f_{22}(n,k;n',\ell)+g_{22}(n,k;n',\ell)\Id_{k\leq \ell} = \widehat R(n,k;n',\ell).
\]
By the anti-symmetry of the kernel, the case $k>\ell$ follows immediately. The upper bound for $K_{22}$ arises from the pairs of poles $(z,w) = (1/s,t)$, $(z,w) = (1/r,t)$, $(z,w) = (1/s,r)$, together with the fact that $r < 1/s$.
\end{proof}

\subsection{Convergence of Fredholm Pfaffian}

\begin{prop}\label{Fredconv}
    For any fixed $\sqrt{q}\in (0,1),$ $s\in (\sqrt{q},1)$, $r\in (0,s)\cup(s,1/s)$, $y \in \mathbb{Z}_{\geq 0}$, and $0<\epsilon \ll \min\left({\sqrt{q}}^{-1}- s^{-1}, \frac{1-\sqrt{q}/s}{\sqrt{q}+\sqrt{q}^{-1}}\right)$, $\mathrm{Pf}(J - {K})$ is analytic for $t \in (1/s-\epsilon,1/s+\epsilon)$ and it has the limit
    \begin{equation}
        \lim_{t\rightarrow 1/s} \mathrm{Pf}(J - {K})_{\ell^2(T)} = \mathrm{Pf}(J - \widehat{K})_{\ell^2(T)},
    \end{equation}
    where $T =\{1\}\times \{d+1,d+2, \dots \}\cup \{2\}\times \{d+y+1,d+y+2, \dots \}$.
\end{prop}

\begin{proof}
Let $\eta = \epsilon + 1/s$. By lemmas \ref{lem: K11 and K12} and \ref{lem: K22}, we have the following upper bound for the $nn'-$block of the kernel:
\begin{equation}
    \begin{aligned}
        |K(n,k;n',\ell)| &=\Bigg|\begin{pmatrix}
        K_{11}(n,k;n',\ell) & K_{12}^{reg}(n,k;n',\ell)\\
        K_{21}^{reg}(n,k;n',\ell) & K_{22}(n,k;n',\ell)
    \end{pmatrix}\Bigg|\\
        &\leq \begin{pmatrix}
        Ce^{(k+\ell)\log (\sqrt{q}+\epsilon)} & Ce^{k\log(\sqrt{q}+\epsilon) + \ell \log (\eta)}\\
        Ce^{k \log(\eta) + \ell\log(\sqrt{q}+ \epsilon)} & Ce^{(k+\ell) \log (\eta)}
    \end{pmatrix}
    \end{aligned}
\end{equation}
for some $C>0$ independent of $k,\ell$.

The special term in $K_{12}$ for $n,n' \in \{1,2\}$ becomes
\begin{equation}
\begin{aligned}
    &V(1,k;1,\ell) = V(2,k;2,\ell) = \begin{pmatrix}
        0 & 0 \\
        0 & 0
    \end{pmatrix},\\
    &V(1,k;2,\ell) = \begin{pmatrix}
        0 & - t^{\ell - k}\Id_{k\leq \ell} \\
        0 & 0
    \end{pmatrix}, \quad V(2,k;1,\ell) = \begin{pmatrix}
        0 & 0 \\
        t^{k - \ell}\Id_{\ell \leq k} & 0
    \end{pmatrix},\\
\end{aligned}
\end{equation}
and we have the upper bound
\begin{equation}
        |V(n,k;n',\ell)| \leq \begin{pmatrix}
            0 & Ce^{(\ell-k)\log(\eta)}\Id_{k\leq \ell}\Id_{n<n'}\\
            Ce^{(k-\ell)\log(\eta)}\Id_{\ell\leq k}\Id_{n>n'} & 0
        \end{pmatrix}.
\end{equation}
Consider the conjugation matrix 
\begin{equation}
        M(n,k)=
        \text{diag}(e^{u_nk}, e^{-u_n k}).
\end{equation}
We define the conjugated kernel to be 
\begin{equation}
    MKM(n,k;n',\ell) = 
    \begin{pmatrix}
        e^{u_{n}k} &0\\
        0 & e^{-u_{n} k}
    \end{pmatrix}K(n,k;n',\ell)\begin{pmatrix}
        e^{u_{n'}\ell} &0\\
        0 & e^{-u_{n'} \ell}
    \end{pmatrix}.
\end{equation}
It is clear that under this gauge transform we have $\Pf(J-K)_{\ell^2(T)} = \Pf(J - MKM)_{\ell^2(T)}$.
We also have the upper bound
\begin{equation}\label{eq: upper bounds for K}
    |(MKM)(n,k;n',\ell)| \leq \begin{pmatrix}
            Ce^{k(\log(\sqrt{q}+\epsilon) + u_n) + \ell(\log(\sqrt{q}+\epsilon) + u_{n'})} & Ce^{k(\log(\sqrt{q}+\epsilon) + u_n) + \ell(\log(\eta) - u_{n'})}\\
            Ce^{k(\log(\eta)-u_n) + \ell(\log(\sqrt{q}+\epsilon) + u_{n'})} & Ce^{k(\log(\eta) - u_n)) + \ell (\log(\eta) - u_{n'})}
        \end{pmatrix},
\end{equation}
and 
\begin{equation}\label{eq: upper bounds for V}
\begin{aligned}
    &|MVM(n,k;n',\ell)| \leq \begin{pmatrix}
            0 & A\\
            B & 0
    \end{pmatrix},
\end{aligned}
\end{equation}
where the $12$-entry is bounded by
\begin{equation}
    A = Ce^{k(u_n-(u_n+u_{n'})/2) + \ell((u_n + u_{n'})/2 - u_{n'}) + (\ell - k)(\log(\eta) - (u_n + u_{n'})/2)}\Id_{k\leq \ell} \Id_{n<n'},
\end{equation}
and the $21$-entry is bounded by
\begin{equation}
    B = Ce^{\ell(u_{n'}-(u_n+u_{n'})/2) + k((u_n + u_{n'})/2 - u_{n}) + (k - \ell)(\log(\eta) - (u_n + u_{n'})/2)}\Id_{\ell\leq k}\Id_{n>n'}.
\end{equation}
Therefore, we need
\begin{equation}
    \begin{aligned}
        &\log(\eta) <u_n, u_{n'} < -\log(\sqrt{q}+\epsilon),\\
    \end{aligned}
\end{equation}
and if $n < n',$ we need in addition need
\begin{equation}
    \begin{aligned}
        u_n<u_{n'}.
    \end{aligned}
\end{equation}
Combining all these restrictions together, it suffices to choose $u_n$, $u_{n'}$ satisfying
\begin{equation}\label{eq: range for conjugation}
\begin{aligned}
    &\log(\eta)<u_n < u_{n'}< -\log({\sqrt{q}+\epsilon}) \text{ for } n<n',\quad
    \log(\eta)<u_{n'} < u_{n} < -\log(\sqrt{q}+\epsilon)\text{ for } n'<n.\\
\end{aligned}
\end{equation}
Then condition on $\epsilon$ guaranties that $\log(\eta) < -\log(\sqrt{q}+\epsilon)$.
Using upper bounds \eqref{eq: upper bounds for K} and \eqref{eq: upper bounds for V} and the point-wise convergence result from Lemmas \ref{lem: K11 and K12} and \ref{lem: K22}, we get the desired convergence of Fredholm Pfaffian by dominated convergence theorem and Hadamard's bound, see \cite[Lemma 5.1]{zeng}. 
\end{proof}

The following Corollary is an immediate consequence of Proposition \ref{prop:maxkernel} and Proposition \ref{Fredconv}. We state it without proof.
\begin{cor}\label{eq: finite N high}
    For any fixed $\sqrt{q} \in (0,1),$ $s\in (\sqrt{q},1),$ $r\in (0,s)\cup (s,1/s),$ we have
    \begin{equation}\label{eq: high stat prelimit}
        \begin{aligned}
            \Pb\left( G_{r,s}^{stat}(N+2,N+2) \leq d \right) = \sum_{y = 0}^{\infty} \Pf\left(J - \widehat{K}\right)_{\ell^{2}(T)}(1-sr)(sr)^y,
        \end{aligned}
    \end{equation}
    where $\widehat{K}$ is defined by \eqref{eq: stat K11}, \eqref{eq: stat K12}, \eqref{eq: stat K22} and $T = \{1\}\times \{d+1,d+2,\dots\} \cup \{2\} \times\{d+y+1,\dots\}$.
\end{cor}

\begin{rmk}
    A key simplification is that our formula \eqref{eq: high stat prelimit} avoids kernel resolvents, which exists in \cite{zeng,zeng2026} and \cite[Theorem 2.7]{BeteaFerrariOccelli2020}. In previous work, these resolvent terms can be replaced by rewriting the expressions as differences of two Fredholm Pfaffians, but the resulting formulas remain complicated. Our approach does not need to extract divergent residues from the kernel. This is possible because, in the flipped models \eqref{eq: weigths flipped} and \eqref{eq: product stationary weigths flipped}, $K_{11}$ has poles only at $z=1/\sqrt{q}$ and $w=\sqrt{q}$, whose fast decay controls the convergence of the Fredholm Pfaffian.
\end{rmk}

\section{Asymptotics of the kernel under high-density regime}\label{paper:hd-asymptotics}
In this section, we use the previously obtained finite $N$ distribution formulas, apply the critical scaling to discover the universal limiting distribution for half-space KPZ fixed point.
We define a constant 
\begin{equation}
    c_0 = \frac{2\sqrt{q}(1+\sqrt{q})}{(1-\sqrt{q})^3}.
\end{equation} 
Throughout this paper, we let $\rho = (c_0N/2)^{1/3}$
Consider the following scaling of parameters.
\begin{equation}\label{eq: scale}
    \begin{aligned}
        &k = \frac{2\sqrt{q}}{1-\sqrt{q}}N + \rho X, \quad \ell = \frac{2\sqrt{q}}{1-\sqrt{q}}N + \rho Y,\quad d= \frac{2\sqrt{q}}{1-\sqrt{q}}N + \rho\tilde{d},&\quad y = \rho\tilde{y}\\
&s = 1 + \rho^{-1}\tilde{s},\quad r = 1 + \rho^{-1}\tilde{r},\quad
z = 1 + \rho^{-1}\zeta,\quad w = 1+ \rho^{-1}\omega.\\
    \end{aligned}
\end{equation}
The condition $s\in (\sqrt{q},1)$ and $r\in (0,s)\cup(s,1/s)$ translate into the condition $\s<0$ and $\rr\in (-\infty, \s)\cup(\s,-\s)$.
It is clear that
\begin{equation}
    \begin{aligned}
        &\lim_{N\rightarrow \infty} \rho\widehat{A}(1,z) = \A(1,\zeta), \quad \lim_{N\rightarrow \infty} \rho^2\widehat{A}(2,z) = \A(2,\zeta),\\
        &\lim_{N\rightarrow \infty} \rho \widehat{B}(1,z) = \B(1,\zeta), \quad \lim_{N\rightarrow \infty} \rho^2\widehat{B}(2,z) = \B(2,\zeta).
    \end{aligned}
\end{equation}

% \textcolor{red}{to be removed}
% The scaling of each entry is as follows:
% \begin{equation}
%     \begin{aligned}
%         &((c_0N/2)^{1/3})^{n+n'+2}K_{11}(n,k;n',\ell) , \quad ((c_0N/2)^{1/3})^{n-n'+1}K_{12}(n,k;n',\ell),\\
%         &((c_0N/2)^{1/3})^{n'-n+1}K_{21}(n,k;n',\ell), \quad ((c_0N/2)^{1/3})^{-n-n'}K_{22}(n,k;n',\ell).
%     \end{aligned}
% \end{equation}

\begin{lem}\label{lem: K11,K12 steepest}
    For any given $u>0$, the following limit holds uniformly over $X,Y \in [-u,u]$,
    \begin{equation}
        \begin{aligned}
            &\lim_{N\rightarrow \infty}\rho^{n+n'+2}\widehat{K}_{11}(n,k;n',\ell) = \widetilde{K}_{11}(n,X;n',Y),\quad \lim_{N\rightarrow \infty}\rho^{n-n'+1}\widehat{K}_{12}(n,k;n',\ell) = \widetilde{K}_{12}(n,X;n',Y).\\
        \end{aligned}
    \end{equation}
    Moreover, for any $\delta>-\s,$ there exist constants $C>0$ independent of $X,Y$, $N_0 \in \mathbb{N}$, and $\epsilon>0$ satisfying $0<\epsilon -\s <\delta$ such that for any $X,Y \geq -u$ and $N\geq N_0,$ we have
    \begin{equation}
        \begin{aligned}
            &|\rho^{n+n'+2}\widehat{K}_{11}(n,k;n',\ell)| \leq Ce^{-\delta(X+Y)},\\
            &|\rho^{n-n'+1}\widehat{K}_{12}(n,k;n',\ell)| \leq Ce^{-\delta X -(\tilde{s}-\epsilon)Y} + \Id_{X\leq Y}\Id_{n=1,n'=2}Ce^{(-\s+\epsilon)(Y-X)}.
        \end{aligned}
    \end{equation}
\end{lem}
\begin{proof}
We perform the standard steepest descent method. We provide a sketch and the remaining details of the proof can be found in \cite[Section 6.1]{zeng}.
Let \begin{equation}
    \begin{aligned}
        h_0(x) &= \log(1-\sqrt{q}/x) - \log (1-\sqrt{q}x) - \frac{2\sqrt{q}}{1-\sqrt{q}}\log (x).\\
    \end{aligned}
\end{equation} We choose steepest descent path for $w-$contour, denoted as $\mathcal{C}_1$, as in \cite[(6.23)]{zeng} and the $z-$contour, denoted as $\mathcal{C}_2$, as in \cite[(6.26)]{zeng} with the modification of a $\mathcal{O}(\rho^{-1})$ neighborhood around the critical point $z = w =1$ as explained in \cite[Figure 6]{zeng}. Applying the scaling and the change of variables yield
\begin{equation}
    \begin{aligned}
        &\rho^{n-n'+1}\widehat{K}_{12}(n,k;n',\ell) = -\rho^{n-n'+1}\int\limits_{\mathcal{C}_1}\frac{dw}{2\pi\I}\int\limits_{\mathcal{C}_2} \frac{dz}{2\pi\I} e^{N(h_0(z) - h_0(w)) - \rho(X\log(z) - Y\log(w))}\\
        &\quad\quad\quad\quad\frac{\widehat{A}(n,z)(zw-1)(z-r)}{\widehat{A}(n',w)w(z^2-1)(z-w)(w-r)} - \Id_{\{X\leq Y\}}\Id_{\{n = 1,n' = 2\}}e^{\s(X-Y) + \mathcal{O}(\rho^{-1}(X-Y))}\\
        =& \,\mathcal{O}(e^{-\epsilon N}) -\!\!\!\!\!\int\limits_{{}_{\tilde{r},-\s,\s}\wcu\, {}_{}}\!\!\!\! \frac{d\omega}{2\pi\I}\!\int\limits_{{ }_{0, \omega}\zcd } \!\frac{d\zeta}{2\pi\I} e^{\frac{\zeta^3}{3} - \frac{\omega^3}{3}  - X\zeta + Y\omega +\mathcal{O}\left(\rho^{-1}\right)} \left(\frac{\A(n,\zeta)}{\A(n',\omega)}\frac{(\zeta + \omega)(\zeta - \tilde{r})}{(2\zeta)(\omega - \tilde{r})(\zeta - \omega)}+\mathcal{O}\left(\rho^{-1}\right)\right)\\
        &- \Id_{X\leq Y} \Id_{n=1,n'=2}e^{\s(X-Y) + \mathcal{O}(\rho^{-1}(X-Y))},
    \end{aligned}
\end{equation}
for some $\epsilon >0$. The  term $\mathcal{O}(e^{-\epsilon N})$ comes from the decay of integration over contours that stay away from the critical point, i.e., $\mathcal{C}_1\times \mathcal{C}_2 \cap (B(1,\eta))^c$ for some fixed $\eta >0.$ The integration over $\mathcal{C}_1\times \mathcal{C}_2 \cap B(1,\eta)$ contributes to the double integral in the asymptotic. 

Next, we show the upper bound of $\widehat{K}_{12}$. Since the $z-$contour only needs to enclose the pole $z=1/\sqrt{q}$, we may deform the $z-$contour to pass through any point $\delta>\max\{\rr,\s,-\s\}$ provided that it does not cross the $w-$contour. In contrast, the $w-$contour must enclose poles $w = r,s,1/s,\sqrt{q}$. Under the scaling, these poles become $\omega = \rr,\s,-\s$, while the pole at $\sqrt{q}$ is sent to $-\infty.$ Therefore, the pole at $\omega = -\s$ provides the dominant contribution to the upper bound. The error term $e^{\mathcal{O}(\rho^{-1}(X-Y))}$ can be absorbed into $\epsilon$ after choosing a sufficiently large $N_0$. The asymptotic limit and upper bound for $\rho^{n+n'+2}\widehat{K}_{11}$ can be obtained analogously.
\end{proof}

\begin{lem}\label{lem: limit and bound of K22}
    For any $n,n' \in \{1,2\}$ and any given $u > 0$, the following limit holds uniformly over $X,Y \in [-u, u]$,
    \begin{equation}
        \begin{aligned}
            &\lim_{N\rightarrow \infty}\rho^{-n-n'}\widehat{K}_{22}(n,k;n',\ell) = \widetilde{h}_{22}(n,X;n',Y) + \widetilde{R}_{22}(n,X;n',Y)\Id_{X\leq Y} - \widetilde{R}_{22}(n',Y;n,X)\Id_{X> Y}.
        \end{aligned}
    \end{equation}
    Moreover, there exist a constant $C>0$ independent of $X,Y$, $N_0 \in \mathbb{N}$, and $\epsilon>0$ such that for any $X,Y \geq -u$ and $N\geq N_0$, we have
    \begin{equation}
        |\rho^{-n-n'}K_{22}(n,k;n',\ell)|\leq Ce^{-(\tilde{s} -\epsilon)(X+Y)}.
    \end{equation}
\end{lem}

\begin{proof}
    Recall the definition of $\widehat{F}$ in \eqref{eq: 22 integrand} and the definition of $\widetilde{F}$ in \eqref{eq: integrand of 22 limit}. Using the steepest descent method explained in Lemma \ref{lem: K11,K12 steepest}, we get the following double integral limits
    \begin{equation}
        \begin{aligned}
        &\lim_{N\rightarrow \infty}-\rho^{-n-n'}\oint\limits_{\Gamma_{\sqrt{q}}}\frac{dw}{2\pi\I}\!\!\!\!\! \oint\limits_{\Gamma_{1/\sqrt{q},1/r,1/s,s}}\!\!\!\!\! \frac{dz}{2\pi\I} \widehat{F}(z,w) = \int\limits_{{}_{}\wcu\, {}_{\tilde{r},-\s,\zeta,0}}\! \frac{d\omega}{2\pi\I}\!\int\limits_{{ }_{}\zcd {}_{\s,-\tilde{r}, -\tilde{s}} } \!\frac{d\zeta}{2\pi\I} \widetilde{F}(\zeta,\omega),\\
        &\lim_{N\rightarrow \infty}-\rho^{-n-n'}\!\!\!\!\!\oint\limits_{\Gamma_{r,s,1/s}}\!\!\!\!\!\frac{dw}{2\pi\I}\!\!\oint\limits_{\Gamma_{1/\sqrt{q}}} \frac{dz}{2\pi\I} \widehat{F}(z,w) = \int\limits_{\Gamma_{\rr,\s,-\s}}\! \frac{d\omega}{2\pi\I}\!\int\limits_{{ }_{\s,-\tilde{r}, -\tilde{s},\omega}\zcd {}_{} } \!\frac{d\zeta}{2\pi\I}  \widetilde{F}(\zeta,\omega).
        \end{aligned}
    \end{equation}
    We need exact evaluation of all other poles. For two pairs of poles $(z,w) = (1/s,r)$ and $(z,w) = (1/r,s),$ we have
    \begin{equation}
    \begin{aligned}
        &\lim_{N\rightarrow \infty}\rho^{-n-n'}\widehat{C}(n,k;n',\ell)= \widetilde{G}(n,X;n',Y),
    \end{aligned}
\end{equation}
where $\widetilde{G}$ is defined in \eqref{eq: definitino of G tilde}.

For two pairs of poles $(z,w) = (s,r)$ and $(z,w) = (1/r,1/s),$ we have
\begin{equation}
    \begin{aligned}
        &\lim_{N\rightarrow \infty}\rho^{-n-n'}\widehat{D}(n,k;n',\ell)= \widetilde{O}(n,X;n',Y)
    \end{aligned}
\end{equation}
where $\widetilde{O}$ is defined in \eqref{eq: def of O}.
As a result, we get the $$\lim_{N\rightarrow \infty} \rho^{-n-n'}\left(\widehat{h}_{22} + \widehat{C} + \widehat{D}\right)(n,k;n',\ell) = \widetilde{h}_{22}(n,X;n',Y).$$
Next, we discuss the asymptotics of $\widehat{R}$. We first explain the most difficult term, i.e., $\widehat{R}(2,k;2,\ell)$. 
The following three terms admit direct asymptotic limits:
\begin{equation}
    \begin{aligned}
        & \frac{\rho^{-4} s^2r^{\ell-k+1}}{(s-r)^2(1-sr)^2} \rightarrow \frac{e^{\rr(Y-X)}}{(\s - \rr)^2(\s+\rr)^2},\quad 
        \frac{\rho^{-4}s^{2-k+\ell}}{(s-r)(1-sr)^2(1-s^2)} \rightarrow \frac{e^{\s(Y - X)}}{(\s - \rr)(\s+\rr)^2(-2\s)},\\
        & \frac{-\rho^{-4} rs^{2+k-\ell}}{(s-r)^2(1-sr)(1-s^2)} \rightarrow -\frac{e^{\s(X-Y)}}{(\s-\rr)^2(\s+\rr)(2\s)}.
    \end{aligned}
\end{equation}

Next, we claim that
\begin{equation}\label{eq: n=2,n'=2}
    \begin{aligned}
        \lim_{N\rightarrow \infty}\frac{\rho^{-4} s^{1+k-\ell}}{(s-r)(1-sr)(1-s^2)}\left(ks - \frac{N\sqrt{q}}{(1-\sqrt{q}/s)} - \frac{N\sqrt{q}}{(1-\sqrt{q}s)}s^2\right) = \frac{e^{\s(X-Y)}}{(\s-\rr)(\s+\rr)(2\s)}\left(X -\s^2\right).
    \end{aligned}
\end{equation}
To prove this, we rewrite LHS of \eqref{eq: n=2,n'=2} as $ J(s) ( I(s) +Xs),$
where
\begin{equation}
    \begin{aligned}
        I(s) := \rho^{-1}\left(\frac{2\sqrt{q}N}{1-\sqrt{q}}s - \frac{N\sqrt{q}}{(1-\sqrt{q}/s)} - \frac{N\sqrt{q}}{(1-\sqrt{q}s)}s^2\right),\quad
        J(s) := \frac{\rho^{-3}s^{1+k-\ell}}{(s-r)(1-sr)(1-s^2)}.
    \end{aligned}
\end{equation}
It is easy to see that 
\begin{equation}\label{eq: limit of J}
    J(s) \rightarrow \frac{e^{\s(X-Y)}}{(\s-\rr)(\s+\rr)(2\s)}.
\end{equation}
Computing derivatives of $I(s)$ gives $I(1) = 0$, $I'(1) = 0$, and $$I''(1) = \rho^{-1}\frac{-2N(1+\sqrt{q})\sqrt{q}}{(1-\sqrt{q})^3} = -\rho^{-1}c_0N,$$
which implies that
\begin{equation}\label{eq: limit of I}
\begin{aligned}
    I(1+\s\rho^{-1}) &=  \frac{-c_0N}{2}\frac{\s^2}{(c_0N/2)} + \mathcal{O}((c_0 N/2)^{-1/3})\rightarrow -\s^2.
\end{aligned}
\end{equation}
Combining \eqref{eq: limit of J} and \eqref{eq: limit of I} gives
\eqref{eq: n=2,n'=2}.
Same computations show that
\begin{equation}
\begin{aligned}
    &\lim_{N\rightarrow \infty}\frac{\rho^{-3}s^{1-k+\ell}}{(s-r)(1-sr)(1-s^2)}\rho^{-1}\left(-(1+k)s + \frac{N\sqrt{q}}{(1-\sqrt{q}/s)}+\frac{N\sqrt{q}s^2}{(1-\sqrt{q}s)}\right)\\
    &= \frac{e^{\s(Y-X)}}{(\s-\rr)(\s+\rr)(2\s)}\left(-X+\s^2\right).
\end{aligned}
\end{equation}
Similar arguments justify that
\begin{equation}
    \lim_{N\rightarrow \infty} (c_0N/2)^{-1}\widehat{R}(1,k;2,\ell) = \widetilde{R}(1,X;2,Y), \quad \lim_{N\rightarrow \infty} (c_0N/2)^{-1}\widehat{R}(2,k;1,\ell) = \widetilde{R}(2,X;1,Y).
\end{equation}
The limit of $(c_0N/2)^{-2/3}\widehat{R}(1,k;1,\ell)$ is immediate.
Furthermore, the upper bound for $\widehat{K}_{22}$ comes from poles $(z,w) = (s,r)$, $(z,w) = (1/r,1/s)$ and terms like $s^{k-\ell},$ $s^{\ell-k}$. The error term $e^{\mathcal{O}((c_0N/2)^{-1/3})}$ is absorbed into $e^{\epsilon(x+Y)}$ for some $\epsilon>0.$
\end{proof}

\begin{lem}\label{lem: asymp pfaffian}
    For any $d, y\in \Z_{\geq 0}$ such that $d = \frac{2\sqrt{q}}{1-\sqrt{q}}N + (c_0N/2)^{1/3}\tilde{d}$ and $y = (c_0N/2)^{1/3}\tilde{y}$, we have
    \begin{equation}
        \lim_{N\rightarrow \infty} \mathrm{Pf}(J - \widehat{K})_{\ell^2(T)} = \mathrm{Pf}(J - \widetilde{K})_{\mathbb{L}^2(\widetilde{T})},
    \end{equation}
    where $T = \{1\} \times \{d+1, \dots\} \cup \{2\} \times \{d+y+1, \dots\}$ and $\widetilde{T} = \{1\}\times (\tilde{d},\infty) \cup \{2\} \times (\tilde{d} + \tilde{y},\infty)$.
\end{lem}

\begin{proof}
    By Lemma \ref{lem: K11,K12 steepest} and Lemma \ref{lem: limit and bound of K22}, we have the entry-wise limit of the kernel and for any $u>0,$ there exists $\epsilon >0$ and $\kappa > -\s + \epsilon$ such that for the $nn'-$block, the following upper bound holds for all $X,Y \geq -u$
    \begin{equation}\label{eq: asymp kernel upper bound}
        \begin{aligned}
        &\Bigg|\begin{pmatrix}
        \rho^{n+n'+2}\widehat{K}_{11}(n,k;n',\ell) & \rho^{n-n'+1}\widehat{K}_{12}(n,k;n',\ell)\\
        \rho^{n'-n+1}\widehat{K}_{21}(n,k;n',\ell) & \rho^{-n-n'}\widehat{K}_{22}(n,k;n',\ell)
    \end{pmatrix}\Bigg| \\
    &\leq \begin{pmatrix}
        Ce^{-\kappa(X+Y)} & \begin{gathered}Ce^{-\kappa X -(\s - \epsilon)Y}\\
            {}+ Ce^{-(\s-\epsilon)(Y-X)}\Id_{X\leq Y}\Id_{n=1,n'=2}\end{gathered}\\[4pt]
        \begin{gathered}Ce^{-(\s - \epsilon)X-\kappa Y}\\
            {}+ Ce^{-(\s-\epsilon)(X-Y)}\Id_{Y\leq X}\Id_{n= 2, n' = 1}\end{gathered} & Ce^{-(\s-\epsilon)(X+Y)}
    \end{pmatrix}.
    \end{aligned}
\end{equation}
Consider the conjugation matrix $\widetilde{M}(n,X) = \text{diag}(e^{\widetilde{u}_n X}, e^{-\widetilde{u}_n X})$. We get the following upper bound on the conjugated kernel
\begin{equation}
    \begin{aligned}
        |\widetilde{M}\widehat{K}\widetilde{M}|(n,X;n',Y) \leq \begin{pmatrix}
            Ce^{(-\kappa + \widetilde{u}_n)X + (-\kappa + \widetilde{u}_{n'})Y} & Ce^{(-\kappa +\widetilde{u}_n)X -(\s-\epsilon + \widetilde{u}_{n'}) Y} + \mathcal{A}\\
            Ce^{-(\s -\epsilon+\widetilde{u}_n)X-(\kappa - \widetilde{u}_{n'})Y} + \mathcal{B} & Ce^{-(\s-\epsilon + \widetilde{u}_n)X -(\s-\epsilon + \widetilde{u}_{n'})Y}
        \end{pmatrix},
    \end{aligned}
\end{equation}
where
\begin{equation}
    \begin{aligned}
        &\mathcal{A} = Ce^{X(\widetilde{u}_n-(\widetilde{u}_n +\widetilde{u}_{n'})/2) + Y((\widetilde{u}_n +\widetilde{u}_{n'})/2 - \widetilde{u}_{n'}) + (Y-X)(-\s+\epsilon - (\widetilde{u}_{n}+\widetilde{u}_{n'})/2)}\Id_{X\leq Y}\Id_{n=1, n'=2},\\
        &\mathcal{B} = Ce^{Y(\widetilde{u}_{n'}-(\widetilde{u}_n +\widetilde{u}_{n'})/2) + X((\widetilde{u}_n +\widetilde{u}_{n'})/2 - \widetilde{u}_{n}) + (X-Y)(-\s+\epsilon - (\widetilde{u}_{n}+\widetilde{u}_{n'})/2)}\Id_{Y\leq X}\Id_{n=2, n'=1}.
    \end{aligned}
\end{equation}
It suffices to choose 
\begin{equation}
\begin{aligned}
-\s+\epsilon<\widetilde{u}_n, \widetilde{u}_{n'}< \kappa, \quad \text{and}\quad \widetilde{u}_n < \widetilde{u}_{n'}\,\,\, \text{ if } n<n', \quad\quad  \widetilde{u}_{n'} < \widetilde{u}_{n}\,\,\, \text{ if } n'<n.
\end{aligned}
\end{equation}
For shorter notation, we write $\xi = \rho^{-1} =(c_0N/2)^{-1/3}$. For each $m\in \mathbb{N},$ we define $E_{\xi}^{(j_1,\dots,j_m)} := E_{\xi}^{j_1}\times\cdots\times E_{\xi}^{j_m}$ where $E_{\xi}^{j_i} := \{\tilde{d} + a\xi : a \in \mathbb{N}\}$ if $j_i = 1$ and $E_{\xi}^{j_i} := \{\tilde{d} + \tilde{y} + a\xi : a \in \mathbb{N}\}$ if $j_i = 2$. Then we have
\begin{equation}
    \begin{aligned}
        &\Pf\left( J - \widehat{K}\right)_{\ell^2(T)}= 1 + \sum_{m = 1}^{\infty}\frac{(-1)^m}{m!} \sum_{j_1,\dots, j_m = 1}^2 \sum_{x_1 = d_{j_1}}^{\infty}\cdots \sum_{x_m = d_{j_m}}^{\infty} \Pf\left(\widehat{K}(j_u,x_u; j_v,x_v)\right)_{u,v =1}^{m}\\
        &=1 + \sum_{m = 1}^{\infty}\frac{(-1)^m}{m!} \!\!\!\!\!\sum_{j_1,\dots, j_m = 1}^2 \int_{ \widetilde{d}_{j_1}}^{\infty}\!\!\!\!\!\!\cdots \!\!\int_{ \widetilde{d}_{j_m}}^{\infty}\sum_{\Vec{y} \in E_{\xi}^m} \!\!\Pf\left(\overline{K}(j_u,X_u; j_v,X_v)\right)_{u,v =1}^{m}\Id_{[y_1-\xi,y_1)\times\cdots\times [y_m-\xi, y_m)}(X_1,\dots,X_m) \prod_{i=1}^m dX_i\\
        &\leq 1 + \sum_{m = 1}^{\infty}\frac{(-1)^m}{m!} \!\!\!\!\!\sum_{j_1,\dots, j_m = 1}^2\int_{ \widetilde{d}_{j_1}}^{\infty}\cdots \int_{ \widetilde{d}_{j_m}}^{\infty} C^m(2m)^{m/2}\prod_{i=1}^m e^{-\sum_{i = 1}^m\theta X_i} \prod_{i=1}^m dX_i <\infty,
    \end{aligned}
\end{equation}
where $\overline{K}$ is the rescaled kernel as in the LHS of \eqref{eq: asymp kernel upper bound} and $\theta = \min\{\kappa-\widetilde{u}_n, \kappa-\widetilde{u}_{n'}, |\widetilde{u}_n - \widetilde{u}_{n'}|/2, \s-\epsilon+ \min(\widetilde{u}_n,\widetilde{u}_{n'})\} > 0$.
Hence, by dominated convergence theorem and Hadamard's bound, we get the convergence of Fredholm Pfaffian.
\end{proof}

\begin{thm}\label{thm: high-density final}
Fix $\sqrt{q} \in (0,1)$. Let $\mathcal{X}$ be an independent geometric random variable with rate $rs$. For any $\tilde{s} <0$ and $\tilde{r} \in (-\infty, \tilde{s})\cup (\tilde{s}, -\tilde{s})$, we have
\begin{equation}
\begin{aligned}
    \lim_{N\rightarrow \infty}\!\!\Pb\!\Bigg( \frac{G(N+1,N)-\frac{2\sqrt{q}}{1-\sqrt{q}}N}{(c_0N/2)^{1/3}} \leq \tilde{d}, &\frac{G(N+2,N)-\mathcal{X}-\frac{2\sqrt{q}}{1-\sqrt{q}}N}{(c_0N/2)^{1/3}} \leq \tilde{d} \Bigg)\\
    &= -(\rr + \s)\int_{0}^{\infty}\Pf(J- \widetilde{K})_{\mathbb{L}^2(\widetilde{T})} e^{(\rr + \s)\tilde{y}}d\tilde{y},
\end{aligned}
\end{equation}
where $\widetilde{T} = \{1\} \times (\tilde{d}, \infty) \cup \{2\} \times (\tilde{d}+\tilde{y}, \infty).$
\end{thm}
\begin{proof}
     Assume the scaling of parameters in \eqref{eq: scale}.  By Lemma \ref{lem: asymp pfaffian}, we get that
     \begin{equation}
         \begin{aligned}
             &\lim_{N\rightarrow \infty}\sum_{\tilde{y}\in\{m\rho^{-1}:m\in\mathbb{N}\}}\Pf(J- K)_{\ell^2(T)} {(-(\rr+\s)\rho^{-1})}\left(1+{(\rr+\s)}{\rho^{-1}} + \mathcal{O}(\rho^{-2})\right)^{\tilde{y}\rho }\\
             &= -(\rr+\s)\int_{0}^{\infty}\Pf\left( J - \widetilde{K}\right)_{\mathbb{L}^2(\widetilde{T})}e^{(\rr+\s)\tilde{y}}d\tilde{y}.
         \end{aligned}
     \end{equation}
\end{proof}

\begin{proof}[Proof of Theorem~\ref{thm: high-density phase theorem}]
    This is an immediate consequence of Corollary \ref{eq: finite N high} and Theorem \ref{thm: high-density final}.
\end{proof}

\section{Asymptotics for maximal-current regime}\label{sec: max current}

In this section, we compute the special case when $\s = 0,$ $\rr<0$ or equivalently, $s = 1$ and $r\in (0,1)$ for the finite $N$ distribution. We will continue with the same letter for each function but add a subscript $M$ to denote the maximal-current phase. Let
\begin{equation}
    \begin{aligned}
        &\widetilde{A}_M(1,\zeta) = \zeta, \quad \widetilde{A}_M(2,\zeta) = \zeta^2,\quad
        \widetilde{B}_M(1,\zeta) = -\zeta, \quad \widetilde{B}_M(2,\zeta) = \zeta^2,\\
    \end{aligned}
\end{equation}

\begin{lem}\label{lem: maximal-current limit of K22}
For any fixed $\rr<0,$ we have the following limit:
\begin{equation}
    \begin{aligned}
        \lim_{\s \rightarrow 0}\widetilde{K}_{11} (n,X;n',Y) = \widetilde{K}_{11}^M(n,X;n',Y),\quad
        \lim_{\s \rightarrow 0}\widetilde{K}_{12} (n,X;n',Y) =\widetilde{K}_{12}^M  (n,X;n',Y).
    \end{aligned}
\end{equation}
    \begin{equation}
        \lim_{\s \rightarrow 0}\widetilde{K}_{22}(n,X;n',Y)= \widetilde{K}_{22}^M(n,X;n',Y).
    \end{equation}
Fix any $u>0.$ For any $\delta>0,$ there exists $C>0$ independent of $X,Y$, $\delta >\epsilon >0$ such that for all $|\s| < \epsilon$ and for all $X,Y > - u$, the following upper bounds hold
\begin{equation}
    \begin{aligned}
        &\big|\widetilde{K}_{11} (n,X;n',Y) \big| \leq Ce^{-\delta(X+Y)}, \quad \big|\widetilde{K}_{12} (n,X;n',Y) \big| \leq Ce^{-\delta X + \epsilon Y} + C\Id_{\{n=1,n'=2, X\leq Y\}}e^{\epsilon(Y-X)}, \\ &\big|\widetilde{K}_{22} (n,X;n',Y) \big| \leq Ce^{\epsilon(X+Y)}.
    \end{aligned}
\end{equation}
\end{lem}
\begin{proof}
Limits of $\widetilde{K}_{11}$ and $\widetilde{K}_{12}$ are straightforward since taking $\s \rightarrow 0$ does not cross the contour and does not introduce undefined terms. Limit of $\widetilde{K}_{22}$ follows from a direct evaluation. We need to group all terms involving $(2\s)^{-1}$ to obtain the limit.
\end{proof}

\begin{prop}\label{prop: max curr limit}
For any $\tilde{r} <0$, we have
    \begin{equation}\label{eq: high density s to 1 limit}
    \lim_{\s\rightarrow 0^{-}}-(\rr + \s)\int_{0}^{\infty}\Pf(J- \widetilde{K})_{\mathbb L^2(\widetilde{T})} e^{(\rr + \s)\tilde{y}}d\tilde{y} = -\rr\int_{0}^{\infty}\Pf(J- \widetilde{K}^{M})_{\mathbb L^2(\widetilde{T})} e^{\rr\tilde{y}}d\tilde{y},
\end{equation}
where $\widetilde{K}^{M} = \begin{pmatrix}
    \widetilde{K}_{11}^M & {\widetilde{K}_{12}^M}\\
    \widetilde{K}_{21}^M & \widetilde{K}_{22}^M\\
\end{pmatrix}$ and $\widetilde{T} = \{1\} \times (\tilde{d}, \infty) \cup \{2\} \times (\tilde{d}+\tilde{y}, \infty)$.
\end{prop}

\begin{proof}
By Lemma \ref{lem: maximal-current limit of K22}, we get the point-wise convergence of each kernel entry and their upper bounds. Then by
Hadamard's bound and dominated convergence theorem, we get the limit of the Fredholm Pfaffian. The convergence of integral is then justified by dominated convergence theorem because $\rr <0$ and $|\Pf(J - \widetilde{K})|\leq 1$ for all $\s <0.$
\end{proof}

We want to upgrade the above result to the convergence of stationary geometric LPP under critical scaling.
Recall the definition of the stationary model \eqref{eq:stat_wts}. Parameters under the maximal-current regime require $s = 1$, $r<1$ and $X\sim \text{Geo}(r)$. Let $G_{r,1}^{stat}(N,N)$ denote the LPP value from $(1,1)$ to $(N,N)$. We claim the following proposition.

\begin{prop}\label{prop: maximal-current asymptotic equal}
For any $\tilde{d} \in \R$, $\rr <0,$ we have
\begin{equation}
    \lim_{\s\rightarrow 0}\lim_{N\rightarrow \infty} P\left(\frac{G_{r,s}^{stat}(N,N)-\frac{2\sqrt{q}}{1-\sqrt{q}}N}{(c_0N/2)^{1/3}} \leq \tilde{d}\right) = \lim_{N\rightarrow \infty} \lim_{s\rightarrow 1}P\left(\frac{G_{r,s}^{stat}(N,N)-\frac{2\sqrt{q}}{1-\sqrt{q}}N}{(c_0N/2)^{1/3}} \leq \tilde{d}\right).
\end{equation}
\end{prop}

The $s\rightarrow 1$ limits of $\widehat{K}_{11}$ and $\widehat{K}_{12}$ is straightforward and we state them without proof. We will use a superscript $M$ to denote the maximal-current phase, i.e., $s = t = 1$.
\begin{equation}\label{eq: maximal current K11 and K12}
    \begin{aligned}
        \widehat{K}_{11}^M(n,k;n',\ell):= &-\oint\limits_{\Gamma_{\sqrt{q}}} \frac{dw}{2\pi\I} \oint\limits_{\Gamma_{1/\sqrt{q}}} \frac{dz}{2\pi\I}  \frac{w^{\ell-1}}{z^{k}}\frac{H(z)(1-1/z)^n(1-w)^{n'}}{H(w)}\frac{(zw-1)(z-r)(1-wr)}{(z^2-1)(1-w^2)(z-w)},\\
        \widehat{K}_{12}^M(n,k; n', \ell) := &-\oint\limits_{\Gamma_{r,1,\sqrt{q}}} \frac{dw}{2\pi\I} \oint\limits_{\Gamma_{1/\sqrt{q}}} \frac{dz}{2\pi\I}\frac{w^{\ell-1}}{z^{k}} \frac{H(z)(1-1/z)^n}{H(w)(1-1/w)^{n'}}\frac{(zw-1)(z-r)}{(z^2-1)(z-w)(w-r)}\\
        &- \Id_{k\leq \ell}\Id_{n = 1, n' = 2}.\\
    \end{aligned}
\end{equation}

We define a few auxiliary functions for $22$-entry.
\begin{equation}
    \begin{aligned}
        \mathsf{f}_1(k,\ell) &:= \frac{r^{\ell-k+1}}{(1-r)^4} + \frac{r^{\ell+1}}{H(r)(1-r)^2}\left(\frac{-1}{(1-r)^2} - \frac{k}{1-r} + \frac{2N\sqrt{q}}{(1-\sqrt{q})(1-r)}\right),\\
        \mathsf{f}_2(k,\ell) &:= -\left(\frac{1}{(1-r)^2} - \frac{1+\ell}{1-r} + \frac{2N\sqrt{q}}{(1-\sqrt{q})(1-r)}\right)\frac{r^{k+1}}{H(r)(1-r)^2} + \frac{(1+r)r^{k+1}}{H(r)(1-r)^4},\\
        \mathsf{f}_3(k,\ell) &:= \left(\frac{1}{(1-r)^2} - \frac{1+\ell}{1-r} + \frac{2N\sqrt{q}}{(1-\sqrt{q})(1-r)}\right)\left(\frac{-k}{1-r} + \frac{2N\sqrt{q}}{(1-\sqrt{q})(1-r)} + \frac{r}{(1-r)^2}\right),\\
        \mathsf{f}_4(k,\ell) &:=\frac{k^2}{(1 - r)(-1 + r)}
- \frac{4kN\sqrt{q}}{(1 - \sqrt{q})(1 - r)(-1 + r)}
+ \frac{4N^2 q}{(1 - \sqrt{q})^2 (1 - r)(-1 + r)}
\\
&\quad+ \frac{r}{(1 - r)^2 (-1 + r)}
- \frac{2kr}{(1 - r)^2 (-1 + r)}
+ \frac{4N\sqrt{q}\, r}{(1 - \sqrt{q})(1 - r)^2 (-1 + r)}
+ \frac{2r^2}{(1 - r)^3 (-1 + r)}.
    \end{aligned}
\end{equation}

\begin{equation}
    \begin{aligned}
        &\mathcal{R}_{22}^M(k,\ell) := \mathsf{f}_1(k,\ell) + \mathsf{f}_2(k,\ell)+ \mathsf{f}_3(k,\ell)+ \mathsf{f}_4(k,\ell).
\end{aligned}
\end{equation}

\begin{equation}
    \begin{aligned}
        &\mathcal{R}_{12}^M(k,\ell) := -\frac{r^{\ell+1}}{H(r)(1-r)^3} + \frac{r^{\ell-k+1}}{(1-r)^3} \\
        &\quad + \left(\frac{-1}{(1-r)^2} + \frac{1+\ell}{1-r} - \frac{2N\sqrt{q}}{(1-\sqrt{q})(1-r)}\right)\frac{1}{1-r}
        +\frac{1-2k}{(1-r)^2}\\
        &+ \frac{4N\sqrt{q}}{(1-\sqrt{q})(1-r)^2} + \frac{2r}{(1-r)^3}\\
        &\quad +\left(\frac{-1}{(1-r)^2} + \frac{1+\ell}{1-r} - \frac{2N\sqrt{q}}{(1-\sqrt{q})(1-r)}\right)\frac{-r^k}{(1-r)H(r)}
        - \frac{(1+r)r^k}{(1-r)^3H(r)}.
    \end{aligned}
\end{equation}

\begin{equation}
    \begin{aligned}
        &\mathcal{R}_{21}^M(k,\ell) :=\frac{r^{\ell-k}}{-(1-r)^3} + \frac{r^{\ell}}{H(r)(1-r)}\left(\frac{1}{(1-r)^2} + \frac{k}{(1-r)} - \frac{2N\sqrt{q}}{(1-\sqrt{q})(1-r)}\right)\\
        &+ \frac{r^{k+1}}{H(r)(1-r)^3} + \left(\frac{k}{1-r} - \frac{2N\sqrt{q}}{(1-\sqrt{q})(1-r)} - \frac{r}{(1-r)^2}\right)\frac{1}{(1-r)}.\\
        &\mathcal{R}_{11}^M(k,\ell) := \frac{r^{\ell}}{H(r)(1-r)^2} - \frac{r^{\ell-k}}{(1-r)^2}+\frac{1}{(1-r)^2} - \frac{r^k}{H(r)(1-r)^2}.
    \end{aligned}
\end{equation}

\begin{lem}\label{lemma: maximal current K22}
    For any $\sqrt{q}\in (0,1)$, $0<r <1$, we have
    \begin{equation}
        \begin{aligned}
            \widehat{K}_{22}^M(n,k;n',\ell) := \lim_{s\rightarrow 1}\widehat{K}_{22}(n,k;n',\ell)  = \mathcal{B}(n,k;n',\ell) +\Id_{k\leq \ell}\mathcal{A}(n,k;n',\ell) - \Id_{k> \ell}\mathcal{A}(n',\ell;n,k),
        \end{aligned}
    \end{equation}
    where 
    \begin{equation}
        \begin{aligned}
            &\mathcal{A}(n,k;n',\ell) :=\Id_{n = n' = 2}\mathcal{R}_{22}^M(k,\ell) + \Id_{n =2, n' = 1}\mathcal{R}_{21}^M(k,\ell) + \Id_{n =1, n' = 2}\mathcal{R}_{12}^M(k,\ell) + \Id_{n =1, n' = 1}\mathcal{R}_{11}^M(k,\ell),\\
            & \mathcal{B}(n,k;n',\ell) = -\oint \limits_{\Gamma_{\sqrt{q}}}\frac{dw}{2\pi\I} \oint \limits_{\Gamma_{1/r,1}}\frac{dz}{2\pi\I} \mathcal{E}(z,w) - \oint \limits_{\Gamma_{\sqrt{q},r,1}}\frac{dw}{2\pi\I} \oint \limits_{\Gamma_{1/\sqrt{q}}}\frac{dz}{2\pi\I} \mathcal{E}(z,w),\\        \end{aligned}
    \end{equation}
    where
    \begin{equation}
        \mathcal{E}(z,w) = \frac{w^{\ell-1}}{z^{k}} \frac{H(z)}{H(w)(1-1/w)^{n'}(1-z)^{n}}\frac{(zw-1)}{(z-w)(1-zr)(w-r)}.
    \end{equation}
\end{lem}

\begin{proof}
    We explain the most difficult case $n = n' = 2$. All other cases of $n,n'$ follow similar computations. It is easy to check that
    \begin{equation}
    \lim_{s\rightarrow 1}\widehat{h}_{22}(2,k;2,\ell) = \mathcal{B}(2,k;2,\ell).
\end{equation}
Consider the remaining terms. It is clear that
    \begin{equation}
        \begin{aligned} &\Id_{n=n'=2}\left(\widehat{C}(2,k;2,\ell) + \widehat{D}(2,k;2,\ell)\right)\\
&= \frac{
\frac{sr}{sr-1}\frac{s^2}{s^2-1}
\left(r^{\ell}s^k-r^ks^{\ell}\right)
}
{H(s)H(r)(s-r)(1-sr)}
+
\frac{
\frac{sr}{sr-1}H(s)
}
{H(r)(s-r)^2(1-s^2)}
\left(
\frac{r^{\ell}}{s^{k-1}}
-
\frac{r^k}{s^{\ell-1}}
\right)\\
&\rightarrow \mathsf{f}_1(k,\ell) + \mathsf{f}_2(k,\ell) - \frac{r^{\ell-k+1}}{(1-r)^4} \text{ as }s\rightarrow 1.
        \end{aligned}
    \end{equation}
By direct evaluations, we also get that
\begin{equation}
    \begin{aligned}
        &\lim_{s\rightarrow 1} \widehat{R}(2,k;2,\ell) = \mathsf{f}_3 + \mathsf{f}_4 +  \frac{r^{\ell-k+1}}{(1-r)^4}.
    \end{aligned}
\end{equation}
Therefore, 
\begin{equation}
    \lim_{s\rightarrow 1}  \Id_{n=n'=2}\left(\widehat{C}(2,k;2,\ell) + \widehat{D}(2,k;2,\ell)\right) + \widehat{R}(2,k;2,\ell) = \mathcal{R}_{22}^M(k,\ell).
\end{equation}
\end{proof}
We use the notation $\widehat{K}^M =  \begin{pmatrix}
    \widehat{K}_{11}^M & \widehat{K}_{12}^M\\
    \widehat{K}_{21}^M & \widehat{K}_{22}^M
\end{pmatrix}$.

\begin{prop}\label{prop: maximal-current asymptotic}
For any $\sqrt{q} \in (0,1)$, $r\in (0,1),$ we have
    \begin{equation}\label{eq: max curr finite N pre limit}
    \begin{aligned}
    \Pb\left( G_{r,1}^{stat}({N+2,N+2}) \leq d\right) =  \sum_{y=0}^{\infty}\Pf\left(J - \widehat{K}^M\right)_{\ell^{2}(T)}(1-r)r^y.\\
\end{aligned}
\end{equation}
\end{prop}

\begin{proof}
    Given \eqref{eq: maximal current K11 and K12} and Lemma \ref{lemma: maximal current K22}, the convergence of Pfaffian, $\lim_{s \rightarrow 1} \Pf\left( J - \widehat{K} \right) = \Pf\left( J - \widehat{K}^M \right),$ can be proved following the same procedure as in Proposition \ref{Fredconv}.
\end{proof}

Next, we find the asymptotic limits of the kernel. Recall the scaling of the parameters in \eqref{eq: scale}.

\begin{lem}\label{lem: max current K22}
    For any $\rr <0$, we have
    \begin{equation}
    \begin{aligned}
    \lim_{N\rightarrow \infty} \rho^{n+n'+2}\widehat{K}_{11}^M(n,k;n',\ell) = \widetilde{K}_{11}^M(n,X;n',Y),\quad
            \lim_{N\rightarrow \infty} \rho^{n-n'+1}\widehat{K}_{12}^M(n,k;n',\ell) = \widetilde{K}_{12}^M(n,X;n',Y)
    \end{aligned}
\end{equation}
\begin{equation}
    \begin{aligned}
        \lim_{N\rightarrow \infty} \rho^{-n-n'} \widehat{K}_{22}^M(n,k;n',\ell) = \widetilde{K}_{22}^M(n,X;n',Y).
    \end{aligned}
    \end{equation}
    Fix any $u>0$. For any $\delta>0,$ there exist a constant $C>0$ independent of $X,Y$, $N_0 \in \mathbb{N}$, and $\delta>\epsilon>0$ such that for any $X,Y \geq -u$ and $N\geq N_0$, we have
    \begin{equation}
    \begin{aligned}
    &|\rho^{n+n'+2}\widehat{K}_{11}^M(n,k;n',\ell)| \leq Ce^{-\delta(X+Y)},\\
    &|\rho^{n-n'+1}\widehat{K}_{12}^M(n,k;n',\ell)| \leq Ce^{-\delta X + \epsilon Y} + \Id_{X\leq Y}\Id_{n=1,n'=2}Ce^{\epsilon (Y-X)},\\
    &|\rho^{-n-n'}\widehat{K}_{22}^M(n,k;n',\ell)|\leq Ce^{\epsilon(X+Y)}.
    \end{aligned}
\end{equation}
\end{lem}

\begin{proof}
We explain the most difficult case $n = n' = 2$. Limits and upper bounds of all other cases of $n,n'$ follow similar arguments. Under the scaling, we get that
\begin{equation}
    \begin{aligned}
        &\lim_{N\rightarrow \infty}\rho^{-4}\mathsf{f}_1(k,\ell)= \frac{e^{\rr (Y-X)}}{\rr^4} + e^{-\frac{\rr^3}{3} + Y\rr}\left(\frac{-1}{\rr^4} +\frac{X}{\rr^3}\right),\quad
        \lim_{N\rightarrow \infty} \rho^{-4}\mathsf{f}_2(k,\ell)= e^{-\frac{\rr^3}{3} + X\rr} \left(\frac{1}{\rr^4}- \frac{Y}{\rr^3}\right),\\
        &\lim_{N\rightarrow \infty} \rho^{-4}\mathsf{f}_3(k,\ell)= \frac{1}{\rr^4} + \frac{1}{\rr^3}(X+Y) + \frac{XY}{\rr^2},\quad \lim_{N\rightarrow \infty}\rho^{-4}\mathsf{f}_4(k,\ell)= -\frac{X^2}{\rr^2} - \frac{2X}{\rr^3} - \frac{2}{\rr^4},
    \end{aligned}
\end{equation}
which explains the limit of $\mathcal{R}_{22}$. The dominant terms of $\rho^{-4}\mathcal{R}_{22}$ are $r^{\ell-k}$, $Xr^{\ell}/H(r)$, $Yr^{k}/H(r)$, $XY$ and $X^2$, which are bounded by $Ce^{\epsilon(X+Y)}$ for some $C,\epsilon>0.$ The constant $\epsilon$ also covers the error term $e^{\mathcal{O}(\rho^{-1}(X+Y))}$.

\end{proof}

\begin{proof}[Proof of Proposition~\ref{prop: maximal-current asymptotic equal}]
Lemma \ref{lem: max current K22} gives the point-wise limit of the kernel entries and their upper bounds. Applying Hadamard's bound and dominated convergence theorem to the right hand side of \eqref{eq: max curr finite N pre limit} gives the same limit as the right hand side of \eqref{eq: high density s to 1 limit}.
\end{proof}

\begin{proof}[Proof of Theorem ~\ref{thm: maximal-current phase theorem}] This is an immediate consequence of Proposition \ref{prop: max curr limit} and Proposition \ref{prop: maximal-current asymptotic equal}. 
\end{proof}

\section{Diagonal distribution of half-space stationary geometric LPP in low-density regime}\label{section: low}

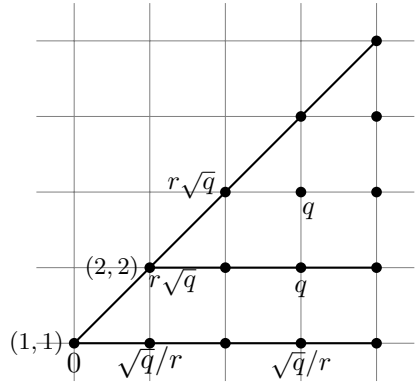
\begin{figure}[H]
\centering
\begin{tikzpicture}
\draw[very thin, gray] (-0.5,-0.5) grid (4.5,4.5);
    \foreach \x in {1,...,4} {
        \foreach \y in {0,1,...,\x} {
            \fill (\x,\y) circle (2pt);
        }
    }
    \fill (0,0) circle (2pt);
    \draw[thick,-] (0,0) -- (4,0);
    \draw[thick,-] (0,0) -- (4,4);
    \draw[thick,-] (1,1) -- (4,1);
    \fill (0,0) circle (1pt);
    \node at (0,-0.25) {$0$};
     
    \node at (1,-0.25) {$\sqrt{q}/r$};
    \node at (-0.5,0) {\footnotesize{$(1,1)$}};
    \node at (0.5,1) {\footnotesize{$(2,2)$}};
    \node at (1.3,0.81) {\small{$r\sqrt{q}$}};
    \fill (1,1) circle (1pt);
    \node at (3,-0.28) {\small{$\sqrt{q}/r$}};
    \node at (3,0.72) {\small{$q$}};
    \node at  (1.55,2.1) {\small{$r\sqrt{q}$}};
    \node at (3.1,1.75) {\small{$q$}};
\end{tikzpicture}
\caption{Geometric LPP model as described in $\eqref{eq: product stationary weigths}$.}
\label{fig: two param stationary prev}
\end{figure}

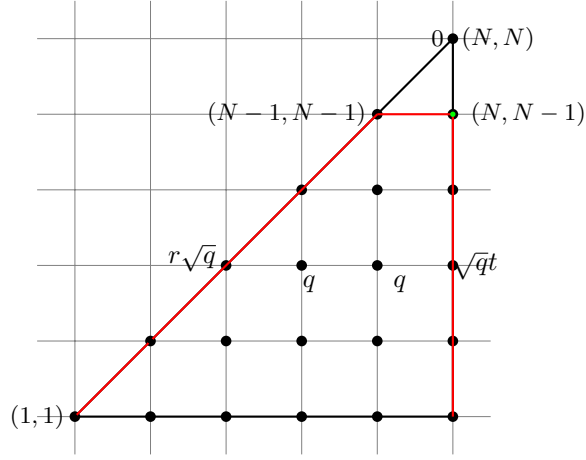
\begin{figure}[htbp]
\centering
\begin{tikzpicture}
\draw[very thin, gray] (-0.5,-0.5) grid (5.5,5.5);
    \foreach \x in {1,...,5} {
        \foreach \y in {0,1,...,\x} {
            \fill (\x,\y) circle (2pt);
        }
    }
    \fill (0,0) circle (2pt);
    \draw[thick,-] (0,0) -- (5,0);
    \draw[thick,-] (0,0) -- (5,5);

    \draw[thick,-] (5,5) -- (5,0);
    \draw[thick,red] (4,4) -- (5,4);
    \draw[thick,red] (5,4) -- (5,0);
    \draw[thick,red] (0,0) -- (4,4);
    \fill (0,0) circle (1pt);

    \fill[green] (5,4) circle (1pt);

    \node at (-0.5,0) {\footnotesize{$(1,1)$}};
    \node at (5.6,5) {\footnotesize{$(N,N)$}};
    \node at (2.8,4) {\footnotesize{$(N-1,N-1)$}};
    \node at (6,4) {\footnotesize{$(N,N-1)$}};
    \node at (4.8,5) {\footnotesize{$0$}};

    \fill (1,1) circle (1pt);
    \node at (5.3,2) {\small{$\sqrt{q}t$}};
    \fill (1,1) circle (1pt);

    \node at (4.3,1.75) {\small{$q$}};
    \fill (1,1) circle (1pt);

    \node at  (1.55,2.1) {\small{$r\sqrt{q}$}};
    \node at (3.1,1.75) {\small{$q$}};
\end{tikzpicture}
\caption{Product stationary geometric LPP model as shown in \eqref{eq: product stationary weigths flipped}.}
\label{fig: two param stationary flipped}
\end{figure}

Consider the half-space geometry $\mathbb{H}=\{(i,j)\in\Z^2: 1\leq j\leq i\}$. Let $q_0 = r,$ $q_1= 1/r$, and $q_2=\dots = q_{N+1} = \sqrt{q}$. We define the product stationary half-space geometric LPP model and specify its weights as follows:
\begin{equation} \label{eq: product stationary weigths}
  \omega_{i, j} = \begin{cases}
    0, & \textrm{if}\ i=j=1, \\
    \mathrm{Geo}\left( \sqrt{q}/r \right), & \textrm{if }i \geq 2, j=1, \\
    \mathrm{Geo}\left( r\sqrt{q} \right), & \textrm{if } 2\leq i=j\leq N, \\
    \mathrm{Geo}(q), &\textrm{otherwise},
  \end{cases}
\end{equation}
where parameters satisfy $\sqrt{q} \in (0,1),$ and $r\in (\sqrt{q},1/\sqrt{q})$. Let $G^{stat}_{r}(N,N)$ denote the LPP from $(1,1)$ to $(N,N)$ of the above model \eqref{eq: product stationary weigths}. Since the distribution of $G^{stat}_{r}(N,N)$ remains invariant under permutations of $q_0$ and $q_1$, it suffices to restrict to $r\in (\sqrt{q},1].$

Following a similar idea as before, we permute the parameters and define a new model. Let $q_0 = r,$ $q_1=\dots = q_{N-1} = \sqrt{q}$, and $q_ = t$. We define weights as follows:
\begin{equation} \label{eq: product stationary weigths flipped}
  \omega_{i, j} = \begin{cases}
    0, & \textrm{if}\ i=j=N, \\
    \mathrm{Geo}\left( t\sqrt{q} \right), & \textrm{if }i = N, 1\leq j\leq N-1, \\
    \mathrm{Geo}\left( r\sqrt{q} \right), & \textrm{if } 1\leq i=j\leq N-1, \\
    \mathrm{Geo}(q), &\textrm{otherwise},
  \end{cases}
\end{equation}
Let $\G(N,N)$ denote the LPP from $(1,1)$ to $(N,N)$ of the model \eqref{eq: product stationary weigths flipped}. The following identity in distribution is immediate and the Fredholm Pfaffian formula is taken from \cite[Proposition 4.2]{BBCS}.

\begin{prop} \label{thm: low kernel}
    The distribution function of $G_{r}^{stat}(N,N) $ is given by:
    \begin{equation}
        \Pb(G_{r}^{stat}(N,N) \leq d) = \lim_{t\rightarrow 1/r}\Pb\left(  \G(N,N-1) \leq d\right).
    \end{equation}
    The CDF is given by a Fredholm Pfaffian formula
    \begin{equation}
        \Pb\left(G_{N,N-1} \leq d \right) = \Pf(J- \KK)_{\ell^2(T)},
    \end{equation}
    where
    \begin{equation}
    \begin{aligned}
        &\Pf\left( J - \KK\right)_{\ell^2(T)} = 1 + \sum_{m = 1}^{\infty}\frac{(-1)^m}{m!}  \sum_{x_1 = d+1}^{\infty}\cdots \sum_{x_m = d+1}^{\infty} \Pf\left(K(x_u, x_v)\right)_{u,v =1}^{m},\\
        &T = \{d+1,d+2, \dots \}.
    \end{aligned}
    \end{equation}
    The $2\times 2$ matrix correlation kernel $\KK$ is given by:
    \begin{equation}\label{eq: one param kernel}
        \begin{aligned}
            \KK_{11}(k,\ell) &= \oint \oint \frac{dz}{2\pi\I}\frac{dw}{2\pi\I}\frac{V(z)V(w)}{z^{k}w^{\ell}}\frac{(z-w)(z-r)(w-r)}{(z^2-1)(w^2-1)(zw-1)},\\
            \KK_{12}(k,\ell) &= \oint \oint \frac{dz}{2\pi\I}\frac{dw}{2\pi\I}\frac{w^{\ell-1}}{z^{k}}\frac{V(z)}{V(w)}\frac{(zw-1)(z-r)}{(z^2-1)(z-w)(w-r)} = -\KK_{21}(l,k),\\
            \KK_{22}(k,\ell) &= \oint \oint \frac{dz}{2\pi\I}\frac{dw}{2\pi\I}\frac{z^{k-1}w^{\ell-1}}{V(z)V(w)}\frac{(z-w)}{(zw-1)(z-r)(w-r)},\\
        \end{aligned}
    \end{equation}
    where 
    \begin{equation}
        V(z) = H(z)(1-t/z), \quad H(z) = \frac{(1-\sqrt{q}/z)^{N-1}}{(1-\sqrt{q}z)^{N-1}},
    \end{equation}
    and where the contours are positively oriented circles centered around the origin satisfying the following conditions:
\begin{itemize}
   \item for $\KK_{11}$, $1 < |z|, |w| < \min \{1/t,1/\sqrt{q}\}$;
   \item for $\KK_{12}$, $\max \{ r,t,\sqrt{q}\} < |w|<1,$ and $1 < |z| < \min\{1/t,1/\sqrt{q}\}$;
   \item for $\KK_{22}$, $\max \{r,t,\sqrt{q}\} < |w|, |z|$ and $1 < |zw|$.
  \end{itemize}
\end{prop}

\begin{lem}\label{lem: one param kernel reformulated}
For any $\sqrt{q} \in (0,1),$ $r,t\in (\sqrt{q},1/\sqrt{q}),$ the kernel $\mathcal{K}$ can be reformulated as follows:
    \begin{equation}\label{eq: low-density kernel}
        \begin{aligned}
            \KK_{11}(k,\ell) = &-\oint\limits_{\Gamma_{\sqrt{q}}} \frac{dw}{2\pi\I} \oint \limits_{\Gamma_{1/\sqrt{q}}}\frac{dz}{2\pi\I}\frac{w^{\ell-1}}{z^{k+1}}\frac{H(z)}{H(w)}\frac{(z-t)(1-tw)(zw-1)(z-r)(1-wr)}{(z^2-1)(1-w^2)(z-w)},\\
            \KK_{12}(k,\ell) = &-\oint\limits_{\Gamma_{\sqrt{q},t,r}} \frac{dw}{2\pi\I} \oint \limits_{\Gamma_{1/\sqrt{q}}}\frac{dz}{2\pi\I}\frac{w^{\ell}}{z^{k+1}}\frac{H(z)}{H(w)}\frac{(z-t)(zw-1)(z-r)}{(w-t)(z^2-1)(z-w)(w-r)} = -\KK_{21}(l,k),\\
            \KK_{22}(k,\ell) 
            = &-\oint \limits_{\Gamma_{\sqrt{q}}} \!\! \frac{dw}{2\pi\I}\!\!\! \oint \limits_{\Gamma_{1/\sqrt{q},1/r,1/t}}\!\!\!\!\!\!\!\!\! \frac{dz}{2\pi\I} \,\mathcal{U}(z,w)-\oint \limits_{\Gamma_{r,t}}\!\! \frac{dw}{2\pi\I}\!\!\! \oint \limits_{\Gamma_{1/\sqrt{q}}}\!\!\! \frac{dz}{2\pi\I}\,\mathcal{U}(z,w)\\
            &+\frac{1}{(t-r)(1-tr)}\left(\frac{r^{\ell}t^k}{H(r)H(t)} - \frac{t^{\ell}r^{k}}{H(r)H(t)}\right)
            + \left(\frac{t^{\ell-k}}{(1-tr)(t-r)} - \frac{r^{\ell-k}}{(1-tr)(t-r)}\right)\Id_{k\leq \ell}\\
            &+ \left(-\frac{t^{k-\ell}}{(1-tr)(t-r)} + \frac{r^{k-\ell}}{(1-tr)(t-r)}\right)\Id_{k> \ell},
\end{aligned}
\end{equation}
where 
\begin{equation}
    \mathcal{U}(z,w) = \frac{w^{\ell}}{z^{k}}\frac{H(z)}{H(w)}\frac{(zw-1)}{(1-tz)(w-t)(z-w)(1-zr)(w-r)}.
\end{equation}
\end{lem}

\begin{proof}
For $\KK_{11}$, we do a change of variable $w \rightarrow w^{-1}$ and then deform the $w$-contour to enclose $\sqrt{q}$ and deform the $z$-contour to enclose $1/\sqrt{q}$. For $\KK_{12},$ we deform contours to include $z,w$ poles. For $\KK_{22}$, we first change the variable $z\rightarrow z^{-1}$ and then move the $w$-contour inside the $z$-contour, which requires adding the residue at $w=z$.
We see that 
\begin{equation}
    \begin{aligned}
        &\text{Res}(\KK_{22},w=z)=\oint \limits_{\Gamma_{t,r,0}} \frac{dz}{2\pi\I} \frac{z^{\ell-k}(1-z^2)}{(1-tz)(z-t)(1-zr)(z-r)}\\
        &= \left(\frac{t^{\ell-k}}{(1-tr)(t-r)} - \frac{r^{\ell-k}}{(1-tr)(t-r)}\right)\Id_{k\leq \ell} + \left(-\frac{t^{k-\ell}}{(1-tr)(t-r)} + \frac{r^{k-\ell}}{(1-tr)(t-r)}\right)\Id_{k> \ell},\\
        &\left(- \oint \limits_{\Gamma_{ t}} \frac{dw}{2\pi\I} \oint \limits_{\Gamma_{1/r}} \frac{dz}{2\pi\I} - \oint \limits_{\Gamma_{ r}} \frac{dw}{2\pi\I} \oint \limits_{\Gamma_{1/t}} \frac{dz}{2\pi\I}\right)\frac{w^{\ell}}{z^{k}}\frac{H(z)}{H(w)}\frac{(zw-1)}{(1-tz)(w-t)(z-w)(1-zr)(w-r)}\\
        &= \frac{1}{(t-r)(1-tr)}\left(\frac{r^{\ell}t^k}{H(r)H(t)} - \frac{t^{\ell}r^{k}}{H(r)H(t)}\right).\\
 \end{aligned}
\end{equation}
\end{proof}

We define the following functions for the $t\rightarrow 1/r$ limit.
\begin{equation}
    \begin{aligned}
        \wKK_{11}(k,\ell)&:=-\oint\limits_{\Gamma_{\sqrt{q}}} \frac{dw}{2\pi\I} \oint \limits_{\Gamma_{1/\sqrt{q}}}\frac{dz}{2\pi\I}\frac{w^{\ell-1}}{z^{k+1}}\frac{H(z)}{H(w)}\frac{(z-r^{-1})(1-r^{-1}w)(zw-1)(z-r)(1-wr)}{(z^2-1)(1-w^2)(z-w)},\\
        \wKK_{12}(k,\ell) &:= -\oint\limits_{\Gamma_{\sqrt{q},r,1/r}} \frac{dw}{2\pi\I} \oint \limits_{\Gamma_{1/\sqrt{q}}}\frac{dz}{2\pi\I}\frac{w^{\ell}}{z^{k+1}}\frac{H(z)}{H(w)}\frac{(z-r^{-1})(zw-1)(z-r)}{(w-r^{-1})(z^2-1)(z-w)(w-r)} = -\wKK_{21}(\ell,k),\\
    \end{aligned}
\end{equation}
\begin{equation}
    \begin{aligned}
        \widehat{\mathsf{F}}(k,\ell) &:= \frac{kr^{1+k-\ell}}{1-r^2}-\frac{kr^{1-k+\ell}}{1-r^2} + \left(\frac{(1-N)r^{k-\ell}\sqrt{q}}{1-r^2}- \frac{(1-N)r^{\ell-k}\sqrt{q}}{1-r^2}\right)\left(\frac{1}{1-\sqrt{q}/r} + \frac{r^2}{1-\sqrt{q}r}\right),\\
        \widehat{\mathcal{Q}}(k,\ell) &:=  -\oint \limits_{\Gamma_{\sqrt{q}}} \frac{dw}{2\pi\I} \oint \limits_{\Gamma_{1/\sqrt{q},1/r,r}} \frac{dz}{2\pi\I}\widehat{\mathcal{U}}(z,w) -\oint \limits_{\Gamma_{r,1/r}} \frac{dw}{2\pi\I} \oint \limits_{\Gamma_{1/\sqrt{q}}} \frac{dz}{2\pi\I}\widehat{\mathcal{U}}(z,w),\\
    \end{aligned}
\end{equation}
where
\begin{equation}
    \widehat{\mathcal{U}}(z,w) = \frac{w^{\ell}}{z^{k}}\frac{H(z)}{H(w)}\frac{(zw-1)}{(1-r^{-1}z)(w-r^{-1})(z-w)(1-zr)(w-r)}.
\end{equation}

\begin{lem}
    For any $\sqrt{q}\in (0,1)$, $r\in (\sqrt{q},1],$ we have the following limits of kernel:
    \begin{equation}
        \begin{aligned}
            &\lim_{t\rightarrow 1/r} \KK_{11}(k,\ell) = \wKK_{11}(k,\ell), \quad \lim_{t\rightarrow 1/r} \KK_{12}(k,\ell) = \wKK_{12}(k,\ell),\\
            &\lim_{t\rightarrow 1/r} \KK_{22}(k,\ell) = \widehat{\mathcal{Q}}(k,\ell) + \widehat{\mathsf{F}}(k,\ell) \Id_{k \leq \ell} - \widehat{\mathsf{F}}(\ell,k) \Id_{k > \ell} .
        \end{aligned}
    \end{equation}
\end{lem}

\begin{proof}
It is straightforward to take $t\rightarrow r^{-1}$ limit for $\KK_{11}$ and $\KK_{12}$. $\widehat{\mathcal{Q}}(k,\ell)$ captures the limit of residues at $w=\sqrt{q}$ and $z = 1/\sqrt{q}$.
The special residue term $z=w$ in $\KK_{22}$ together with residues $\text{Res}(w=t,z=1/r)$ and $\text{Res}(w=r,z=1/t)$ give the following limit. When $k\leq \ell$,
\begin{equation}
    \begin{aligned}
        &\lim_{t\rightarrow 1/r}\frac{t^{\ell-k}}{(1-tr)(t-r)} - \frac{r^{\ell-k}}{(1-tr)(t-r)}+  \frac{1}{(t-r)(1-tr)}\left(\frac{r^{\ell}t^k}{H(r)H(t)} - \frac{t^{\ell}r^{k}}{H(r)H(t)}\right)\\
        &= \lim_{t\rightarrow 1/r}\frac{t^{\ell-k}}{(1-tr)(t-r)} -\frac{1}{(t-r)(1-tr)}\frac{t^{\ell}r^{k}}{H(r)H(t)}\\
        &\qquad + \lim_{t\rightarrow 1/r}  \frac{-r^{\ell-k}}{(1-tr)(t-r)}+  \frac{1}{(t-r)(1-tr)}\frac{r^{\ell}t^k}{H(r)H(t)}\\
        &= \widehat{\mathsf{F}}(k,\ell).
    \end{aligned}
\end{equation}
\end{proof}

\begin{lem}
    Fix any $\sqrt{q} \in (0,1)$, $r\in (\sqrt{q},1],$ and any $0<\epsilon \ll 1$ satisfying $\epsilon < \frac{1-\sqrt{q}/r}{2\sqrt{q} + 1/r}$. There exists a positive constant $C$ independent of $k,\ell$ such that the following upper bounds hold:
    \begin{equation}
        \begin{aligned}
            &|\KK_{11}(k,\ell)| \leq Ce^{(k+\ell)\log(\sqrt{q} +\epsilon)}, \quad |\KK_{12}(k,\ell)| \leq Ce^{k\log(\sqrt{q} +\epsilon) + \ell\log(1/r + \epsilon)},\\
            &|\KK_{21}(k,\ell)| \leq Ce^{\ell\log(\sqrt{q} +\epsilon) + k \log(1/r+\epsilon)}, \quad |\KK_{22}(k,\ell)| \leq Ce^{(k+\ell)\log(1/r + \epsilon)}.\\
        \end{aligned}
    \end{equation}
    Moreover, we have the limit
    \begin{equation}
        \lim_{t\rightarrow 1/r} \Pf(J- \KK)_{\ell^2(T)} = \Pf(J - \wKK)_{\ell^2(T)},
    \end{equation}
    where $T = \{d+1,d+2,\dots\}$.
\end{lem}
\begin{proof}
    Upper bounds for the kernel follow from similar analysis as in Lemma \ref{lem: K11 and K12} and Lemma \ref{lem: K22}. The condition on $\epsilon$ ensures that $-\log(\sqrt{q} +\epsilon) > \log(1/r+\epsilon)$, which allows us to apply Hadamard's bound and dominated convergence theorem for the convergence of Fredholm Pfaffian.
\end{proof}

\subsection{Asymptotic analysis under low-density regime}\label{section: low asymp}
Recall the scaling of the parameters $k,\ell,r,z,w,d$ in \eqref{eq: scale}. Recall the definition of $\rho = (c_0N/2)^{1/3}$.
\begin{lem}\label{lem: low pfaffian convergence}
    For any $u>0$ and any $\rr \leq 0,$ the following limits hold uniformly over $X,Y \in [-u,u]$,
    \begin{equation}
        \begin{aligned}
            &\lim_{N\rightarrow \infty} \rho^4 \wKK_{11}(k,\ell) = \tKK_{11}(X,Y),\quad \lim_{N\rightarrow \infty} \rho \wKK_{12}(k,\ell) = \tKK_{12}(X,Y),\\
            &\lim_{N\rightarrow \infty} \rho \wKK_{21}(k,\ell) = \tKK_{21}(X,Y),\quad \lim_{N\rightarrow \infty} \rho^{-2} \wKK_{22}(k,\ell) = \tKK_{22}(X,Y).\\
        \end{aligned}
    \end{equation}
    For any $\kappa >-\rr,$ there exists a positive constant $C>0$ independent of $X,Y$, $N_0 \in \mathbb{N}$, and $\epsilon>0$ with $0<-\rr + \epsilon < \kappa$ such that for any $X,Y \geq -u$ and $N\geq N_0$, we have
    \begin{equation}
        \begin{aligned}
            &|\rho^4 \wKK_{11}(X,Y)| \leq Ce^{-\kappa(X+Y)}, \quad |\rho \wKK_{12}(X,Y)| \leq Ce^{-\kappa X + (-\rr + \epsilon)Y},\\
            &|\rho \wKK_{21}(X,Y)| \leq Ce^{ (-\rr + \epsilon)X -\kappa Y}, \quad |\rho^{-2} \wKK_{22}(X,Y)| \leq Ce^{(-\rr + \epsilon)(X+Y)}.
        \end{aligned}
    \end{equation}
    Moreover, for any $d\in \Z_{\geq 0}$ such that $d = \frac{2\sqrt{q}}{1-\sqrt{q}}N + (c_0N/2)^{1/3}\tilde{d}$, we have
    \begin{equation}
        \begin{aligned}
            \lim_{N\rightarrow \infty}\Pf\left(J - \wKK\right)_{\ell^2(T)} = \Pf\left(J - \tKK\right)_{\mathbb{L}^2(\widetilde{T})}, 
        \end{aligned}
    \end{equation}
    where $\tKK = \begin{pmatrix}
        \tKK_{11} & \tKK_{12}\\
        \tKK_{21} & \tKK_{22}
    \end{pmatrix}$, $T = \{d+1,d+2,\dots\},$ and $\widetilde{T}=(\tilde{d}, \infty)$.
\end{lem}

\begin{proof}
    Asymptotics of $\wKK_{11}$, $\wKK_{12}$ and $\widehat{\mathcal{Q}}$ can be obtained via the steepest descent method as in Lemma \ref{lem: K11,K12 steepest}. The special term is $\widehat{\mathsf{F}}$.
    Consider $\varepsilon(r) = kr -N \sqrt{q}\left( \frac{1}{1-\sqrt{q}/r} + \frac{r^2}{1-\sqrt{q}r} \right)$ and recall that $k = \frac{2\sqrt{q}}{1-\sqrt{q}}N + \rho X$. 
    Then we have $\varepsilon(1) =\rho X$, $\varepsilon'(1) = \rho X$ and
    $\varepsilon''(1) = -\frac{2N(1+\sqrt{q})\sqrt{q}}{(1-\sqrt{q})^3} = -Nc_0.$
    Then we have 
    \begin{equation}
        \varepsilon(r) = \rho X + X\rr -\frac{1}{2}Nc_0\frac{\rr^2}{(c_0N/2)^{2/3}} + \mathcal{O}(N^{-1/3}).
    \end{equation}
    Therefore,
    \begin{equation}
        \lim_{N\rightarrow \infty} \rho^{-2}\widehat{\mathcal{F}}(X,Y) = \lim_{N\rightarrow\infty}\rho^{-2}\frac{\varepsilon(r)}{1-r^2}\left(r^{k-\ell} - r^{\ell-k}\right) = \widetilde{\mathsf{F}}(X,Y).
    \end{equation}
    For $\wKK_{11}$, poles $z = 1/\sqrt{q},$ and $w = \sqrt{q}$ are sent to infinity, resulting in the $z$-contour ($w$-contour resp.) can stay on the right of any $\kappa >0$ (the left of any $-\kappa <0$). For $\wKK_{12}$, the $w$-contour has to enclose poles $w = r,1/r$, which become $\omega = \rr,-\rr$ after the scaling. The dominant term comes from the pole $\omega = -\rr$. For $\wKK_{22}$, the dominant term comes from $\widehat{\mathsf{F}}$, which is upper bounded by $Cr^{-(k+\ell)}$. The convergence of Pfaffian is justified by Hadamard's bound and dominated convergence theorem as in \cite[Lemma 6.3]{zeng}. 
\end{proof}

\begin{proof}[Proof of Theorem ~\ref{thm: low-density phase theorem}] This is an immediate consequence of Proposition \ref{thm: low kernel} and Lemma \ref{lem: low pfaffian convergence}.
\end{proof}

\section{Matching of log-gamma with geometric in the low-density phase}
\label{sec:matching-low-density}
In this section, we show that our low-density formula agrees with the
formula derived from the stationary half-space log-gamma polymer in
\cite{zeng2026}. We compare the two finite $N$ representations after degenerating
both models to half-space exponential LPP. The geometric and log-gamma
degenerations yield kernels $K^{\mathrm{gexp}}$ and $K^{\mathrm{lexp}}$,
respectively. Although these kernels differ entrywise, their Fredholm
Pfaffians are related by a first-order differential operator in $d$.

We retain the notation $t$ and $r$ for the first-row and diagonal
parameters, respectively, but now use them for the exponential model;
$\alpha$ denotes its bulk parameter. We initially assume
$t,r,\alpha>0$, $t\neq r$, and $N\geq1$. The stationary specialization
$t\to-r$ is handled in Section \ref{section: low}.

Let $\Gamma$ be a compact, positively oriented, piecewise smooth contour
in the open left half-plane, enclosing $-t,-r,-\alpha$ and excluding $0$.
Then $Z+W\neq0$ for all $Z,W\in\Gamma$. Unless otherwise indicated,
all double contour integrals below are over $\Gamma\times\Gamma$. We first define two kernels $K^{\mathrm{gexp}}$ and $K^{\mathrm{lexp}}$.
The proof for the following lemma is provided in the appendix.

\begin{lem}[Geometric-to-exponential degeneration]
\label{lem:geometric-exp-kernel}
Fix $\alpha, t,r >0$ $t\neq r$ and $t,r < \alpha$. The geometric-to-exponential limit of the low-density geometric kernel is the
skew-symmetric $2\times 2$ matrix kernel $K^{\mathrm{gexp}}$ with entries
\begin{align}
K^{\mathrm{gexp}}_{11}(X,Y)
={}&\oint_{\Gamma}\frac{dZ}{2\pi\mathrm i}
\oint_{\Gamma}\frac{dW}{2\pi\mathrm i}\,
 e^{XZ+YW}
 \left(\frac{\alpha-Z}{\alpha+Z}\right)^{N-1}
 \left(\frac{\alpha-W}{\alpha+W}\right)^{N-1}
 \frac{(t-Z)(t-W)(r-Z)(r-W)(Z-W)}{4ZW(Z+W)},
\label{eq:Kshift11}\\
K^{\mathrm{gexp}}_{12}(X,Y)
={}&-\oint_{\Gamma}\frac{dZ}{2\pi\mathrm i}
\oint_{\Gamma}\frac{dW}{2\pi\mathrm i}\,
 e^{XZ+YW}
 \left(\frac{\alpha-Z}{\alpha+Z}\right)^{N-1}
 \left(\frac{\alpha-W}{\alpha+W}\right)^{N-1}
 \frac{(t-Z)(r-Z)(Z-W)}{2Z(t+W)(r+W)(Z+W)},
\label{eq:Kshift12}\\
K^{\mathrm{gexp}}_{21}(X,Y)
={}&-K^{\mathrm{gexp}}_{12}(Y,X),
\label{eq:Kshift21}\\
K^{\mathrm{gexp}}_{22}(X,Y)
={}&\oint_{\Gamma}\frac{dZ}{2\pi\mathrm i}
\oint_{\Gamma}\frac{dW}{2\pi\mathrm i}\,
 e^{XZ+YW}
 \left(\frac{\alpha-Z}{\alpha+Z}\right)^{N-1}
 \left(\frac{\alpha-W}{\alpha+W}\right)^{N-1}
 \frac{Z-W}{(Z+t)(W+t)(Z+W)(Z+r)(W+r)}
\notag\\
&\quad+\mathbf 1_{\{Y>X\}}
 \frac{e^{-t(Y-X)}-e^{-r(Y-X)}}{r^2-t^2}
 -\mathbf 1_{\{X>Y\}}
 \frac{e^{-t(X-Y)}-e^{-r(X-Y)}}{r^2-t^2}.
\label{eq:Kshift22}
\end{align}
\end{lem}

{
We next take the zero-temperature limit of the half-space log-gamma kernel.  The scaling below is the finite-$N$ version of the log-gamma-to-exponential scaling used in \cite[Section~17]{zeng2026}.  We include the short kernel computation in the appendix, since the precise finite-$N$ degeneration needed here is not stated there in the form used below.
}

\begin{lem}[Log-gamma-to-exponential degeneration]
\label{lem:loggamma-exp-kernel}
Fix $\alpha,t,r>0$, with $t\neq r$ and $t,r<\alpha$.  In the half-space log-gamma kernel of \cite[Theorem~5.10]{imamura2022solvable}, specialize the diagonal and bulk parameters so that the  parameters are $\varepsilon r$, $\varepsilon t$, and $N-1$ copies of o$\varepsilon\alpha$, and rescale the threshold and kernel variables by
\[
\tau=\frac d\varepsilon,\qquad X=\frac x\varepsilon,\qquad Y=\frac y\varepsilon,
\]
while writing the spectral variables as $Z=\varepsilon z$ and $W=\varepsilon w$.  After the Pfaffian-preserving conjugation induced by
\[
D_\varepsilon=\operatorname{diag}(\varepsilon^{1/2},\varepsilon^{-1/2}),
\qquad \det D_\varepsilon=1,
\]
the rescaled log-gamma kernel converges, locally uniformly for $x,y$ in compact subsets of $\mathbb R$, to the skew-symmetric matrix kernel $K^{\mathrm{lexp}}$. In the following formulas, we rename the limiting variables $(x,y,z,w)$ as $(X,Y,Z,W)$, respectively. The kernel entries are
\begin{align}
K^{\mathrm{lexp}}_{11}(X,Y)
={}&
\oint_{\Gamma}\frac{dZ}{2\pi\mathrm{i}}
\oint_{\Gamma}\frac{dW}{2\pi\mathrm{i}}\,
e^{XZ+YW}
\frac{r-Z}{r+Z}\frac{t-Z}{t+Z}
\left(\frac{\alpha-Z}{\alpha+Z}\right)^{N-1}
\frac{r-W}{r+W}\frac{t-W}{t+W}
\left(\frac{\alpha-W}{\alpha+W}\right)^{N-1}
\frac{1}{4ZW}\frac{Z-W}{Z+W},
\label{eq:Kexp11}\\[1ex]
K^{\mathrm{lexp}}_{12}(X,Y)
={}&
-\oint_{\Gamma}\frac{dZ}{2\pi\mathrm{i}}
\oint_{\Gamma}\frac{dW}{2\pi\mathrm{i}}\,
e^{XZ+YW}
\frac{r-Z}{r+Z}\frac{t-Z}{t+Z}
\left(\frac{\alpha-Z}{\alpha+Z}\right)^{N-1}
\frac{r-W}{r+W}\frac{t-W}{t+W}
\left(\frac{\alpha-W}{\alpha+W}\right)^{N-1}
\frac{1}{4Z}\frac{Z-W}{Z+W},
\label{eq:Kexp12}\\[1ex]
K^{\mathrm{lexp}}_{21}(X,Y)
={}&
-\oint_{\Gamma}\frac{dZ}{2\pi\mathrm{i}}
\oint_{\Gamma}\frac{dW}{2\pi\mathrm{i}}\,
e^{XZ+YW}
\frac{r-Z}{r+Z}\frac{t-Z}{t+Z}
\left(\frac{\alpha-Z}{\alpha+Z}\right)^{N-1}
\frac{r-W}{r+W}\frac{t-W}{t+W}
\left(\frac{\alpha-W}{\alpha+W}\right)^{N-1}
\frac{1}{4W}\frac{Z-W}{Z+W},
\label{eq:Kexp21}\\[1ex]
K^{\mathrm{lexp}}_{22}(X,Y)
={}&
\oint_{\Gamma}\frac{dZ}{2\pi\mathrm{i}}
\oint_{\Gamma}\frac{dW}{2\pi\mathrm{i}}\,
e^{XZ+YW}
\frac{r-Z}{r+Z}\frac{t-Z}{t+Z}
\left(\frac{\alpha-Z}{\alpha+Z}\right)^{N-1}
\frac{r-W}{r+W}\frac{t-W}{t+W}
\left(\frac{\alpha-W}{\alpha+W}\right)^{N-1}
\frac{1}{4}\frac{Z-W}{Z+W}.
\label{eq:Kexp22}
\end{align}
\end{lem}

The following proposition also explains how the shift argument is embedded within a Fredholm Pfaffian formula via modifying the kernel.
\begin{prop}
\label{prop:exact-matching-finite-N}
For every fixed $N$ and all $t,r,\alpha>0$ and $d\in\mathbb R_{>0}$,
\begin{equation}
\operatorname{Pf}(J-K^{\mathrm{gexp}})_{\mathbb L^2(d,\infty)}
=
\left(1+\frac{1}{t+r}\,\partial_d\right)
\operatorname{Pf}(J-K^{\mathrm{lexp}})_{\mathbb L^2(d,\infty)}.
\label{eq:exact-matching-main}
\end{equation}
\end{prop}
We will prove the above identity directly from the Fredholm Pfaffian
representations, without using their interpretation as distribution
functions of the same exponential LPP model. Our aim is to identify the
limiting formula in Theorem~\ref{thm: low-density phase theorem} with that
in \cite[Theorem~1.8]{zeng2026}. Comparing the asymptotic expressions
directly is difficult because divergent terms (coming from the $t\rightarrow -r$ limit) cancel differently in the
two representations which breaks the nice kernel structure.

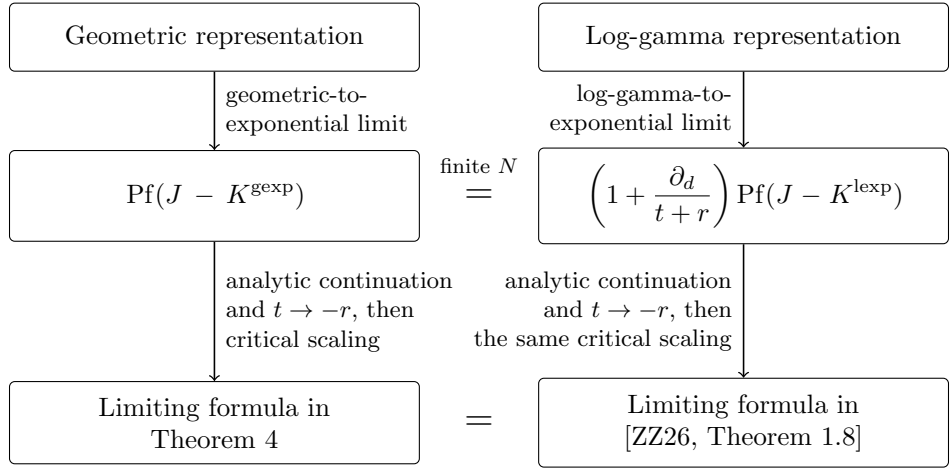
\begin{figure}[htbp]
\centering
\begin{tikzpicture}[
  box/.style={
    draw,
    rounded corners=2pt,
    align=center,
    text width=5cm,
    minimum height=0.9cm,
    inner sep=6pt,
    font=\small
  },
  route/.style={
    >=to,
    ->,
    semithick
  },
  annotation/.style={
    font=\footnotesize,
    align=left,
    inner sep=4pt
  }
]

% Original representations.
\node[box] (geom) at (0,0)
  {Geometric representation};

\node[box] (lg) at (7,0)
  {Log-gamma representation};

% Exact matching at the finite-N exponential level.
\node[box, minimum height=1.2cm] (shift) at (0,-2.1)
  {$\displaystyle
    \operatorname{Pf}(J-K^{\mathrm{gexp}})$};

\node[box, minimum height=1.2cm] (exp) at (7,-2.1)
  {$\displaystyle
    \left(1+\frac{\partial_d}{t+r}\right)
    \operatorname{Pf}(J-K^{\mathrm{lexp}})$};

\node[font=\Large] at (3.5,-2.1) {$=$};
\node[font=\scriptsize] at (3.5,-1.7) {finite $N$};

% Limiting formulas.
\node[box] (ours) at (0,-5.1)
  {Limiting formula in\\
   Theorem~\ref{thm: low-density phase theorem}};

\node[box] (theirs) at (7,-5.1)
  {Limiting formula in\\
   \cite[Theorem~1.8]{zeng2026}};

\node[font=\Large] at (3.5,-5.1) {$=$};

% Degenerations to the exponential formulas.
\draw[route] (geom.south) --
  node[midway, right, annotation]
    {geometric-to-\\exponential limit}
  (shift.north);

\draw[route] (lg.south) --
  node[midway, left, annotation, align=right]
    {log-gamma-to-\\exponential limit}
  (exp.north);

% Stationary specialization, followed by critical scaling.
\draw[route] (shift.south) --
  node[midway, right, annotation]
    {analytic continuation\\
     and $t\to-r$, then\\
     critical scaling}
  (ours.north);

\draw[route] (exp.south) --
  node[midway, left, annotation, align=right]
    {analytic continuation\\
     and $t\to-r$, then\\
     the same critical scaling}
  (theirs.north);

\end{tikzpicture}
\caption{Comparison through the exact finite-$N$ identity
\eqref{eq:exact-matching-main}. Applying analytic continuation to $t=-r$ and the same critical scaling yields the same asymptotic. }
\label{fig:finite-N-matching-strategy}
\end{figure}

Adapting a similar approach of \cite[Section~18]{zeng2026}, we construct a
transformation operator that establishes an exact identity between the
two forms of Fredholm Pfaffian formulas at finite $N$, before the stationary
specialization. The analytic continuations established in Section \ref{section: low} and in
\cite{zeng2026} then extend this identity to the stationary limit
$t\to-r$. Applying the same critical scaling to both sides gives identical
limits. The asymptotic of the Fredholm Pfaffian using exponential kernels gives the same limit as Theorem~\ref{thm: low-density phase theorem};
the log-gamma-side limit is the formula in
\cite[Theorem~1.8]{zeng2026}. Figure~\ref{fig:finite-N-matching-strategy}
summarizes this argument.

\subsection{Operator factorization of the exponential kernel}
\label{subsec:operator-factorization-exp}

On the contour space
\[
\mathcal H_{\Gamma}:=\mathbb L^2(\Gamma,|dZ|)
\]
we use the bilinear contour pairing
\begin{equation}
\langle f,g\rangle_{\Gamma}
:=\oint_{\Gamma}f(Z)g(Z)\frac{dZ}{2\pi\mathrm i}.
\label{eq:contour-pairing}
\end{equation}
Let
\begin{equation}
\mathbb L^2_d:=\mathbb L^2((d,\infty),dX)\otimes\mathbb C^2,
\label{eq:L2d-J}
\end{equation}
\begin{equation}
\vartheta(Z):=
\left(\frac{\alpha-Z}{\alpha+Z}\right)^{N-1}
\frac{(t-Z)(r-Z)}{(t+Z)(r+Z)}\frac{1}{2Z}.
\label{eq:theta-def}
\end{equation}
Define the operators
$\mathsf F:\mathcal H_{\Gamma}\to\mathbb L^2_d$ and
$\mathsf C:\mathcal H_{\Gamma}\to\mathcal H_{\Gamma}$ by
\begin{align}
(\mathsf Fh)(X)
&:=\oint_{\Gamma}\vartheta(Z)e^{XZ}
 \begin{pmatrix}1\\-Z\end{pmatrix}h(Z)\frac{dZ}{2\pi\mathrm i},
\label{eq:F-def}\\
\mathsf C(Z,W)&:=-\frac{Z-W}{Z+W}.
\label{eq:C-def}
\end{align}
Since $\Gamma$ is compact and contained in the left half-plane, $\mathsf F$ is
Hilbert--Schmidt. The kernel of $\mathsf C$ is continuous on
$\Gamma\times\Gamma$, so $\mathsf C$ is Hilbert--Schmidt, and in particular
bounded.

Let $\mathfrak S_2(\mathcal H_\Gamma)$ and
$\mathfrak S_1(\mathcal H_\Gamma)$ denote the classes of
Hilbert--Schmidt and trace-class operators on $\mathcal H_\Gamma$,
equipped with the Hilbert--Schmidt norm $\|\cdot\|_{\mathfrak S_2}$
and the trace norm $\|\cdot\|_{\mathfrak S_1}$, respectively.
For any $A,B\in\mathfrak S_2(\mathcal H_\Gamma)$,
their product $AB$ belongs to $\mathfrak S_1(\mathcal H_\Gamma)$ because
\[
\|AB\|_{\mathfrak S_1}
\leq
\|A\|_{\mathfrak S_2}\,
\|B\|_{\mathfrak S_2}.
\]

For the operator-theoretic justification of the Fredholm manipulations, one may
apply the standard determinant-one Pfaffian conjugation
\begin{equation}
D_{\gamma}(X):=\operatorname{diag}(e^{\gamma X},e^{-\gamma X}),
\qquad 0<\gamma<-\sup_{Z\in\Gamma}\Re Z.
\label{eq:gauge-conjugation}
\end{equation}
Because $\det D_{\gamma}(X)=1$, this conjugation leaves every term of the
Fredholm Pfaffian series unchanged, while the spatial kernels acquire
exponential decay in $X+Y$.
\begin{lem}
\label{lem:Kexp-factorization}
The exponential kernel admits the factorization
\begin{equation}
K^{\mathrm{lexp}}=-\mathsf F\mathsf C\mathsf F^{\mathsf T}.
\label{eq:Kexp-factorization}
\end{equation}
Moreover, if
\begin{equation}
\mathsf M_d:=\mathsf F^{\mathsf T}J^{-1}\mathsf F,
\label{eq:Md-def-abstract}
\end{equation}
then
\begin{equation}
\mathsf M_d(Z,W)
=\frac{Z-W}{Z+W}\,
 \vartheta(Z)\vartheta(W)e^{d(Z+W)}.
\label{eq:Md-explicit}
\end{equation}
Consequently,
\begin{equation}
\operatorname{Pf}(J-K^{\mathrm{lexp}})_{\mathbb L^2(d,\infty)}
=\det(\Id+\mathsf C\mathsf M_d)_{\mathcal H_{\Gamma}}^{1/2}.
\label{eq:Pfexp-det}
\end{equation}
\end{lem}

\begin{proof}
From~\eqref{eq:F-def}--\eqref{eq:C-def},
\begin{align*}
-\mathsf F\mathsf C\mathsf F^{\mathsf T}(X,Y)
={}&\oint_{\Gamma}\frac{dZ}{2\pi\mathrm i}
\oint_{\Gamma}\frac{dW}{2\pi\mathrm i}\,
\vartheta(Z)\vartheta(W)e^{XZ+YW}
\frac{Z-W}{Z+W}
\begin{pmatrix}
1&-W\\-Z&ZW
\end{pmatrix}.
\end{align*}
Substitution of~\eqref{eq:theta-def} gives exactly
\eqref{eq:Kexp11}--\eqref{eq:Kexp22}, and proves
\eqref{eq:Kexp-factorization}.

Since $J^{-1}=-J$, for $Z,W\in\Gamma$,
\begin{align*}
\mathsf M_d(Z,W)
&=\int_d^{\infty}
\vartheta(Z)e^{XZ}(1,-Z)
J^{-1}
\vartheta(W)e^{XW}\begin{pmatrix}1\\-W\end{pmatrix}dX\\
&=(W-Z)\vartheta(Z)\vartheta(W)
\int_d^{\infty}e^{X(Z+W)}dX\\
&=\frac{Z-W}{Z+W}\vartheta(Z)\vartheta(W)e^{d(Z+W)},
\end{align*}
which is~\eqref{eq:Md-explicit}. The kernels of $\mathsf C$ and
$\mathsf M_d$ are bounded on the compact set $\Gamma\times\Gamma$; hence both
operators are Hilbert--Schmidt, and $\mathsf C\mathsf M_d$ is trace class.
Then by the identity $\det(\Id + AB) = \det(\Id + BA)$, we get
\begin{align*}
\operatorname{Pf}(J+\mathsf F\mathsf C\mathsf F^{\mathsf T})^2
&=\det(\Id+J^{-1}\mathsf F\mathsf C\mathsf F^{\mathsf T})\\
&=\det(\Id+\mathsf C\mathsf F^{\mathsf T}J^{-1}\mathsf F)\\
&=\det(\Id+\mathsf C\mathsf M_d).
\end{align*}
Both sides tend to $1$ as $d\to+\infty$, which fixes the square-root branch and
proves~\eqref{eq:Pfexp-det}.
\end{proof}

\subsection{The spatial transformation}
\label{subsec:spatial-transformation}

For $a>0$, let
\begin{equation}
P_a(X,Y):=\mathbf 1_{\{Y>X\}}e^{-a(Y-X)}.
\label{eq:Pa-def}
\end{equation}
Let
\begin{equation}
R(X,Y):=P_tP_r=\mathbf 1_{\{Y>X\}}
\frac{e^{-t(Y-X)}-e^{-r(Y-X)}}{r-t}, \qquad S:=\partial_XR.
\label{eq:R-kernel}
\end{equation}
Set
\begin{equation}
\mathsf L:=
\begin{pmatrix}
-(\partial_X+t)(\partial_X+r)&0\\
S&-R
\end{pmatrix}.
\label{eq:L-def}
\end{equation}

We introduce the notation
\begin{equation}
Q:=\frac{1}{t+r}R,
\qquad
Q(X,Y)=\mathbf 1_{\{Y>X\}}
\frac{e^{-t(Y-X)}-e^{-r(Y-X)}}{r^2-t^2}.
\label{eq:Q-kernel}
\end{equation}

\begin{lem}
\label{lem:LJL}
Let $A_0:=\mathsf LJ\mathsf L^{\mathsf T}$. Then, in the sense of
distributional kernels,
\begin{equation}
A_0
=\begin{pmatrix}
0&I\\
-I&-(Q-Q^{\mathsf T})
\end{pmatrix}
=J-
\begin{pmatrix}
0&0\\0&Q-Q^{\mathsf T}
\end{pmatrix},
\label{eq:A0-explicit}
\end{equation}
and
\begin{equation}
A_0^{-1}
=\begin{pmatrix}
-(Q-Q^{\mathsf T})&-I\\
I&0
\end{pmatrix}.
\label{eq:A0-inverse}
\end{equation}
The normalized Fredholm Pfaffian of $A_0$ equals $1$.
\end{lem}

\begin{proof}
Write $D:=(\partial_X+t)(\partial_X+r)$. Matrix multiplication gives
\begin{equation}
\mathsf LJ\mathsf L^{\mathsf T}
=\begin{pmatrix}
0&DR^{\mathsf T}\\
-RD^{\mathsf T}&RS^{\mathsf T}-SR^{\mathsf T}
\end{pmatrix}.
\label{eq:LJL-blocks}
\end{equation}
We first prove $DR^{\mathsf T}=\Id$. Write
\begin{equation}
R^{\mathsf T}(X,Y)
=\mathbf 1_{\{X>Y\}}q(X-Y),
\qquad
q(v):=\frac{e^{-tv}-e^{-rv}}{r-t},
\quad v\geq0.
\label{eq:RT-q}
\end{equation}
Since $q(0)=0$, the first distributional derivative produces no
Dirac mass. Therefore,
\begin{align}
(\partial_X+r)R^{\mathsf T}(X,Y)
&=\mathbf 1_{\{X>Y\}}
  \bigl(q'(X-Y)+rq(X-Y)\bigr)
\notag\\
&=\mathbf 1_{\{X>Y\}}e^{-t(X-Y)}.
\label{eq:first-resolvent-derivative}
\end{align}
Applying $\partial_X+t$ once more gives
\begin{equation}
(\partial_X+t)
\bigl(\mathbf 1_{\{X>Y\}}e^{-t(X-Y)}\bigr)
=\delta(X-Y).
\label{eq:second-resolvent-derivative}
\end{equation}
Thus $DR^{\mathsf T}=\Id$, and taking algebraic transposes yields
$RD^{\mathsf T}=\Id$.

It remains to compute the $22$-entry. Since
\begin{equation}
S(X,Y)=\mathbf 1_{\{Y>X\}}
\frac{t e^{-t(Y-X)}-r e^{-r(Y-X)}}{r-t},
\label{eq:S-kernel}
\end{equation}
we have
\begin{align}
(RS^{\mathsf T})(X,Y)
={}&\frac{1}{(r-t)^2}
\int_{\max\{X,Y\}}^{\infty}
\bigl(e^{-t(Z-X)}-e^{-r(Z-X)}\bigr)
\notag\\
&\hspace{30mm}\times
\bigl(t e^{-t(Z-Y)}-r e^{-r(Z-Y)}\bigr)\,dZ.
\label{eq:RST-integral}
\end{align}

If $Y\geq X$, set $u=Z-Y$. Expanding the integrand and integrating
each exponential gives
\begin{align*}
(RS^{\mathsf T})(X,Y)
&=\frac{1}{(r-t)^2}
\int_0^\infty
\bigl(e^{-t(u+Y-X)}-e^{-r(u+Y-X)}\bigr)
\bigl(t e^{-tu}-r e^{-ru}\bigr)\,du\\
&=\frac{1}{(r-t)^2}
\left[
e^{-t(Y-X)}
\left(\frac12-\frac{r}{t+r}\right)
+
e^{-r(Y-X)}
\left(\frac12-\frac{t}{t+r}\right)
\right]\\
&=-\frac{e^{-t(Y-X)}-e^{-r(Y-X)}}{2(r^2-t^2)}
=-\frac12Q(X,Y).
\end{align*}
This calculation includes $X=Y$, where the expression vanishes.

If $X>Y$, set $u=Z-X$. The same calculation now gives
\begin{align*}
(RS^{\mathsf T})(X,Y)
&=\frac{1}{(r-t)^2}
\int_0^\infty
\bigl(e^{-tu}-e^{-ru}\bigr)
\bigl(t e^{-t(u+X-Y)}-r e^{-r(u+X-Y)}\bigr)\,du\\
&=\frac{1}{(r-t)^2}
\left[
e^{-t(X-Y)}
\left(\frac12-\frac{t}{t+r}\right)
+
e^{-r(X-Y)}
\left(\frac12-\frac{r}{t+r}\right)
\right]\\
&=\frac{e^{-t(X-Y)}-e^{-r(X-Y)}}{2(r^2-t^2)}
=\frac12Q^{\mathsf T}(X,Y).
\end{align*}
Combining the two cases yields
\begin{equation}
RS^{\mathsf T}=-\frac12(Q-Q^{\mathsf T}).
\label{eq:RST}
\end{equation}
Taking algebraic transposes, we also obtain
\begin{equation}
SR^{\mathsf T}=\frac12(Q-Q^{\mathsf T}).
\label{eq:SRT}
\end{equation}
Substituting these identities into~\eqref{eq:LJL-blocks}
proves~\eqref{eq:A0-explicit}. Direct multiplication then
gives~\eqref{eq:A0-inverse}.

Finally, set
\[
K_0:=J-A_0
=\begin{pmatrix}
0&0\\
0&Q-Q^{\mathsf T}
\end{pmatrix}.
\]
By definition, the normalized Fredholm Pfaffian of $A_0=J-K_0$ is
\[
\operatorname{Pf}(J-K_0)_{\mathbb L^2(d,\infty)}
=
1+\sum_{n=1}^{\infty}\frac{(-1)^n}{n!}
\int_{(d,\infty)^n}
\operatorname{Pf}\!\left[
K_0(X_i,X_j)
\right]_{i,j=1}^{n}
\,dX_1\cdots dX_n.
\]
For every $n\geq1$, the $2n\times2n$ skew-symmetric matrix
$\bigl[K_0(X_i,X_j)\bigr]_{i,j=1}^{n}$ has zero rows and columns
at all positions corresponding to the first component.
Its Pfaffian therefore vanishes. Every positive-order
term in the Fredholm series is consequently zero, and hence
\[
\operatorname{Pf}(J-K_0)_{\mathbb L^2(d,\infty)}=1.
\]
\end{proof}

\paragraph{Notation for the next two lemmas.}
For $Z\in\Gamma$ and $X>d$, define
\begin{equation}
e_Z(X):=e^{XZ}\begin{pmatrix}1\\-Z\end{pmatrix},
\qquad
n(Z):=(t-Z)(r-Z),
\qquad
p(Z):=(t+Z)(r+Z).
\label{eq:e-n-p}
\end{equation}
Thus $e_Z$ is vector-valued, whereas $e^{Z\,\cdot}$ denotes the scalar
function $X\mapsto e^{XZ}$. We use the bilinear pairing
\[
\langle f,g\rangle_{d}
:=\int_d^\infty f(X)^{\mathsf T}g(X)\,dX,
\]
without complex conjugation and with the transpose omitted for
scalar-valued functions.

\begin{lem}
\label{lem:kernel-conjugation}
The two matrix kernels satisfy the exact distributional operator identity
\begin{equation}
J-K^{\mathrm{gexp}}
=\mathsf L(J-K^{\mathrm{lexp}})\mathsf L^{\mathsf T}
=A_0+(\mathsf L\mathsf F)\mathsf C(\mathsf L\mathsf F)^{\mathsf T}.
\label{eq:kernel-operator-identity}
\end{equation}
Consequently,
\begin{equation}
\operatorname{Pf}(J-K^{\mathrm{gexp}})_{\mathbb L^2(d,\infty)}
=
\det\!\left(\Id+\mathsf C\widetilde{\mathsf M}_d\right)_{\mathcal H_\Gamma}^{1/2},
\qquad
\widetilde{\mathsf M}_d
:=(\mathsf L\mathsf F)^{\mathsf T}A_0^{-1}(\mathsf L\mathsf F).
\label{eq:Pfshift-det}
\end{equation}
\end{lem}

\begin{proof}
The elementary identity
\begin{equation}
P_a(e^{Z\,\cdot})(X)=\frac{e^{XZ}}{a-Z}
\label{eq:Pa-exponential}
\end{equation}
implies
\begin{equation}
R(e^{Z\,\cdot})(X)=\frac{e^{XZ}}{n(Z)},
\qquad
\partial_XR(e^{Z\,\cdot})(X)=\frac{Z e^{XZ}}{n(Z)}.
\label{eq:R-exponential}
\end{equation}
Therefore
\begin{equation}
\mathsf Le_Z(X)
=e^{XZ}
\begin{pmatrix}
-p(Z)\\[2mm]2Z/n(Z)
\end{pmatrix}.
\label{eq:LeZ}
\end{equation}
The $Z$-column of $\mathsf F$ is $\vartheta(Z)e_Z$. Substituting
\eqref{eq:theta-def} and~\eqref{eq:LeZ} into
$-(\mathsf L\mathsf F)\mathsf C(\mathsf L\mathsf F)^{\mathsf T}$ gives,
respectively,
\begin{align*}
&\frac{n(Z)n(W)}{4ZW}\frac{Z-W}{Z+W},
\qquad
-\frac{n(Z)}{2Zp(W)}\frac{Z-W}{Z+W},\\
&-\frac{n(W)}{2Wp(Z)}\frac{Z-W}{Z+W},
\qquad
\frac{1}{p(Z)p(W)}\frac{Z-W}{Z+W},
\end{align*}
multiplied by
$e^{XZ+YW}((\alpha-Z)/(\alpha+Z))^{N-1}
((\alpha-W)/(\alpha+W))^{N-1}$.
These are exactly the double-contour parts of
\eqref{eq:Kshift11}--\eqref{eq:Kshift22}. By
\eqref{eq:A0-explicit}, the remaining $22$-entry of
$J-\mathsf L(J-K^{\mathrm{lexp}})\mathsf L^{\mathsf T}$ is
$Q-Q^{\mathsf T}$, precisely the indicator term in
\eqref{eq:Kshift22}. This proves~\eqref{eq:kernel-operator-identity}
in the sense of distributional kernels.

To justify the Fredholm identity, define the integral operator
$\mathsf U:\mathcal H_\Gamma\to\mathbb L^2(d,\infty)$ by the explicit kernel
\[
\mathsf U(X,Z)
:=\vartheta(Z)e^{XZ}
\begin{pmatrix}
-p(Z)\\[2mm]2Z/n(Z)
\end{pmatrix}.
\]
Thus $\mathsf U=\mathsf L\mathsf F$ on the contour columns, and apply Lemma \ref{lem:LJL} gives
\[
J-K^{\mathrm{gexp}}=A_0+\mathsf U\mathsf C\mathsf U^{\mathsf T}.
\]

Set $\eta=-\sup_{Z\in\Gamma}\Re Z>0$, choose $0<\gamma<\eta$, and recall
\[
D_\gamma(X)=\operatorname{diag}(e^{\gamma X},e^{-\gamma X}).
\]
Define the gauged kernel and contour operator by
\[
K_\gamma(X,Y)
:=D_\gamma(X)K^{\mathrm{gexp}}(X,Y)D_\gamma(Y),
\qquad
\mathsf U_\gamma(X,Z):=D_\gamma(X)\mathsf U(X,Z).
\]
Let $E_\gamma$ denote multiplication by $e^{-\gamma X}$ on scalar
$\mathbb L^2(d,\infty)$, and set
\[
B_\gamma:=E_\gamma(Q-Q^{\mathsf T})E_\gamma,
\qquad
A_{0,\gamma}:=
\begin{pmatrix}
0&\Id\\
-\Id&-B_\gamma
\end{pmatrix}.
\]
Since $Q=(t+r)^{-1}P_tP_r$, we have
\[
E_\gamma Q E_\gamma
=\frac{1}{t+r}(E_\gamma P_t)(P_rE_\gamma).
\]
For every $a>0$, direct integration gives
\[
\|E_\gamma P_a\|_{\mathfrak S_2}^2
=\frac{e^{-2\gamma d}}{4a\gamma},
\qquad
\|P_aE_\gamma\|_{\mathfrak S_2}^2
=\frac{e^{-2\gamma d}}{4\gamma(a+\gamma)}.
\]
Thus $E_\gamma Q E_\gamma$ is a product of two Hilbert--Schmidt
operators and is trace class. Its transpose is also trace class, so
$B_\gamma\in\mathfrak S_1$.

The coefficients of $\mathsf U(X,Z)$ are bounded on $\Gamma$, and the
spatial exponential factors in the two components of $\mathsf U_\gamma$
are $e^{X(Z+\gamma)}$ and $e^{X(Z-\gamma)}$. Since $\gamma<\eta$, both
$\mathsf U$ and $\mathsf U_\gamma$ are Hilbert--Schmidt. Their algebraic
transposes are Hilbert--Schmidt as well: the contour-pairing density
with respect to arclength has constant nonzero modulus.
Also, $\|P_a\|_{\mathrm{op}}\le a^{-1}$, so $Q$ and the inverse
$A_0^{-1}$ from Lemma~\ref{lem:LJL} are bounded. In particular,
$\mathsf U^{\mathsf T}A_0^{-1}\mathsf U$ is well defined and trace class.
Because $\mathsf C$ is bounded, the kernel identity above gives
\[
J-K_\gamma=A_{0,\gamma}
+\mathsf U_\gamma\mathsf C\mathsf U_\gamma^{\mathsf T},
\qquad
K_\gamma=
\begin{pmatrix}
0&0\\
0&B_\gamma
\end{pmatrix}
-\mathsf U_\gamma\mathsf C\mathsf U_\gamma^{\mathsf T}
\in\mathfrak S_1\bigl(\mathbb L^2(d,\infty)\bigr).
\]

We also verify convergence of the Fredholm Pfaffian series. The explicit
kernels imply, with $\delta=\min\{\gamma,\eta-\gamma\}>0$,
\[
|(K_\gamma)_{ab}(X,Y)|
\le C_d e^{-\delta(X-d)}e^{-\delta(Y-d)},
\qquad X,Y>d,\quad a,b\in\{1,2\}.
\]
For $d\ge0$, one may take $C_d=Ce^{-2\delta d}$ with $C$ independent
of $d$. Hadamard's inequality therefore yields
\[
\frac{1}{m!}\int_{(d,\infty)^m}
\left|\operatorname{Pf}
[K_\gamma(X_i,X_j)]_{i,j=1}^m\right|
\prod_{i=1}^m dX_i
\le
\frac{(C_d/(2\delta))^m(2m)^{m/2}}{m!},
\]
which is summable over $m\ge1$. Moreover, finite-dimensional Pfaffian
congruence and $\det D_\gamma(X)=1$ give
\[
\operatorname{Pf}[K_\gamma(X_i,X_j)]_{i,j=1}^m
=
\left(\prod_{i=1}^m\det D_\gamma(X_i)\right)
\operatorname{Pf}[K^{\mathrm{gexp}}(X_i,X_j)]_{i,j=1}^m
=
\operatorname{Pf}[K^{\mathrm{gexp}}(X_i,X_j)]_{i,j=1}^m.
\]
Hence both Fredholm series converge absolutely and
\[
\operatorname{Pf}(J-K^{\mathrm{gexp}})_{\mathbb L^2(d,\infty)}
=
\operatorname{Pf}(J-K_\gamma)_{\mathbb L^2(d,\infty)}.
\]
This argument uses the gauge only at the kernel level, not as a bounded
similarity transformation.

The operator $A_{0,\gamma}$ has bounded inverse
\[
A_{0,\gamma}^{-1}
=
\begin{pmatrix}
-B_\gamma&-\Id\\
\Id&0
\end{pmatrix},
\qquad
J^{-1}A_{0,\gamma}
=
\begin{pmatrix}
\Id&B_\gamma\\
0&\Id
\end{pmatrix}.
\]
In particular, $J^{-1}A_{0,\gamma}-\Id$ is trace class and has square
zero, and therefore
\[
\det(J^{-1}A_{0,\gamma})_{\mathbb L^2(d,\infty)}=1.
\]
The Fredholm Pfaffian--determinant identity, multiplicativity of
Fredholm determinants, and Sylvester's identity now give
\begin{align*}
\operatorname{Pf}(J-K^{\mathrm{gexp}})_{\mathbb L^2(d,\infty)}^2
&=
\det\!\left(J^{-1}(J-K_\gamma)\right)_{\mathbb L^2(d,\infty)}\\
&=
\det(J^{-1}A_{0,\gamma})_{\mathbb L^2(d,\infty)}
\det\!\left(
\Id+A_{0,\gamma}^{-1}
\mathsf U_\gamma\mathsf C\mathsf U_\gamma^{\mathsf T}
\right)_{\mathbb L^2(d,\infty)}\\
&=
\det\!\left(
\Id+\mathsf C\mathsf U_\gamma^{\mathsf T}
A_{0,\gamma}^{-1}\mathsf U_\gamma
\right)_{\mathcal H_\Gamma}.
\end{align*}
These are ordinary Fredholm determinants: the rectangular operators
$A_{0,\gamma}^{-1}\mathsf U_\gamma$ and
$\mathsf C\mathsf U_\gamma^{\mathsf T}$ are Hilbert--Schmidt, so both
products used in Sylvester's identity are trace class.

Finally, write $\mathsf U=\binom{\mathsf U_1}{\mathsf U_2}$.
Substitution of the explicit kernels into the absolutely convergent
compression integrals cancels the gauge factors and gives
\begin{align*}
\mathsf U_\gamma^{\mathsf T}A_{0,\gamma}^{-1}\mathsf U_\gamma
&=
-\mathsf U_1^{\mathsf T}(Q-Q^{\mathsf T})\mathsf U_1
-\mathsf U_1^{\mathsf T}\mathsf U_2
+\mathsf U_2^{\mathsf T}\mathsf U_1\\
&=
\mathsf U^{\mathsf T}A_0^{-1}\mathsf U
=\widetilde{\mathsf M}_d.
\end{align*}
Consequently,
\[
\operatorname{Pf}(J-K^{\mathrm{gexp}})_{\mathbb L^2(d,\infty)}^2
=
\det\!\left(\Id+\mathsf C\widetilde{\mathsf M}_d\right)_{\mathcal H_\Gamma}.
\]
\end{proof}

\subsection{Compression to the contour space}
\label{subsec:compression-contour}

\begin{lem}
\label{lem:Mtilde-computation}
The compressed contour kernel in~\eqref{eq:Pfshift-det} satisfies
\begin{equation}
\widetilde{\mathsf M}_d
=\mathsf M_d+\frac{1}{t+r}\,\partial_d\mathsf M_d.
\label{eq:Mtilde-identity}
\end{equation}
\end{lem}

\begin{proof}
Since $\operatorname{Re}(Z+W)<0$ for $Z,W\in\Gamma$,
\begin{equation}
\langle e^{Z\,\cdot},e^{W\,\cdot}\rangle_d
=\int_d^\infty e^{X(Z+W)}\,dX
=-\frac{e^{d(Z+W)}}{Z+W}.
\label{eq:exp-pairing}
\end{equation}
Moreover, using $Q=R/(t+r)$ and
$R(e^{W\,\cdot})=n(W)^{-1}e^{W\,\cdot}$,
\begin{equation}
\bigl\langle
e^{Z\,\cdot},(Q-Q^{\mathsf T})e^{W\,\cdot}
\bigr\rangle_d
=
-\frac{e^{d(Z+W)}}{(t+r)(Z+W)}
\left(\frac{1}{n(W)}-\frac{1}{n(Z)}\right).
\label{eq:Q-pairing}
\end{equation}
Combining~\eqref{eq:A0-inverse}, \eqref{eq:LeZ},
\eqref{eq:exp-pairing}, and~\eqref{eq:Q-pairing}, we obtain
\begin{align*}
&\left\langle
\mathsf Le_Z,A_0^{-1}\mathsf Le_W
\right\rangle_d\\
&\quad=
\frac{e^{d(Z+W)}}{Z+W}
\left[
\frac{p(Z)p(W)}{t+r}
\left(\frac{1}{n(W)}-\frac{1}{n(Z)}\right)
-\frac{2Wp(Z)}{n(W)}
+\frac{2Zp(W)}{n(Z)}
\right].
\end{align*}
The bracket simplifies without a full expansion. Since
$p(Z)-n(Z)=2(t+r)Z$, and similarly for $W$,
\begin{align*}
&\frac{p(Z)p(W)}{t+r}
\left(\frac{1}{n(W)}-\frac{1}{n(Z)}\right)
-\frac{2Wp(Z)}{n(W)}
+\frac{2Zp(W)}{n(Z)}\\
&\qquad=
\frac{p(Z)}{n(W)}
\left(\frac{p(W)}{t+r}-2W\right)
+
\frac{p(W)}{n(Z)}
\left(2Z-\frac{p(Z)}{t+r}\right)\\
&\qquad=
\frac{p(Z)-p(W)}{t+r}
=
(Z-W)\left(1+\frac{Z+W}{t+r}\right).
\end{align*}
Multiplying by $\vartheta(Z)\vartheta(W)$ therefore yields
\begin{align*}
\widetilde{\mathsf M}_d(Z,W)
&=
\vartheta(Z)\vartheta(W)e^{d(Z+W)}
\frac{Z-W}{Z+W}
\left(1+\frac{Z+W}{t+r}\right)\\
&=
\mathsf M_d(Z,W)
+\frac{1}{t+r}\partial_d\mathsf M_d(Z,W),
\end{align*}
which proves~\eqref{eq:Mtilde-identity}.
\end{proof}

\subsection{A rank-two derivative identity}
\label{subsec:rank-two-identity}

The final step converts the derivative of the contour kernel into the derivative
of the Fredholm Pfaffian. 

\begin{lem}
\label{lem:rank-two-det-identity}
Let
\begin{equation}
\tau(d):=\det(\Id+\mathsf C\mathsf M_d)_{\mathcal H_{\Gamma}}^{1/2}
=\operatorname{Pf}(J-K^{\mathrm{lexp}})_{\mathbb L^2(d,\infty)}.
\label{eq:tau-def}
\end{equation}
Then
\begin{equation}
\det\!\left(\Id+\mathsf C\mathsf M_d
+\frac{1}{t+r}\mathsf C\partial_d\mathsf M_d\right)_{\mathcal H_{\Gamma}}^{1/2}
=\tau(d)+\frac{1}{t+r}\tau'(d).
\label{eq:rank-two-final}
\end{equation}
\end{lem}

\begin{proof}
First, all contour determinants in the statement are well defined. Indeed,
$\mathsf C$ and $\mathsf M_d$ are Hilbert--Schmidt, so
$\mathsf C\mathsf M_d$ is trace class. Furthermore,
\begin{equation}
\partial_d\mathsf M_d(Z,W)
=(Z-W)\vartheta(Z)\vartheta(W)e^{d(Z+W)}
=a_d(Z)b_d(W)-b_d(Z)a_d(W),
\label{eq:Md-rank-two}
\end{equation}
where
\begin{equation}
a_d(Z):=Z\vartheta(Z)e^{dZ},
\qquad
b_d(Z):=\vartheta(Z)e^{dZ}.
\label{eq:ad-bd}
\end{equation}
Thus $\partial_d\mathsf M_d$ has rank at most two, and
$\mathsf C\partial_d\mathsf M_d$ is trace class.

We first work for $d$ sufficiently large. Since $\Gamma$ is a compact
subset of the left half-plane, choose $\eta>0$ such that
$\Re Z\leq-\eta$ on $\Gamma$. From~\eqref{eq:C-def}
and~\eqref{eq:Md-explicit},
\[
|\mathsf M_d(Z,W)|
=|\vartheta(Z)\vartheta(W)e^{d(Z+W)}\mathsf C(Z,W)|
\leq \|\vartheta\|_{L^\infty(\Gamma)}^2 e^{-2\eta d}
|\mathsf C(Z,W)|.
\]
Consequently,
\begin{equation}
\|\mathsf C\mathsf M_d\|_{\mathfrak S_1}
\leq \|\mathsf C\|_{\mathfrak S_2}
     \|\mathsf M_d\|_{\mathfrak S_2}
\leq \|\mathsf C\|_{\mathfrak S_2}^2
     \|\vartheta\|_{L^\infty(\Gamma)}^2 e^{-2\eta d}.
\label{eq:CM-trace-bound}
\end{equation}

Thus $\Id+\mathsf C\mathsf M_d$ is invertible for all sufficiently large $d$ by
the Neumann series. Define
\begin{equation}
\mathsf G_d:=(\Id+\mathsf C\mathsf M_d)^{-1}\mathsf C.
\label{eq:Gd-def}
\end{equation}

We next prove that $\mathsf G_d$ is skew-symmetric. The push-through identity
\begin{equation}
(\Id+\mathsf C\mathsf M_d)^{-1}\mathsf C
=\mathsf C(\Id+\mathsf M_d\mathsf C)^{-1}
\label{eq:push-through}
\end{equation}
follows from
$(\Id+\mathsf C\mathsf M_d)\mathsf C
=\mathsf C(\Id+\mathsf M_d\mathsf C)$. Since
$\mathsf C^{\mathsf T}=-\mathsf C$ and
$\mathsf M_d^{\mathsf T}=-\mathsf M_d$,
\begin{align*}
\mathsf G_d^{\mathsf T}
&=\mathsf C^{\mathsf T}
\bigl((\Id+\mathsf C\mathsf M_d)^{\mathsf T}\bigr)^{-1}\\
&=-\mathsf C(\Id+\mathsf M_d\mathsf C)^{-1}
=-\mathsf G_d.
\end{align*}

Let
\begin{equation}
\mathsf V_d:=(a_d\ \ b_d):\mathbb C^2\longrightarrow\mathcal H_{\Gamma},
\qquad
 J=\begin{pmatrix}0&1\\-1&0\end{pmatrix}.
\label{eq:U-Omega}
\end{equation}
Then~\eqref{eq:Md-rank-two} gives
\begin{equation}
\partial_d\mathsf M_d=\mathsf V_d J\mathsf V_d^{\mathsf T}.
\label{eq:Mprime-U}
\end{equation}
Since $\mathsf G_d$ is skew-symmetric, the $2\times2$ matrix
$\mathsf V_d^{\mathsf T}\mathsf G_d\mathsf V_d$ is skew-symmetric. Thus there
is a scalar
\begin{equation}
\varrho(d):=a_d^{\mathsf T}\mathsf G_db_d
\label{eq:rho-def}
\end{equation}
such that
\begin{equation}
\mathsf V_d^{\mathsf T}\mathsf G_d\mathsf V_d
=\varrho(d) J.
\label{eq:UGU}
\end{equation}

Let $\lambda:=(t+r)^{-1}$. Factoring $\Id+\mathsf C\mathsf M_d$ and applying
Sylvester's determinant identity gives
\begin{align}
&\det\!\left(\Id+\mathsf C\mathsf M_d
+\lambda\mathsf C\partial_d\mathsf M_d\right)
\notag\\
&\quad=\det(\Id+\mathsf C\mathsf M_d)
\det\!\left(\Id_2+\lambda J
\mathsf V_d^{\mathsf T}\mathsf G_d\mathsf V_d\right)
\notag\\
&\quad=\det(\Id+\mathsf C\mathsf M_d)
\det(\Id_2+\lambda\varrho(d) J^2)
\notag\\
&\quad=\det(\Id+\mathsf C\mathsf M_d)(1-\lambda\varrho(d))^2.
\label{eq:det-rank-two}
\end{align}
On the other hand, differentiation of the Fredholm determinant yields
\begin{align}
\frac{\tau'(d)}{\tau(d)}
&=\frac12\operatorname{Tr}
\left((\Id+\mathsf C\mathsf M_d)^{-1}
\mathsf C\partial_d\mathsf M_d\right)
\notag\\
&=\frac12\operatorname{tr}_{\mathbb C^2}
\left( J\mathsf V_d^{\mathsf T}\mathsf G_d\mathsf V_d\right)
=\frac12\varrho(d)\operatorname{tr}( J^2)
=-\varrho(d).
\label{eq:tau-log-derivative}
\end{align}
For large $d$, all determinants tend to $1$, so the square-root branch in
\eqref{eq:det-rank-two} is the one satisfying
\begin{equation}
\det\!\left(\Id+\mathsf C\mathsf M_d
+\lambda\mathsf C\partial_d\mathsf M_d\right)^{1/2}
=\tau(d)(1-\lambda\varrho(d)).
\label{eq:sqrt-branch}
\end{equation}
Using~\eqref{eq:tau-log-derivative}, the right-hand side equals
$\tau(d)+\lambda\tau'(d)$, which proves~\eqref{eq:rank-two-final} for all
sufficiently large $d$.

To remove the restriction on $d$, translate $X=d+x$ and $Y=d+y$. On the fixed
space $\mathbb L^2((0,\infty))\otimes\mathbb C^2$, the Fredholm Pfaffians on both sides
of~\eqref{eq:rank-two-final} are analytic functions of the complex parameter
$d$: the translated kernels depend on $d$ through factors $e^{d(Z+W)}$, and
the conjugated Fredholm series converge locally uniformly. Hence the identity,
proved on a nonempty open half-plane of large $d$, extends to every $d>0$ by
the identity theorem.
\end{proof}

\begin{proof}[Proof of Proposition~\ref{prop:exact-matching-finite-N}]
Combining Lemmas~\ref{lem:kernel-conjugation} and
\ref{lem:Mtilde-computation} gives
\begin{align*}
\operatorname{Pf}(J-K^{\mathrm{gexp}})_{\mathbb L^2(d,\infty)}
&=\det\!\left(\Id+\mathsf C\mathsf M_d
+\frac{1}{t+r}\mathsf C\partial_d\mathsf M_d\right)^{1/2}.
\end{align*}
Lemma~\ref{lem:rank-two-det-identity} then yields
\begin{align*}
\operatorname{Pf}(J-K^{\mathrm{gexp}})_{\mathbb L^2(d,\infty)}
&=\operatorname{Pf}(J-K^{\mathrm{lexp}})_{\mathbb L^2(d,\infty)}
+\frac{1}{t+r}\partial_d
\operatorname{Pf}(J-K^{\mathrm{lexp}})_{\mathbb L^2(d,\infty)},
\end{align*}
which is exactly~\eqref{eq:exact-matching-main}.
\end{proof}

\begin{rmk}
The matching is established at finite $N$, before any stationarity limit or KPZ scaling is taken. The result follows from
the exact operator identity~\eqref{eq:kernel-operator-identity}, the contour
compression~\eqref{eq:Mtilde-identity}, and the rank-two derivative identity
\eqref{eq:rank-two-final}.
\end{rmk}

\section{Matching the high-density and low-density formulas at the boundary}
\label{sec:boundary-matching}

In this section, we prove that the high-density formula \eqref{paper:hd-limit} converges to the
low-density formula \eqref{paper:ld-limit} at both $\rr=-\s$ and $\rr=\s$, for fixed
$\s<0$. We retain the notation of Theorems~\ref{thm: high-density phase theorem} and~\ref{thm: low-density phase theorem}:
$\widetilde K$ is the two-label high-density kernel in \eqref{paper:hd-kernel}--\eqref{paper:hd-residue}, and
$\tKK$ is the low-density kernel in \eqref{paper:ld-kernel}--\eqref{paper:ld-residue}.

For $\tilde{d}\in\R$ and $\tilde{y}\ge0$, let
\[
\widetilde T=\{1\}\times(\tilde{d},\infty)\ \cup\ \{2\}\times(\tilde{d}+\tilde{y},\infty).
\]
We denote the right-hand sides of \eqref{paper:hd-limit} and \eqref{paper:ld-limit}, respectively, by
\begin{align*}
P^{\mathrm{HD}}_{\rr,\s}(\tilde{d})
&:=-(\rr+\s)\int_0^\infty e^{(\rr+\s)\tilde{y}}
\Pf(J-\widetilde K)_{\mathbb{L}^2(\widetilde T)}\,\mathrm{d}\tilde{y},\\
F^{\mathrm{LD}}_{\s}(\tilde{d})
&:=\left.\Pf(J-\tKK)_{\mathbb{L}^2(\tilde{d},\infty)}\right|_{\rr=\s}.
\end{align*}
The low-density parameter in the second expression is $\s$. 

\begin{prop}\label{boundary:prop:matching}
Fix $\s<0$. For every $\tilde{d}\in\R$,
\begin{equation}\label{boundary:eq:claim}
\lim_{\rr\uparrow-\s}P^{\mathrm{HD}}_{\rr,\s}(\tilde{d})
=\lim_{\rr\to\s}P^{\mathrm{HD}}_{\rr,\s}(\tilde{d})
=F^{\mathrm{LD}}_{\s}(\tilde{d}).
\end{equation}
The second limit is two-sided. The entries of $\widetilde K$ have removable
singularities at $\rr=\pm\s$, and
$P^{\mathrm{HD}}_{\rr,\s}(\tilde{d})$ extends holomorphically to a complex
neighborhood of $\rr=\s$.
\end{prop}

\subsection{Analyticity and uniform Fredholm estimates}

\begin{lem}[Removal of the apparent poles]\label{boundary:lem:removable}
For every $n,n'\in\{1,2\}$ and fixed $X,Y$, every entry of
$\widetilde K(n,X;n',Y)$ extends holomorphically in $\rr$ to neighborhoods of
$\s$ and $-\s$. The extensions are obtained by combining
$\widetilde G+\widetilde O+\widetilde R$ before taking the limit.
\end{lem}

\begin{proof}
First suppose $X\le Y$. For the following calculation only, set
\begin{equation}\label{boundary:eq:numerator}
f_{nn'}(\rr)=2\s(\s-\rr)^2(\s+\rr)^2
(\widetilde G+\widetilde O+\widetilde R)(n,X;n',Y).
\end{equation}
Multiplying out \eqref{eq: definitino of G tilde}--\eqref{paper:hd-residue} clears every denominator. Thus $f_{nn'}$ is a
finite sum of polynomials times exponentials, and is entire in the complex
variable $\rr$. For each fixed label pair, $f_{nn'}$ is scalar-valued. 

For the four label pairs, the cancellations are as follows. Substitution
in the denominator-free expression gives
\[
\bigl(f_{nn'}(\s)\bigr)_{n,n'=1}^2
=\bigl(f_{nn'}(-\s)\bigr)_{n,n'=1}^2=0.
\]
Differentiation, followed by substitution, gives
\begin{align}
\bigl(f'_{nn'}(\s)\bigr)_{n,n'=1}^2
&=\bigl[(Y-X)+(X-\s^2)-(Y-\s^2)\bigr]
\begin{pmatrix}
0&0\\
4\s^2e^{\s(Y-X)}&2\s e^{\s(Y-X)}
\end{pmatrix}=0,\label{boundary:eq:derivplus}\\
\bigl(f'_{nn'}(-\s)\bigr)_{n,n'=1}^2
&=\bigl[(Y-X)+(X-\s^2)-(Y-\s^2)\bigr]
\begin{pmatrix}
0&4\s^2e^{\s(X-Y)}\\
0&2\s e^{\s(X-Y)}
\end{pmatrix}=0.\label{boundary:eq:derivminus}
\end{align}
These identities follow from
$\partial_{\rr}e^{-\rr^3/3+X\rr}=(X-\rr^2)e^{-\rr^3/3+X\rr}$.
For example, the nonzero terms before cancellation in $f'_{12}(-\s)$ are
\[
4\s^2e^{\s(X-Y)}
\bigl[-(Y-\s^2)+(Y-X)+(X-\s^2)\bigr].
\]
For the label pair $(1,1)$,
\begin{align*}
f_{11}(\rr)=2\s(\s^2-\rr^2)\Bigl[
&e^{\rr(Y-X)}-e^{\s(Y-X)}\\
&-e^{-(\rr^3+\s^3)/3+Y\rr+X\s}
+e^{-(\rr^3+\s^3)/3+X\rr+Y\s}\Bigr].
\end{align*}
The bracket vanishes at both boundaries, so both zeros are double.

For $\rr_0\in\{\s,-\s\}$, Taylor's formula along the line segment from
$\rr_0$ to $\rr$ now gives the explicit extension
\begin{equation}\label{boundary:eq:taylor}
(\widetilde G+\widetilde O+\widetilde R)(n,X;n',Y)
=\frac{\displaystyle\int_0^1(1-t)
f''_{nn'}\bigl(\rr_0+t(\rr-\rr_0)\bigr)\,\mathrm{d} t}
{2\s(\rr+\rr_0)^2}.
\end{equation}
The denominator is nonzero in a sufficiently small disk about $\rr_0$,
since $\s\ne0$. When $\rr = \rr_0$ this equals
$f''_{nn'}(\rr_0)/(16\s\rr_0^2)$.
The formula is holomorphic and agrees with the original sum off the center.
The skew-symmetric extension gives the result for $X>Y$.

For the contour terms, choose a small complex disk about each boundary
$\rr=\pm\s$. The contours can be fixed throughout this disk, with all
poles staying a positive distance from the contours. Thus the integrands
are holomorphic in $\rr$, even when poles merge away from the contours.

The integrands decay exponentially along the
contours and uniformly in $\rr$. Differentiating in $\rr$ contributes only
rational factors, whose poles still stay away from the contours.
The derivatives therefore also have an integrable bound independent
of $\rr$ in the disk. By dominated convergence, we may take the complex
derivative in $\rr$ inside each contour integral. Hence every contour
integral is holomorphic in $\rr$, which completes the proof.
\end{proof}

\begin{lem}[Uniform bounds and regularity of the Fredholm Pfaffians]
\label{boundary:lem:series}
Fix a lower spatial bound $\tilde d_0$ and a sufficiently small closed
complex disk $D$ about either boundary $\tilde{r} = \s$ or $\rr = -\s$.
Then there exist $C,\theta>0$ such that, uniformly for $\rr\in D$,
$m\ge1$, $n_j\in\{1,2\}$, and $X_j\ge\tilde d_0$,
\begin{equation}\label{boundary:eq:pfbound}
\left|\Pf\bigl[\widetilde K(n_j,X_j;n_k,X_k)\bigr]_{j,k=1}^m\right|
\le C^m(2m)^{m/2}\prod_{j=1}^m e^{-\theta X_j}.
\end{equation}
\end{lem}

\begin{proof} Lemma \ref{boundary:lem:removable} shows the point-wise limit of $\widetilde{K}_{i,j}$ for $i,j \in \{1,2\}$.
Then uniform finite-$N$ Pfaffian bound established in the proof of
Lemma~\ref{lem: asymp pfaffian} passes to the limit by entry-wise
convergence and by the fact that taking $\rr = \pm \s$ does not affect the exponential decay coming from $e^{- X\zeta + Y\omega}$. Thus, one can
choose $B>\max(0,-\s)$ and $\kappa > B$ such that for all $\rr \in D$
\begin{align}
|\widetilde K_{11}(n,X;n',Y)|&\le C e^{-\kappa (X+Y)},\nonumber\\
|\text{double integral in }\widetilde K_{12}(n,X;n',Y)|
&\le C e^{-\kappa X+BY},\label{boundary:eq:rawbounds}\\
|\widetilde K_{22}(n,X;n',Y)|&\le C e^{B(X+Y)}.\nonumber
\end{align} Applying the same conjugation argument in Lemma~\ref{lem: asymp pfaffian} shows that there exists $C,\theta > 0$ such that \eqref{boundary:eq:pfbound} holds.
\end{proof}

\begin{lem}[The second-label estimate]\label{boundary:lem:B}
Uniformly for $\rr$ near either boundary, for $\tilde{d}$ in a compact interval,
and for every $\tilde{y}\ge0$,
\begin{equation}\label{boundary:eq:B}
\left|\Pf(J-\widetilde K)_{\mathbb{L}^2(\widetilde T)}
-\Pf\!\bigl(J-\widetilde K(1,\cdot;1,\cdot)\bigr)_{\mathbb{L}^2(\tilde{d},\infty)}\right|
\le C e^{-\theta \tilde{y}},
\end{equation}
where $\theta$ is chosen in Lemma \ref{boundary:lem:series} and $C >0$ is independent of $\tilde{y}$.
\end{lem}

\begin{proof}
The order-$m$ term of the two-label series is
\[
\frac{(-1)^m}{m!}\sum_{n_1,\ldots,n_m=1}^2
\int_{\tilde{d}+\tilde{y}\Id_{n_1=2}}^\infty\!\cdots
\int_{\tilde{d}+\tilde{y}\Id_{n_m=2}}^\infty
\Pf\bigl[\widetilde K(n_j,X_j;n_k,X_k)\bigr]_{j,k=1}^m
\prod_{j=1}^m\,\mathrm{d} X_j.
\]
Absolute convergence permits subtraction term by term. The $m=0$ terms
and every all-label-$1$ term cancel exactly. For a remaining label
sequence containing $k\ge1$ second labels, \eqref{boundary:eq:pfbound} gives
\begin{equation}\label{boundary:eq:klabels}
\int\cdots\int\prod_{j=1}^m e^{-\theta X_j}\,\mathrm{d} X_j
=\bigl(\tfrac1 \theta e^{-\theta \tilde{d}}\bigr)^m e^{-\theta k\tilde{y}}.
\end{equation}
There are $\binom{m}{k}$ such sequences: choose the $k$ positions
occupied by label $2$; every other position must have label $1$. Hence the left-hand side of
\eqref{boundary:eq:B} is at most
\[
\sum_{m=1}^\infty\frac{C^m(2m)^{m/2}}{m!}
\bigl(\tfrac1 \theta e^{-\theta\tilde{d}}\bigr)^m
\left((1+e^{-\theta\tilde{y}})^m-1\right).\]

For $0\le x\le1$, the mean value theorem gives
$(1+x)^m-1\le m2^{m-1}x$. Thus the last display is bounded by
\begin{equation}\label{boundary:eq:Bsum}
\frac{e^{-\theta \tilde{y}}}2
\sum_{m=1}^\infty
\frac{(2C/\theta)^m(2m)^{m/2}}{(m-1)!}e^{-\theta m\tilde{d}}
\le C e^{-\theta \tilde{y}}.
\end{equation}
\end{proof}

\subsection{Matching at \texorpdfstring{$\rr=-\s$}{r = -s}}

\begin{lem}[Identification of the first-label kernel]\label{boundary:lem:identify}
With the analytic extension of Lemma~\ref{boundary:lem:removable},
\begin{equation}\label{boundary:eq:identify}
\left.\widetilde K(1,X;1,Y)\right|_{\rr=-\s}
=\left.\tKK(X,Y)\right|_{\rr=\s}.
\end{equation}
\end{lem}

\begin{proof}
For the $11$ and $12$ entries, \eqref{paper:hd-ab} gives, at $\rr=-\s$,
\begin{align*}
-\A(1,\zeta)\B(1,\omega)(\zeta-\rr)(\omega+\rr)
&=(\zeta^2-\s^2)(\omega^2-\s^2),\\
\frac{\A(1,\zeta)(\zeta-\rr)}{\A(1,\omega)(\omega-\rr)}
&=\frac{\zeta^2-\s^2}{\omega^2-\s^2}.
\end{align*}
There is no propagation term for $n=n'=1$. These are exactly the
integrands in \eqref{paper:ld-kernel}, with low-density parameter $\s$; the contours agree after
irrelevant pole constraints are deleted. Transposition gives the $21$ entry.

The integrand \eqref{eq: integrand of 22 limit} of the $22$ entry becomes
\[
-e^{\zeta^3/3-\omega^3/3-X\zeta+Y\omega}
\frac{\zeta+\omega}
{(\zeta^2-\s^2)(\omega^2-\s^2)(\zeta-\omega)}.
\]
The first contour pair agrees with the first pair in \eqref{paper:ld-residue}. The second
$\omega$ contour is closed and surrounds the poles at $\s,-\s$ once,
as does the contour in \eqref{paper:ld-residue}. The coalescing pole labels do
not create two copies of the contour. The corresponding $\zeta$ contours
can be chosen identically to the right of this loop and all their poles.

Finally, $\widetilde O(1,X;1,Y)=0$. To verify the elementary part explicitly,
use the bracket in the proof of Lemma~\ref{boundary:lem:removable}. Its derivative
at $\rr=-\s$ is
\[
(X-\s^2)\bigl[e^{\s(Y-X)}-e^{\s(X-Y)}\bigr].
\]
The derivative of $\s^2-\rr^2$ there is $2\s$. Thus, for $X\le Y$,
\[
\left.(\widetilde G+\widetilde R)(1,X;1,Y)\right|_{\rr=-\s}
=\frac{X-\s^2}{2\s}
\bigl[e^{\s(Y-X)}-e^{\s(X-Y)}\bigr].
\]
This is precisely $\widetilde{\mathsf{F}}$ in \eqref{paper:ld-residue}, with low-density parameter
$\s$. For $X>Y$ use antisymmetry. This proves the whole matrix identity.
\end{proof}

\begin{proof}[Proof of the first limit in \eqref{boundary:eq:claim}]
Put $\rr=-\s-\varepsilon$, $\varepsilon>0$. Then
\[
P^{\mathrm{HD}}_{-\s-\varepsilon,\s}(\tilde{d})=\varepsilon\int_0^\infty e^{-\varepsilon\tilde{y}}
\Pf(J-\widetilde K)_{\mathbb{L}^2(\widetilde T)}\,\mathrm{d}\tilde{y}.
\]
By \eqref{boundary:eq:B}, replacing this Pfaffian by its first-label series has
error at most
\[
C\varepsilon\int_0^\infty e^{-(\varepsilon+\theta)\tilde{y}}\,\mathrm{d}\tilde{y}
=\frac{C\varepsilon}{\varepsilon+\theta}\longrightarrow0
\]
for some fixed $\theta>0.$
The first-label series is independent of $\tilde{y}$, and the averaging density
has integral one. Lemma~\ref{boundary:lem:series} passes the parameter limit through
that series, and Lemma~\ref{boundary:lem:identify} identifies its limit with \eqref{paper:ld-limit}.
This proves the first limit in \eqref{boundary:eq:claim}.
\end{proof}

The same lemmas show that the low-density Fredholm series at parameter $\s$ is
absolutely convergent and continuous in $\tilde{d}$. This continuity is used
in the finite-size coupling comparison below.

\subsection{Matching at \texorpdfstring{$\rr=\s$}{r = s}}
We put subscripts on $\widetilde K$ to indicate its parameter values.
\begin{lem}[Explicit two-label limit at $\rr=\s$]\label{boundary:lem:equal-kernel}
Fix $\s<0$. The boundary kernel $\widetilde K_{\s,\s}$ is given explicitly
below, for every $n,n'\in\{1,2\}$. For every $\tilde{y}\ge0$,
\begin{equation}\label{eq: s = r limit}
\lim_{\rr\to\s}\Pf(J-\widetilde K_{\rr,\s})_{\mathbb{L}^2(\widetilde T)}
=\Pf(J-\widetilde K_{\s,\s})_{\mathbb{L}^2(\widetilde T)}.
\end{equation}
The convergence is uniform for $\tilde{y}\ge0$ and $\tilde{d}$ in a fixed compact
interval.
\end{lem}

\begin{proof}
We compute the contour terms and the elementary terms separately, then
pass their sum through the Fredholm series.

For contour terms,
\[
\A(n,z)=(z-\s)(z+\s)^{n-1},\qquad
\B(n,z)=(-1)^n(z+\s)(z-\s)^{n-1}.
\]
For each entry, direct substitution gives
\begin{align*}
\widetilde K_{\s,\s;11}(n,X;n',Y)
&=-\int\frac{\,\mathrm{d}\omega}{2\pi\I}
\int\frac{\,\mathrm{d}\zeta}{2\pi\I}
e^{\zeta^3/3-\omega^3/3-X\zeta+Y\omega}\A(n,\zeta)\B(n',\omega)
\frac{(\zeta+\omega)(\zeta-\s)(\omega+\s)}
{4\zeta\omega(\zeta-\omega)},\\
\widetilde K_{\s,\s;12}(n,X;n',Y)
&=-\int\frac{\,\mathrm{d}\omega}{2\pi\I}
\int\frac{\,\mathrm{d}\zeta}{2\pi\I}
e^{\zeta^3/3-\omega^3/3-X\zeta+Y\omega}
\frac{\A(n,\zeta)}{\A(n',\omega)}
\frac{(\zeta+\omega)(\zeta-\s)}
{2\zeta(\omega-\s)(\zeta-\omega)}\\
&\quad-\Id_{n=1,n'=2,X\le Y}e^{-\s(Y-X)},\\
\widetilde K_{\s,\s;21}(n,X;n',Y)
&=-\widetilde K_{\s,\s;12}(n',Y;n,X).
\end{align*}

For the $22$ component, keep both double integrals in \eqref{paper:hd-h22} and replace their integrand \eqref{eq: integrand of 22 limit} by
\[
\left.\widetilde F(\zeta,\omega)\right|_{\rr=\s}
=e^{\zeta^3/3-\omega^3/3-X\zeta+Y\omega}
\frac{\zeta+\omega}
{\A(n',\omega)\B(n,\zeta)(\zeta-\omega)(\zeta+\s)(\omega-\s)}.
\]
These formulas specify every contour entry for all four spatial-label pairs.

As $\rr\to\s$, no pole crosses the contours. For fixed $X,Y$, the
cubic decay of $e^{\zeta^3/3-\omega^3/3}$ along the Airy contours gives
an integrable bound uniform in $\rr$ near $\s$. Dominated convergence
therefore permits taking the limit inside the double contour integrals.

\emph{The elementary terms.} The remaining calculation is the finite
part of $\widetilde G+\widetilde O+\widetilde R$. For $X\le Y$, we have
\begin{align*}
\lim_{\rr\to\s}(\widetilde G+\widetilde O+\widetilde R)(1,X;1,Y)
&=\frac{Y-X}{2\s}
\left[e^{-2\s^3/3+\s(X+Y)}-e^{\s(Y-X)}\right],\\[4pt]
\lim_{\rr\to\s}(\widetilde G+\widetilde O+\widetilde R)(1,X;2,Y)
&=\frac{2\s(Y-X)-1}{8\s^3}
\left[e^{-2\s^3/3+\s(X+Y)}-e^{\s(Y-X)}\right]\\
&\quad+\left[\frac{(X-\s^2)^2}{4\s}
+\frac{X-\s^2}{4\s^2}-\frac12\right]e^{-\s(Y-X)},\\[4pt]
\lim_{\rr\to\s}(\widetilde G+\widetilde O+\widetilde R)(2,X;1,Y)
&=\frac{2\s(Y-X)+1}{8\s^3}
\left[e^{-2\s^3/3+\s(X+Y)}-e^{\s(Y-X)}\right]\\
&\quad+\left[\frac12-\frac{(X-\s^2)(2Y-X-\s^2)}{4\s}
-\frac{X-\s^2}{4\s^2}\right]e^{\s(Y-X)},\\[4pt]
\lim_{\rr\to\s}(\widetilde G+\widetilde O+\widetilde R)(2,X;2,Y)
&=\frac{Y-X}{8\s^3}
\left[e^{-2\s^3/3+\s(X+Y)}-e^{\s(Y-X)}\right]\\
&\quad+\frac{(X-\s^2)^2-2\s}{8\s^2}
\left[e^{-\s(Y-X)}-e^{\s(Y-X)}\right]\\
&\quad-\frac{(Y-X)(X-\s^2)}{4\s^2}e^{\s(Y-X)}.
\end{align*}
These finite parts follow by Taylor expansion. Use the scalar numerator $f_{nn'}$ already
defined in \eqref{boundary:eq:numerator}, and put $\delta=\rr-\s$ for this
Taylor calculation only. The needed expansions are
\begin{align*}
e^{\rr(Y-X)}
&=e^{\s(Y-X)}
\left[1+(Y-X)\delta+\frac{(Y-X)^2}{2}\delta^2\right]+O(\delta^3),\\
e^{-\rr^3/3+X\rr}
&=e^{-\s^3/3+X\s}
\left[1+(X-\s^2)\delta+
\left(\frac{(X-\s^2)^2}{2}-\s\right)\delta^2\right]+O(\delta^3).
\end{align*}
The same second expansion with $X$ replaced by $Y$ treats the remaining
exponentials in \eqref{eq: definitino of G tilde}--\eqref{paper:hd-residue}. Substitute them into the denominator-free
numerators $f_{nn'}$. Their constant and linear terms cancel for every
label pair. For example, the coefficient of $\delta^2$ in the mixed
numerator $f_{12}(\s+\delta)$ is
\begin{align*}
&(1-2\s(Y-X))
\left[e^{\s(Y-X)}-e^{-2\s^3/3+\s(X+Y)}\right]\\
&\qquad+\left[2\s^2(X-\s^2)^2
+2\s(X-\s^2)-4\s^3\right]e^{-\s(Y-X)}.
\end{align*}
Since
\[
2\s(\s-\rr)^2(\s+\rr)^2
=8\s^3\delta^2+O(\delta^3),
\]
division by $8\s^3$ gives exactly the displayed $(1,2)$ value.
The coefficients of $\delta^2$ in the other three numerators, divided by
the same $8\s^3$, give the other three displayed values. Equivalently,
the calculation evaluates every $f''_{nn'}(\s)/(16\s^3)$ explicitly.
The Taylor expansions hold for fixed $X,Y$, while Lemma~\ref{boundary:lem:series}
provides the uniform spatial bound.

Finally, the uniform kernel entry bounds and the Fredholm estimate in Lemma~\ref{boundary:lem:series} together with the dominated convergence theorem
proves \eqref{eq: s = r limit}. Uniformity for $\tilde{y}\ge0$ and for $\tilde{d}$
in a compact interval follows by bounding the integral of the absolute difference of Pfaffians over the common domain $X_j\ge\tilde{d}_0$.
\end{proof}

\begin{lem}[Analyticity of the high-density average at $\rr=\s$]\label{boundary:lem:average}
$P^{\mathrm{HD}}_{\rr,\s}(\tilde{d})$ extends holomorphically in $\rr$ to a
complex neighborhood of $\s$. Its value at $\s$ is
\begin{equation}\label{boundary:eq:average}
P^{\mathrm{HD}}_{\s,\s}(\tilde{d})=(-2\s)\int_0^\infty e^{2\s\tilde{y}}
\left.\Pf(J-\widetilde K)_{\mathbb{L}^2(\widetilde T)}\right|_{\rr=\s}\,\mathrm{d}\tilde{y}.
\end{equation}
\end{lem}

\begin{proof}
By Lemma~\ref{boundary:lem:equal-kernel},
\[
\Pf(J-\widetilde K_{\rr,\s})_{\mathbb{L}^2(\widetilde T)}
\longrightarrow\Pf(J-\widetilde K_{\s,\s})_{\mathbb{L}^2(\widetilde T)}
\]
uniformly for $\tilde{y}\ge0$ and locally uniformly in $\tilde{d}$, with the boundary
kernel given explicitly in that lemma.

Choose $|\rr-\s|<|\s|/2$. Then $\Re(\rr+\s)\le\s<0$ and the full
integrand is bounded by $Ce^{\s\tilde{y}}$. Consequently,
\[
\lim_{\rr\to\s}P^{\mathrm{HD}}_{\rr,\s}(\tilde{d})
=(-2\s)\int_0^\infty e^{2\s\tilde{y}}
 \Pf(J-\widetilde K_{\s,\s})_{\mathbb{L}^2(\widetilde T)}\,\mathrm{d}\tilde{y}.
\]
Dominated convergence proves this identity and joint continuity in
$\rr,\tilde{d}$. The local uniform Fredholm bounds and entrywise holomorphy
also give holomorphy of the average by Morera's theorem and Fubini.
\end{proof}

We now identify the boundary value \eqref{boundary:eq:average} by comparing
$\rr=\s-\varepsilon$ and $\rr=\s+\varepsilon$. Here $r,s$ are the finite-size parameters in Definition~\ref{def: stationary measure}. Write $\mathbf{G}_{r,s}$
for the solution of \eqref{recurrence}--\eqref{weights} with initial condition
\[
\mathbf{G}_{r,s}(1,1)=0,\qquad \mathbf{G}_{r,s}(j+1,1)=\mathcal{I}_{r,s}(j),\quad j\ge0.
\]
We write $Y_{r,s}$ for the independent $\mathrm{Geom}(rs)$ variable in
Definition~\ref{def: stationary measure} and recall the scale $\rho=(c_0N/2)^{1/3}$ from \eqref{eq: scale}.

\begin{lem}[Finite-$N$ comparison with a uniform error]\label{boundary:lem:coupling}
Fix $0<\varepsilon<\min(1,|\s|/4)$ and put
\begin{equation}\label{boundary:eq:parameters}
s=1+\frac{\s}{\rho},\qquad
r_-=1+\frac{\s-\varepsilon}{\rho},\qquad
r_+=1+\frac{\s+\varepsilon}{\rho}.
\end{equation}
For all sufficiently large $N$, the three stationary models admit a
coupling such that, at the common endpoint $(N+2,N+2)$,
\begin{equation}\label{boundary:eq:pathwise}
\begin{aligned}
\mathbf{G}_{r_-,s}(N+2,N+2)-(Y_{s,s}-Y_{r_-,s})
&\le \mathbf{G}_{s,s}(N+2,N+2)\\
&\le \mathbf{G}_{r_+,s}(N+2,N+2)+(Y_{r_+,s}-Y_{s,s}),
\end{aligned}
\end{equation}
and
\begin{equation}\label{boundary:eq:couplingerror}
\Pb\bigl(|Y_{r_\pm,s}-Y_{s,s}|>\rho\sqrt\varepsilon\bigr)
\le C\sqrt\varepsilon.
\end{equation}
Here $C$ and the lower bound on $N$ can be chosen independently of
$\varepsilon$ in the stated interval.
\end{lem}

\begin{proof}
In Definition~\ref{def: stationary measure}, the increments of $R_1$ have law $\mathrm{Geom}(\sqrt q/s)$.
The increments of $R_2$ have law $\mathrm{Geom}(\sqrt q\,s)$. Each walk has
independent increments, and the two walks are independent. Their laws
depend on $s$ alone, so use the same walks in all three models.
Using a common uniform random variable, independent of the
walks, and inverse distribution functions, couple
\[
Y_{r_-,s}\le Y_{s,s}\le Y_{r_+,s},\qquad
Y_{r,s}\sim\mathrm{Geom}(rs).
\]
The initial profile is explicitly
\[
\mathcal{I}_{r,s}(x)=R_2(x)+\max\left\{0,\max_{1\le k\le x}
\bigl(R_1(k)-R_2(k-1)\bigr)-Y_{r,s}\right\},\qquad x\ge1.
\]
For fixed walks, this is a decreasing, $1$-Lipschitz function of
$Y_{r,s}$, simultaneously for every $x$. Therefore
\begin{align*}
0\le \mathcal{I}_{r_-,s}(x)-\mathcal{I}_{s,s}(x)&\le Y_{s,s}-Y_{r_-,s},\\
0\le \mathcal{I}_{s,s}(x)-\mathcal{I}_{r_+,s}(x)&\le Y_{r_+,s}-Y_{s,s}.
\end{align*}
Use common bulk weights. Independently of these and the initial walks,
couple each diagonal weight $\mathrm{Geom}(r\sqrt q)$ increasingly in $r$ by
a common uniform variable at that site. Thus every model retains the
independence and marginal distributions required by \eqref{recurrence}--\eqref{weights}.

The recurrence \eqref{recurrence} is monotone in every weight and in the initial
profile. Also, if two profiles differ everywhere by at most a constant
and their weights agree, their last-passage times differ by at most
that constant: this follows by induction using the $1$-Lipschitz property
of the maximum. For the lower inequality, first lower the diagonal
weights of the middle model to those of $r_-$, then compare the profiles.
For the upper inequality, first raise them to those of $r_+$.
This proves \eqref{boundary:eq:pathwise}.

Since $\E[\mathrm{Geom}(p)]=p/(1-p)$, monotonicity of the coupling gives exactly
\begin{equation}\label{boundary:eq:gapmean}
\frac{\E|Y_{r_\pm,s}-Y_{s,s}|}{\rho}
=\frac{s\varepsilon}{\rho^2(1-r_\pm s)(1-s^2)}.
\end{equation}
For $\rho\ge3|\s|$, direct expansion of \eqref{boundary:eq:parameters} gives
\[
\rho(1-s^2)\ge|\s|,\qquad
\rho(1-r_\pm s)\ge|\s|.
\]
Also $0<r_-<s<r_+<1$. Increase the lower bound on $N$, if needed,
so that $s>\sqrt q$. Then all geometric parameters and all stationary
models above are valid, uniformly in the indicated range of $\varepsilon$.
Consequently \eqref{boundary:eq:gapmean} is at most
$\varepsilon/\s^2$. Markov's inequality at the level
$\rho\sqrt\varepsilon$ proves \eqref{boundary:eq:couplingerror}, with
$C=|\s|^{-2}$. No independence between the two coupled $Y$ variables
is required.
\end{proof}

At $r=s$, Definition~\ref{def: stationary measure} identifies the initial condition with the product
stationary measure in Theorem~\ref{thm: low-density phase theorem}. By Proposition~\ref{prop: stationary measure} and the recurrence \eqref{recurrence},
\begin{equation}\label{boundary:eq:endpoint-shift}
\E\bigl[\mathbf{G}_{s,s}(N+2,N+2)-\mathbf{G}_{s,s}(N,N)\bigr]
=2\left(\frac{\sqrt q/s}{1-\sqrt q/s}
+\frac{s\sqrt q}{1-s\sqrt q}\right)=O(1).
\end{equation}
Indeed, each diagonal step consists of a horizontal increment with law
$\mathrm{Geom}(\sqrt q/s)$ and a diagonal weight with law $\mathrm{Geom}(s\sqrt q)$.
The bound is uniform for the parameters in \eqref{boundary:eq:parameters} when $N$
is sufficiently large. Since the difference is nonnegative, Markov's
inequality gives
\[
\frac{\mathbf{G}_{s,s}(N+2,N+2)-\mathbf{G}_{s,s}(N,N)}{\rho}
\longrightarrow0
\qquad\text{in probability}.
\]
Thus Theorem~\ref{thm: low-density phase theorem} also applies at the endpoint used in \eqref{boundary:eq:pathwise},
with the same centering and scale.

\begin{proof}[Proof of the second limit in \eqref{boundary:eq:claim}]
Fix $\varepsilon>0$. Subtract the common
centering $\frac{2\sqrt q}{1-\sqrt q}N$ in \eqref{boundary:eq:pathwise} and
divide by $\rho$. On the event that the corresponding $Y$ difference
is at most $\rho\sqrt\varepsilon$, the lower and upper comparisons
give the following CDF bounds:
\begin{align}
&\Pb\!\left(
\frac{\mathbf{G}_{r_+,s}(N+2,N+2)-\frac{2\sqrt q}{1-\sqrt q}N}{\rho}
\le\tilde{d}-\sqrt\varepsilon\right)-C\sqrt\varepsilon\nonumber\\
&\quad\le
\Pb\!\left(
\frac{\mathbf{G}_{s,s}(N+2,N+2)-\frac{2\sqrt q}{1-\sqrt q}N}{\rho}
\le\tilde{d}\right)\label{boundary:eq:CDFsqueeze}\\
&\quad\le
\Pb\!\left(
\frac{\mathbf{G}_{r_-,s}(N+2,N+2)-\frac{2\sqrt q}{1-\sqrt q}N}{\rho}
\le\tilde{d}+\sqrt\varepsilon\right)+C\sqrt\varepsilon.\nonumber
\end{align}
First let $N\to\infty$. The two outer probabilities converge by
Theorem~\ref{thm: high-density phase theorem}: their fixed scaled high-density parameters are
$\rr=\s+\varepsilon$ and $\rr=\s-\varepsilon$, respectively.
Both are in the domain of Theorem~\ref{thm: high-density phase theorem}, since
$0<\varepsilon<|\s|/4$ implies
$\s-\varepsilon<\s<\s+\varepsilon<-\s$.
The middle probability converges to the low-density limit by Theorem~\ref{thm: low-density phase theorem} and
\eqref{boundary:eq:endpoint-shift}. Hence
\begin{align}
&P^{\mathrm{HD}}_{\s+\varepsilon,\s}(\tilde{d}-\sqrt\varepsilon)
-C\sqrt\varepsilon\nonumber \le
F^{\mathrm{LD}}_{\s}(\tilde{d})
\label{boundary:eq:limitsqueeze}\le
P^{\mathrm{HD}}_{\s-\varepsilon,\s}(\tilde{d}+\sqrt\varepsilon)
+C\sqrt\varepsilon.\nonumber
\end{align}
Now let $\varepsilon\downarrow0$. Lemma~\ref{boundary:lem:average} gives
joint continuity in the parameter and threshold, so both outer
expressions tend to $P^{\mathrm{HD}}_{\s,\s}(\tilde{d})$ in \eqref{boundary:eq:average}. The errors vanish. The squeeze
therefore proves the identity
\begin{equation}\label{boundary:eq:finalidentity}
(-2\s)\int_0^\infty e^{2\s\tilde{y}}
\left.\Pf(J-\widetilde K)_{\mathbb{L}^2(\widetilde T)}\right|_{\rr=\s}\,\mathrm{d}\tilde{y}
=F^{\mathrm{LD}}_{\s}(\tilde{d}).
\end{equation}
The analytic extension in Lemma~\ref{boundary:lem:average} shows that this is
the ordinary two-sided limit as $\rr\to\s$, completing the proof.
\end{proof}

\appendix

\section*{Appendix}

\begin{proof}[Proof of Lemma~\ref{lem:geometric-exp-kernel}]
Recall the kernel defined in Lemma \ref{lem: one param kernel reformulated}. To distinguish the geometric
parameters from their exponential counterparts, we set
\[
q_\varepsilon=e^{-2\varepsilon\alpha},
\qquad
t_\varepsilon=e^{-\varepsilon t},
\qquad
r_\varepsilon=e^{-\varepsilon r}.
\]
Thus $\sqrt{q_\varepsilon}<t_\varepsilon,r_\varepsilon<1$, as required
by Lemma~\ref{lem: one param kernel reformulated}. It is clear that under this scaling, we have
$\varepsilon\,\mathrm{Geom}(e^{-\varepsilon a})
\Rightarrow\mathrm{Exp}(a)$.

Write $\mathcal K^{(\varepsilon)}$ for the geometric kernel with these
parameters. On the lattice $X=\varepsilon k$, $Y=\varepsilon\ell$, define
\[
\widehat{\mathcal K}_\varepsilon(X,Y)
:=
D_\varepsilon\mathcal K^{(\varepsilon)}(k,\ell)D_\varepsilon,
\qquad
D_\varepsilon=\operatorname{diag}(\varepsilon^{-2},\varepsilon).
\]
Explicitly, the four entries are rescaled by
$\varepsilon^{-4},\varepsilon^{-1},\varepsilon^{-1},\varepsilon^2$,
respectively.

Choose a positively oriented simple closed contour
$\Gamma$ contained in the open left half-plane and enclosing
$-\alpha,-t,-r$. Fix $\eta>0$ such that
$\operatorname{Re}Z\leq-\eta$ on $\Gamma$.
For sufficiently small $\varepsilon$, the maps
$Z\mapsto e^{-\varepsilon Z}$ and $W\mapsto e^{\varepsilon W}$
are injective on a neighborhood of $\Gamma$ and its interior.
Set
\[
\gamma_\varepsilon^z=e^{-\varepsilon\Gamma},
\qquad
\gamma_\varepsilon^w=e^{\varepsilon\Gamma},
\]
with their inherited positive orientations. These contours enclose,
respectively,
\[
\{q_\varepsilon^{-1/2},t_\varepsilon^{-1},r_\varepsilon^{-1}\},
\qquad
\{q_\varepsilon^{1/2},t_\varepsilon,r_\varepsilon\}.
\]
Moreover, $|z|>1>|w|$ on their product, so the pole $z=w$ is never
crossed.

For the $11$-entry, enlarging both contours to these contours adds no
poles; the same is true for the $z$-contour in the $12$-entry.
For the $22$-entry, the explicit double-residue term in
\eqref{eq: low-density kernel} is precisely the contribution of
\[
(z,w)=(t_\varepsilon^{-1},r_\varepsilon),
\qquad
(z,w)=(r_\varepsilon^{-1},t_\varepsilon).
\]
Consequently, the entire $22$-entry can be written as
\[
\mathcal K^{(\varepsilon)}_{22}(k,\ell)
=
-\oint_{\gamma_\varepsilon^w}\frac{dw}{2\pi\mathrm i}
 \oint_{\gamma_\varepsilon^z}\frac{dz}{2\pi\mathrm i}\,
 \mathcal U_\varepsilon(z,w)
+E_\varepsilon(k,\ell),
\]
where $\mathcal U_\varepsilon$ is $\mathcal U$ evaluated at the scaled
parameters, and
\[
E_\varepsilon(k,\ell)
=
\operatorname{sgn}(\ell-k)
\frac{t_\varepsilon^{|\ell-k|}
      -r_\varepsilon^{|\ell-k|}}
     {(1-t_\varepsilon r_\varepsilon)
      (t_\varepsilon-r_\varepsilon)}.
\]

In every double integral, make the change of variables
\[
z=e^{-\varepsilon Z},
\qquad
w=e^{\varepsilon W}.
\]
Both transformed contours are exactly $\Gamma$. Writing
\[
H_\varepsilon(z)
=
\left(
\frac{1-e^{-\varepsilon\alpha}/z}
     {1-e^{-\varepsilon\alpha}z}
\right)^{N-1},
\]
we have
\[
\frac{H_\varepsilon(e^{-\varepsilon Z})}
     {H_\varepsilon(e^{\varepsilon W})}
=
\phi_\varepsilon(Z)\phi_\varepsilon(W),
\qquad
\phi_\varepsilon(Z)
=
\left(
\frac{1-e^{-\varepsilon(\alpha-Z)}}
     {1-e^{-\varepsilon(\alpha+Z)}}
\right)^{N-1}
\longrightarrow
\left(\frac{\alpha-Z}{\alpha+Z}\right)^{N-1},
\]
uniformly on $\Gamma$. Also,
\[
w^\ell z^{-k}=e^{XZ+YW},
\qquad
dz\,dw=-\varepsilon^2e^{\varepsilon(W-Z)}\,dZ\,dW.
\]
The remaining factors satisfy, uniformly on $\Gamma\times\Gamma$,
\[
\begin{aligned}
\varepsilon^{-1}(z-t_\varepsilon)&\longrightarrow t-Z,
&
\varepsilon^{-1}(1-t_\varepsilon w)&\longrightarrow t-W,
\\
\varepsilon^{-1}(w-t_\varepsilon)&\longrightarrow W+t,
&
\varepsilon^{-1}(1-t_\varepsilon z)&\longrightarrow Z+t,
\\
\varepsilon^{-1}(z-w)&\longrightarrow-(Z+W),
&
\varepsilon^{-1}(zw-1)&\longrightarrow W-Z,
\\
\varepsilon^{-1}(z^2-1)&\longrightarrow-2Z,
&
\varepsilon^{-1}(1-w^2)&\longrightarrow-2W,
\end{aligned}
\]
and similarly with $t_\epsilon$ replaced by $r_\epsilon$.
Since $\Gamma$ avoids all poles and
$\operatorname{Re}(Z+W)\leq-2\eta$, every normalized denominator is
uniformly bounded away from zero. Substitution therefore gives the
double-integral parts of
\eqref{eq:Kshift11}, \eqref{eq:Kshift12}, and
\eqref{eq:Kshift22}, with precisely the rescalings specified above.

For the remaining term, use
\[
(1-t_\varepsilon r_\varepsilon)
(t_\varepsilon-r_\varepsilon)
=
\varepsilon^2(r^2-t^2)+O(\varepsilon^3).
\]
It follows that
\[
\varepsilon^2E_\varepsilon(k,\ell)
\longrightarrow
\operatorname{sgn}(Y-X)
\frac{e^{-t|Y-X|}-e^{-r|Y-X|}}{r^2-t^2}.
\]
This expression vanishes continuously at $X=Y$.
Thus $\widehat{\mathcal K}_\varepsilon$ converges to
$K^{\mathrm{gexp}}$, locally uniformly on the rescaled lattice.
The $21$-entry follows from skew-symmetry.

Then we prove the convergence of the Fredholm Pfaffian.
Fix $d\geq0$ and let
\[
d_\varepsilon=\lfloor d/\varepsilon\rfloor,
\qquad
T_\varepsilon=\{d_\varepsilon+1,d_\varepsilon+2,\ldots\}.
\]
The fixed-contour estimates above give, uniformly for small
$\varepsilon$ and $X,Y\geq0$, a bound
$Ce^{-\eta(X+Y)}$ for every rescaled double-integral entry.
Furthermore,
$|\varepsilon^2E_\varepsilon(k,\ell)|\leq C$.
Hence, with
\[
M(X)=\operatorname{diag}(e^{\eta X/2},e^{-\eta X/2}),
\qquad
\mathcal L_\varepsilon(X,Y)
=
M(X)\widehat{\mathcal K}_\varepsilon(X,Y)M(Y),
\]
all entries satisfy
\[
|(\mathcal L_\varepsilon)_{ij}(X,Y)|
\leq Ce^{-\eta(X+Y)/2}.
\]
Hadamard's inequality gives
\[
\left|
\operatorname{Pf}
\bigl[\widehat{\mathcal K}_\varepsilon(X_a,X_b)\bigr]_{a,b=1}^{m}
\right|
\leq
C^m(2m)^{m/2}\prod_{a=1}^{m}e^{-\eta X_a}.
\]

Since $\det D_\varepsilon=\varepsilon^{-1}$, the discrete Fredholm
expansion is exactly
\[
\operatorname{Pf}(J-\mathcal K^{(\varepsilon)})_{\ell^2(T_\varepsilon)}
=
\sum_{m=0}^{\infty}\frac{(-1)^m}{m!}\,
\varepsilon^m
\sum_{k_1,\ldots,k_m\in T_\varepsilon}
\operatorname{Pf}
\bigl[
\widehat{\mathcal K}_\varepsilon
(\varepsilon k_a,\varepsilon k_b)
\bigr]_{a,b=1}^{m}.
\]
For each $m$, local uniform convergence and the exponential bound
justify convergence of this Riemann sum to the corresponding
$m$-fold integral. Moreover, its absolute value, including the factor
$1/m!$, is bounded by
\[
\frac{C_d^m(2m)^{m/2}}{m!},
\]
uniformly in $\varepsilon$, and this bound is summable in $m$.
Dominated convergence therefore yields
\[
\lim_{\varepsilon\downarrow0}
\operatorname{Pf}(J-\mathcal K^{(\varepsilon)})_{\ell^2(T_\varepsilon)}
=
\operatorname{Pf}(J-K^{\mathrm{gexp}})_{L^2(d,\infty)}.
\]
\end{proof}

{
\begin{proof}[Proof of Lemma~\ref{lem:loggamma-exp-kernel}]
We record the degeneration because it also fixes the normalization of the exponential parameters used in this section.  Under the above specialization, the scalar part of the log-gamma kernel in \cite[Theorem~5.10]{imamura2022solvable} (equivalently, the specialization written in \cite[Theorem~5.1]{zeng2026}) contains the factors
\[
\frac{\Gamma(\varepsilon r+Z)}{\Gamma(\varepsilon r-Z)},\qquad
\frac{\Gamma(\varepsilon t+Z)}{\Gamma(\varepsilon t-Z)},\qquad
\left(\frac{\Gamma(\varepsilon\alpha+Z)}{\Gamma(\varepsilon\alpha-Z)}\right)^{N-1},
\]
and the same factors with $Z$ replaced by $W$, together with
\[
Q(Z,W)=\Gamma(-2Z)\Gamma(-2W)
\frac{\sin\!\bigl(\pi(Z-W)\bigr)}{\sin\!\bigl(\pi(Z+W)\bigr)}.
\]
Set $Z=\varepsilon z$ and $W=\varepsilon w$.  Since
\[
\Gamma(\varepsilon c)=\frac1{\varepsilon c}+O(1),
\qquad \varepsilon\downarrow0,
\]
uniformly when $c$ ranges over compact sets bounded away from zero, we obtain
\[
\frac{\Gamma(\varepsilon(a+z))}{\Gamma(\varepsilon(a-z))}
\longrightarrow \frac{a-z}{a+z},
\qquad a\in\{r,t,\alpha\}.
\]
Moreover,
\begin{align*}
Q(\varepsilon z,\varepsilon w)
&=\Gamma(-2\varepsilon z)\Gamma(-2\varepsilon w)
\frac{\sin\!\bigl(\pi\varepsilon(z-w)\bigr)}
     {\sin\!\bigl(\pi\varepsilon(z+w)\bigr)}\\
&=\frac1{4\varepsilon^2zw}\frac{z-w}{z+w}\,(1+O(\varepsilon)),
\end{align*}
uniformly on compact contours avoiding $0$ and $z+w=0$.  Finally,
\[
e^{XZ+YW}=e^{xz+yw},\qquad dZ\,dW=\varepsilon^2\,dz\,dw.
\]
Hence the scalar double-contour kernel converges to
\[
\oint_\Gamma\frac{dz}{2\pi\mathrm i}
\oint_\Gamma\frac{dw}{2\pi\mathrm i}\,
 e^{xz+yw}
 \frac{r-z}{r+z}\frac{t-z}{t+z}
 \left(\frac{\alpha-z}{\alpha+z}\right)^{N-1}
 \frac{r-w}{r+w}\frac{t-w}{t+w}
 \left(\frac{\alpha-w}{\alpha+w}\right)^{N-1}
 \frac1{4zw}\frac{z-w}{z+w}.
\]
This is exactly $K^{\mathrm{lexp}}_{11}(x,y)$.

The matrix kernel in \cite[Theorem~5.10]{imamura2022solvable} is obtained from the scalar kernel by applying $-\partial_Y$, $-\partial_X$, and $\partial_X\partial_Y$ in the $12$, $21$, and $22$ entries, respectively.  Under $X=x/\varepsilon$ and $Y=y/\varepsilon$, these derivatives contribute respectively the factors $-\varepsilon w$, $-\varepsilon z$, and $\varepsilon^2zw$.  At the same time the change of variables in the Fredholm operator contributes the Jacobian $\varepsilon^{-1}$, and conjugation by
$D_\varepsilon=\operatorname{diag}(\varepsilon^{1/2},\varepsilon^{-1/2})$
normalizes the four entries.  The limiting factors $w,z,zw$ cancel the corresponding factors in $1/(zw)$ and give precisely \eqref{eq:Kexp12}--\eqref{eq:Kexp22}; skew-symmetry gives the $21$ entry.

The contour may be chosen, after the rescaling, as a fixed compact contour $\Gamma$ in the open left half-plane enclosing $-r,-t,-\alpha$ and excluding $0$.  On such a contour all the above expansions are uniform.  Therefore the kernel convergence is locally uniform in $x,y$.  The same fixed-contour estimates give exponential decay in $x+y$ on $\mathbb L^2(d,\infty)$, so the Fredholm series are dominated term by term and the corresponding Fredholm Pfaffians converge as well.
\end{proof}
}

\bibliographystyle{amsalpha}
\bibliography{main}

\end{document}